\documentclass[draftclsnofoot,onecolumn,journal,12pt]{IEEEtran}

\usepackage{amssymb,graphicx}
\usepackage{latexsym}
\usepackage{mathtools}
\usepackage{cite}
\usepackage{color}
\usepackage{subfigure}
\usepackage{bm,bbm}
\usepackage{tabularx,booktabs}
\usepackage{ltablex}
\usepackage{mathrsfs}
\usepackage{enumitem}
\usepackage{cases}
\usepackage{dsfont}

\newcolumntype{M}[1]{>{\centering\arraybackslash}m{#1}}

\usepackage[linesnumbered,lined,boxed]{algorithm2e}
\usepackage{setspace}
\SetAlCapSkip{1em}
\SetKwComment{Comment}{$\triangleright$\ }{}
\RequirePackage{amsthm,amsmath}
\RequirePackage[colorlinks,citecolor=blue,urlcolor=blue]{hyperref}
\newcommand*{\fullref}[2]{\hyperref[{#1}]{\ref*{#1} {\it(#2)}}}
\newcommand*{\partialref}[2]{\hyperref[{#1}]{#2}}

\usepackage[deletedmarkup=xout]{changes}

\usepackage{mathrsfs,BOONDOX-cal}
\graphicspath{{./figs/}}

\newcommand{\plainproof}{\emph{Proof.}~}
\renewcommand{\proof}[1]{\emph{Proof of #1.}~}
\renewcommand{\endproof}{{\hfill $\blacksquare$}}

\newtheorem{definition}{\bf Definition}
\newtheorem{theorem}{\bf Theorem}
\newtheorem{lemma}{\bf Lemma}
\newtheorem{proposition}[theorem]{\bf Proposition}
\newtheorem{corollary}{\bf Corollary}

\newtheorem{condition}{\bf Condition}

\makeatletter
\def\env@cases{%
  \let\@ifnextchar\new@ifnextchar
  \left.
  \def\arraystretch{1.2}%
  \array{@{}l@{\,\,}l@{}}%
}%
\makeatother

\begin{document}
\title{A Restless Bandit Model for Energy-Efficient Job Assignments in Server Farms}
\author{Jing Fu, \IEEEmembership{Member, IEEE}, Xinyu Wang, Zengfu Wang, and Moshe Zukerman, \IEEEmembership{Life Fellow, IEEE}
\thanks{This work was supported in part by the Shenzhen Municipal Science and Technology Innovation Committee under Project JCYJ20180306171144091, a grant from the Research Grants Council of the Hong Kong SAR, P. R. China (CityU 11200318), and Dr Jing Fu's start-up funds from RMIT (Ref. 800132109).}
\thanks{Jing Fu is with the School of Engineering, STEM College, RMIT University, VIC3000, Australia (e-mail: jing.fu@rmit.edu.au). }
\thanks{Zengfu Wang is with the Research \& Development Institute of Northwestern Polytechnical University in Shenzhen, Shenzhen 518057, P. R. China, and with the School of Automation, Northwestern Polytechnical University, Xi'an 710072, P. R. China.(e-mail: wangzengfu@nwpu.edu.cn).}
\thanks{Xinyu Wang and Moshe Zukerman are with the Department of Electrical Engineering, City University of Hong Kong, Hong Kong SAR, P. R. China (e-mail: xywang47-c@my.cityu.edu.hk; moshezu@cityu.edu.hk).}}

\maketitle

\begin{abstract}
We aim to maximize the energy efficiency, gauged as average energy cost per job, in a large-scale server farm with various storage or/and computing components modeled as parallel abstracted servers. 
Each server operates in multiple power modes characterized by potentially different service and energy consumption rates.
The heterogeneity of servers and multiple power modes complicate the maximization problem, where optimal solutions are generally intractable. 
Relying on the Whittle relaxation technique, we resort to a near-optimal, scalable job-assignment policy. 
Under a mild condition related to the service and energy consumption rates of the servers, we prove that our proposed policy approaches optimality as the size of the entire system tends to infinity; that is, it is asymptotically optimal. 
For the non-asymptotic regime, we show the effectiveness of the proposed policy through  numerical simulations, where the policy  outperforms all the tested baselines, and we numerically demonstrate its robustness against heavy-tailed job-size distributions. 
\end{abstract}

\begin{IEEEkeywords}
Restless bandit; job-assignment; asymptotic optimality.
\end{IEEEkeywords}

\section{Introduction}
\label{sec:intro}

\IEEEPARstart{T}{he} ever-increasing demand for internet services in recent decades has led to explosive growth in data centers and the markets of computing and storage infrastructures to the so-called Zettabyte Era \cite{cisco2018whitepaper,shehabi2016USreport}. In 2014, 
U. S. data centers were reported to consume around $70$ billion kWh of annual electricity and this consumption is predicted to continue increasing~\cite{shehabi2016USreport}.
Computing and storage components have been considered major contributors to power consumption in data centers~\cite{kliazovich2015energy,dayarathna2016data}.
We study energy-efficient scheduling policies in a large server farm with widely deployed abstracted servers, each of which represents a  physical component used to serve incoming customer requests.

Methodologies applicable to server farm scheduling or network resource allocation have been studied from several perspectives.
Energy-efficient policies were considered in~\cite{lu2013simple,gelenbe2013energy,gebrehiwot2017near} with only identical servers,
and in~\cite{chowdhury2012vineyard,hu2013resolve,esposito2016distributed,feng2017approximation} through static scheduling mechanics without dynamic reuse of released physical resources.
Nonetheless, to meet the various demands of internet customers, 
service providers have launched a large number of networked facilities with highly diverse physical features in backhaul networks, where frequent reuse of released resources is the preferred option because of its efficiency benefits.
The heterogeneity of backhaul network facilities keeps increasing in terms of not only their characteristics of implementing specific functions but also their generations~\cite{guo2014server}.
Advanced virtualization techniques abstract these physical components as network resources in modern Cloud computing platforms \cite{hameed2016survey}. 
It is important to consider the heterogeneity of such abstracted components in the development of scheduling policy because it has substantial implications for the optimized profit of service providers and costs incurred by customers.

Research on the development of dynamic resource allocation methodologies for large-scale server farms (networks) with the reuse of released physical resources has been conducted under certain simplifying assumptions. The work in~\cite{rosberg2014,hyytia2014task,fu2015energy,lin2017game} considered heterogeneous servers but under the assumption of negligible power consumption of idle servers, while in ~\cite{fu2016asymptotic,wei2017data,mishra2018sustainable,fu2020energy} it was assumed that servers either operate at their peak power consumption rates or stay idle. Then, in \cite{akgun2013energy,li2018online,fu2018restless,wang2018energy,fu2020resource}, it was assumed that  servers' power consumption linearly increases in their service rates. Such assumptions of specific power functions  simplify the analysis of the relevant problems.
However, the diversity of computing/storage components prevents a specific function of their power consumption from being always applicable  and there is a need for a methodology applicable to any power consumption behavior. 
Publications focusing on the server, GPU, and storage component power consumption pointed out that real-life power behavior does not support the convexity (or linearity) assumption \cite{lewis2012runtime,bansal2013speed,mei2017survey}. 
In this paper, we do not assume convexity or linearity of power functions but consider more general and practical situations and in this way, our policy solutions can apply to a wider range of practical scenarios. 

We focus on energy-efficient server farms consisting of a large number of abstracted components that are potentially different in terms of service rates, power consumption rates, and their abilities to serve different jobs. 
Optimizing resource allocation in server farms is achieved  in the vein of the \emph{restless multi-armed bandit problem} (RMABP) proposed in \cite{whittle1988restless}.
The RMABP is a special type of Markov decision process (MDP) consisting of parallel \emph{bandit processes}, which are also MDPs evolving with binary actions, referred to as the \emph{active} and \emph{passive} modes. 
The RMABP includes a large number of such bandit processes that are competing for limited opportunities of being evolved in the active mode. 
In \cite{papadimitriou1999complexity}, RMABP was proved to be PSPACE-hard in general.
Whittle \cite{whittle1988restless} proposed the classical \emph{Whittle index policy} through the \emph{Whittle relaxation technique} and conjectured its asymptotic optimality; that is, the Whittle index policy was conjectured to approach optimality as the scale of the entire system tends to infinity.
The Whittle index policy is scalable for large problems and, if it is asymptotically optimal, its performance degradation is bounded and diminishes when the problem size becomes larger and larger.
Nonetheless, in general, Whittle relaxation technique does not ensure either the existence of \emph{Whittle indices}, the main parameters required to construct the Whittle index policy, or the bounded performance degradation.
Whittle indices were originally defined in \cite{whittle1988restless} under a condition, subsequently referred to as \emph{Whittle indexability}.
In \cite{weber1990index}, Weber and Weiss proved the asymptotic optimality of the Whittle index policy with an extra, non-trivial condition that requires the existence of a \emph{global attractor} of a proposed process associated with the RMABP. In the research field of RMABP, the discussions on Whittle indexability and the global attractor remain open questions in the past several decades.

In \cite{fu2020energy}, based on the results in \cite{fu2018restless}, a scalable job-assignment policy was proposed and proved to be asymptotically optimal for a simplified server farm, where 
the servers/components were assumed to have only two power modes (corresponding to two power consumption rates).
We refer to Section~\ref{sec:rWork} for a detailed survey of RMABP and other related work.

The contributions of this paper are listed as follows.  
\begin{itemize}
 
\item
We provide a scalable policy that aims to maximize the energy efficiency of a large-scale system of deployed computing/storage clusters.
This policy always prioritizes physical components (abstracted servers) according to the descending order of their associated \emph{indices}, which are real numbers representing the marginal rewards gained by selecting these components. 
The indices are pre-computed and the complexity of implementing the index policy is only linear in the number of available physical components.
We refer to the policy as \emph{Multiple Power Modes with Priorities} (MPMP), reflecting its applicability to a parallel-server system with multiple power modes.

\item \label{cont:asym_opt}
When job sizes are exponentially distributed, under a mild condition related to the service and energy consumption rates of physical components, we prove that the index policy approaches optimality as the job arrival rates and the number of physical components in each cluster tend to be arbitrarily large proportionately; that is, it is asymptotically optimal.
The asymptotic optimality is appropriate for computing/storage clusters with a rapidly increasing number of physical components.
We prove that the performance deviation between our proposed MPMP policy and the optimal point in the asymptotic regime diminishes exponentially in the size of the problem. It implies that the MPMP policy is already near-optimal for a relatively small system.

Recall that no previous results can be directly applied for scalable, asymptotically optimal policies in the large server farm, where the abstracted computing/storage components operate in multiple power/service modes.
The complexity of the server farm problem requires a new analysis of the entire system, provided in this paper, including discussions on the indexability and the global attractor for proving asymptotic optimality in the continuous-time case.

\item  We numerically demonstrate the effectiveness of MPMP in the general case, 
where MPMP significantly outperforms baseline policies in all the tested cases. 
We further explore its performance with different job-size distributions and numerically show that the energy efficiency of MPMP is not very sensitive to different shapes of job-size distributions.
\end{itemize}

The remainder of the paper is organized as follows.
In Section~\ref{sec:rWork}, we discuss other related work for job assignments and RMABP.
In Section~\ref{sec:model}, a description of the server farm
model is provided, and in Section~\ref{sec:problem},  the underlying stochastic
optimization problem is rigorously defined. 
In Section~\ref{sec:indexability}, we discuss the indexability of the underlying stochastic process and propose the \emph{indices} - the most important parameters for constructing our policy.
In Section~\ref{sec:algo}, we formally define the MPMP policy, and in Section~\ref{sec:asym_opt} we prove its asymptotic optimality.
Section~\ref{sec:simulation} presents extensive numerical results that demonstrate the effectiveness of MPMP in the general case.
The conclusions of this paper are included in Section~\ref{sec:conclusion}.

\section{Other Related Work}\label{sec:rWork}

Job-assignment policies with strict capacity constraints of physical resources have been studied in~\cite{fu2016asymptotic,wang2018energy,fu2020energy}, where  the release and reuse of physical resources were considered.
Following the ideas of restless bandits \cite{whittle1988restless,weber1990index}, the authors of \cite{fu2016asymptotic,wang2018energy,fu2020energy} proposed scalable policies and proved that the policies, which do not necessarily perform well in small systems, approach the optimal solution in large-scale systems. Optimization problems for small systems can be solved by conventional algorithms, which cannot be directly applied in large cases because of high computational complexity.
Nevertheless, these publications assumed either two power modes (that is, fixed power consumption values for busy and idle servers) or power consumption linearly increasing in servers' traffic loads.
As mentioned in Section~\ref{sec:intro}, here, we overcome the weaknesses of past publications and provide general solutions that are applicable to realistic  power functions.

Apart from the job-assignment problems, conventional RMABP has been widely studied and applied to scheduling problems.
For instance, 
in \cite{yu2018deadline}, a set of identical servers/processors were scheduled to serve stochastically identical jobs that keep arriving.
In \cite{borkar2017whittle}, Borkar considered a special type of bandit processes, of which the state variables take binary values and are only partially observable. He proved the \emph{Whittle indexability}, a key property for an RMABP, by generating and analyzing the indexability of an equivalent process of the original bandit process. 
In \cite{wang2019whittle}, Whittle indexability was proved for a channel-selecting problem where each bandit process was associated with a wireless channel and its state variable was defined as the number of successive transmission failures in that channel.
In \cite{abbou2019group}, a group maintenance problem was modeled as a standard RMABP with a detailed analysis of its Whittle indexability.

In \cite{ninomora2001restless}, Ni\~no-Mora proposed the partial conservation law (PCL) and the PCL-indexability for RMABP. The latter was proved to imply (be stronger than) the Whittle indexability. Later in 2002, Ni\~no-Mora \cite{ninomora2002dynamic} identified a set of optimization problems that satisfy PCL-indexability and thus are Whittle indexable. A detailed survey about the Whittle and PCL-indexability was provided in \cite{ninomora2007dynamic}.
In \cite{ninomora2020verification}, Ni\~no-Mora defined the indexability and Whittle indices of a bandit process with continuous state space and proposed a method that verified the indexability and computed the corresponding Whittle indices. 
All in all, these studies have established computational methodologies for verifying Whittle indexability for the general RMABP that aims to maximize/minimize the expected cumulative rewards/costs. 
Our work, in this paper, aims to maximize a long-run average objective that prevents the existing off-the-shelf techniques from being applied directly. Although from \cite{ross1992applied}, optimizing the long-run average rewards/costs of an MDP can usually be translated to a problem that optimizes the expected cumulative rewards/costs of the same process with an attached, real-valued criterion, this real-valued criterion is not known a priori and has a strong impact on the indexability of the underlying bandit processes. We refer to Section~\ref{subsec:indexability} for a detailed explanation for demonstrating the indexability with the long-run average objective for our server farm.

For a general RMABP, to further prove that the Whittle index policy approaches optimality as the number of bandit processes increases to infinity (that is, asymptotic optimality), Weber and Weiss \cite{weber1990index} required another non-trivial condition; namely, there exists a fixed point such that the underlying stochastic process of the RMABP will almost surely enter a nearby neighborhood of the point. Such a point is referred to as a \emph{global attractor} of the process. 
In \cite{ouyang2015downlink,maatouk2020asymptotically}, similar assumptions related to the global attractor were required in the proofs of asymptotic optimality of Whittle index policy in channel selection problems, which are special cases of RMABP.
In \cite{fu2018restless}, for a system consisting of a special type of bandit processes, Fu \emph{et al.} proved that such a global attractor exists and hence the resulting policy is asymptotically optimal. The idea of this technique was later applied to a server farm model in \cite{fu2020energy}.
Nonetheless, due to the complexity of the server farm model considered in this paper, it does not fall exactly in the scope of \cite{fu2018restless} and we cannot directly apply the same conclusions here. 
Recall that modeling the server farm as an RMABP or RMABP-like system cannot ensure the existence of a scalable near-optimal policy with theoretically bounded performance.
The complexity of the problem requires a new thorough analysis of the indexability and global attractor, which have been, in general, open questions in the past several decades.
In this paper, we resort to scalable policies for such a challenging server farm model with a theoretical performance guarantee.

Our job-assignment problem can be modeled as a set of parallel, restless bandit processes coupled with action constraints. A detailed description of our model is provided in Section~\ref{sec:model}.
As discussed, in this paper, we generalize the power functions discussed in \cite{fu2020energy}, which requires a new analysis. 
In particular, we prove that, for a mild condition related to the service and energy consumption rates of physical components (abstracted servers), the MPMP policy is asymptotically optimal with respect to energy efficiency.
To the best of our knowledge, there is no published work applicable to a server farm model with strict capacity constraints and generalized power functions, where a scalable policy is proposed with theoretically guaranteed performance when the system is largely scaled.

\section{Model}\label{sec:model}
For any positive integer $N$, let $[N]$ represent the set $\{1,2,\ldots,N\}$. Let $\mathbb{R}$, $\mathbb{R}_+$ and $\mathbb{R}_0$ represent the sets of all real numbers, positive real numbers, and non-negative real numbers, respectively.
Similarly, $\mathbb{N}$, $\mathbb{N}_+$ and $\mathbb{N}_0$ are the sets of integers, positive integers, and non-negative integers, respectively.

There are $I$ clusters of physical components. These components are identical within the same cluster in terms of their availability of accommodating jobs, hardware features and software profiles. Physical components in different clusters may be different.

Consider $L$ classes of jobs, each of which is characterized by a tuple $(\lambda_{\ell},\mathcal{I}_{\ell})$ for $\ell\in[L]$, where
\begin{itemize}
\item the arrival process of jobs in class $\ell$ follows a Poisson process with rate $\lambda_{\ell}>0$;
\item the set of clusters that are potentially able to serve jobs of class $\ell$ is given by $\mathcal{I}_{\ell}\in 2^{[I]}\backslash\{\emptyset\}$.
\end{itemize}
The components in the clusters $i\in\mathcal{I}_{\ell}$ are referred to as the \emph{available components} for jobs of class $\ell$.
The sizes of all jobs are considered as independently and identically distributed (i.i.d.) random variables with unit mean, and jobs are arriving sequentially with positive inter-arrival time.
We refer to a job of class $\ell\in[L]$ as an $\ell$-job.

Each physical component of cluster $i\in[I]$ has a finite capacity, $C_{i}\in\mathbb{N}_{+}$, as the maximal number of jobs it can serve simultaneously. Define its energy consumption and service rates as functions of carried load: $\varepsilon_i(n)$ and $\mu_i(n)$, $n\in\{0\}\cup[C_i]$. Note that the service units of each physical component are being reused and released dynamically along the timeline.
For the purpose of this paper, assume that $\varepsilon_i(n)$ and $\mu_i(n)$ are finite, non-negative and increasing in $n$ (that is, $0\leq \varepsilon_i(n)\leq \varepsilon_i(n+1)< \infty$ and $0\leq\mu_i(n)\leq \mu_i(n+1)<\infty$) with $\mu_i(0)\equiv 0$, $\varepsilon_i(0)\geq 0$ and, for $n>0$, $\varepsilon_i(n)>\varepsilon_i(0)$ and $\mu_i(n)>0$.

Consider $M_i^0\cdot h$ components in cluster $i$ where $M_i^0$ and $h$ are positive integers. Let $\lambda_{\ell}=\lambda_{\ell}^{0}\cdot h$ for some $\lambda_{\ell}^0\in\mathbb{R}_+$. In this context, $h$ is a constant that signifies the scale of the multi-cluster system, which is referred to as the \emph{scaling parameter}. There are in total $J=h\sum_{i\in[I]}M_i^0$ components in our system, labeled by $j\in[J]$, with each computing cluster $i\in[I]$ consisting of physical components $j_1,j_2,\ldots,j_{hM_i^0}\in[J]$.
For $j\in[J]$, let $i_j$ represent the label of the cluster with $j\in\mathscr{J}_{i_j}$.

\section{Optimization Problem}\label{sec:problem}

At time $t\geq 0$, there are $N_j(t)$ jobs accommodated by physical component $j$, and define $\bm{N}(t) = (N_j(t):j\in[J])$. The $N_j(t)$ is a random  variable and is referred to as the \emph{state variable} of component $j$.
When a job is assigned to a component $j$ at time $t$, the state of component $j$ transitions from $N_j(t)$ to $N_j(t)+1$; when a job on component $j$ is completed, the state decreases to $N_j(t)-1$. Preemption of jobs is not permitted in our system.
The variable $N_j(t)$ is affected by the underlying scheduling policy, denoted by $\phi$. When an $\ell$-job arrives, the scheduling policy selects a component in cluster $i\in\mathcal{I}_{\ell}$ with at least a vacant slot to accommodate this job or block it. 
For fairness, we do not allow the rejection of any job when there is a vacant slot in any of the available components for this job. In other words, rejection of a job occurs if and only if all the available components for the job are fully occupied. 

To indicate the dependency of the counting process $\{\bm{N}(t),t\geq 0\}$ and the scheduling policy $\phi$, we rewrite the state variables as $N^{\phi}_j(t)$ and the state vector $\bm{N}^{\phi}(t)$.

In this context, the set of all possible values of $N^{\phi}_j(t)$ of component $j$ in cluster $i$ is $\mathscr{N}_i\coloneqq \{0,1,\ldots,C_i\}$; that is, $\mathscr{N}_i$ is the \emph{state space} of the process $\{N^{\phi}_j(t),t\geq0\}$ associated with a component $j$ in cluster $i$, abbreviated as the state space of cluster $i$. 
The state vector $\bm{N}^{\phi}(t)$ of the entire system is then taking values in 
\begin{equation}
\mathscr{N}\coloneqq \prod_{i\in[I]}\left(\mathscr{N}_i\right)^{h\cdot M_i^0},
\end{equation}
of which the size is increasing exponentially in $h$.

More precisely, define \emph{action variable} $a^{\phi}_{\ell,j}(\bm{n})\in\{0,1\}$ as a function of a state vector $\bm{n}\in\mathscr{N}$ associated with policy $\phi$: if $a^{\phi}_{\ell,j}(\bm{n}) = 1$, then a newly-arrived job of class $\ell$ will be accommodated by physical component $j$ when $\bm{N}^{\phi}(t)=\bm{n}$ under policy $\phi$; otherwise, the job will not be assigned to component $j$. To take account for the rejection of jobs, define a  \emph{virtual component} for each class $\ell$, labeled by $j_{\ell}\coloneqq J+\ell$; that is, if $a^{\phi}_{\ell,j_{\ell}}(\bm{n}) = 1$, a new job of class $\ell$ will be blocked. 
For $\ell\in[L]$, assume without loss of generality that the virtual component $j_{\ell}$ belongs to a virtual cluster $i_{\ell}\coloneqq I+\ell$. Define a stochastic process $N^{\phi}_{j_{\ell}}(t) \equiv 0$, $t\geq 0$, associated with the virtual component $j_{\ell}$, which has state space $\mathscr{N}_{i_\ell}\coloneqq \{0\}$. 

To mathematically define the feasibility of policy $\phi$, we introduce the following constraints for the action variables:
\begin{equation}\label{eqn:constraint:action}
\sum\limits_{j\in \mathcal{J}_{\ell}\cup\{j_{\ell}\}}a^{\phi}_{\ell,j}(\bm{n})=1,~\forall \ell\in[L], \bm{n}\in\mathscr{N},
\end{equation}
and
\begin{equation}\label{eqn:constraint:capacity}
a^{\phi}_{\ell,j}(\bm{n}) = 0,\text{if } n_j = C_{i_j},~\forall \ell\in[L],j\in\mathcal{J}_{\ell}, \bm{n}\in\mathscr{N},
\end{equation}
where $\mathcal{J}_{\ell}$ is the set of all available components of job class $\ell$, and
cluster $i_j$ is the cluster in which component $j$ is located.
Define that $a^{\phi}_{\ell,j}(\cdot)\equiv 0$ for all $j\notin \mathcal{J}_{\ell}\cup\{j_{\ell}\}$. 
Constraints~\eqref{eqn:constraint:action} ensure that only one component is selected for an arrived job, and Constraints~\eqref{eqn:constraint:capacity} disable a fully-occupied component from accommodating more jobs.

Recall that a rejection of jobs is not permitted if there is a vacant slot on any available component. This can be addressed by introducing an intermediate variable 
$\bar{a}^{\phi}_{\ell}(\bm{n}) $ ($\bm{n}\in\mathscr{N}$) that satisfies
\begin{equation}\label{eqn:constraint:boundary}
\bar{a}^{\phi}_{\ell}(\bm{n}) + \sum\limits_{j\in\mathcal{J}_{\ell}}I(C_{i_j}-n_j) \leq 1,~
\forall \ell\in[L],\bm{n}\in\mathscr{N},
\end{equation}
where $I(x)$ ($x\in\mathbb{R}$) is a Heaviside function with $I(x)=1$ for $x>0$; and $0$ for $x \leq 0$. For such $\bar{a}^{\phi}_{\ell}(\bm{n})$ that takes values in $(-\infty, 1]$, define 
\begin{equation}
a^{\phi}_{\ell,j_{\ell}}(\bm{n}) \coloneqq I(\bar{a}^{\phi}_{\ell}(\bm{n})).
\end{equation}
From \eqref{eqn:constraint:boundary}, when there is at least one available component in $\mathcal{J}_{\ell}$, $\bar{a}^{\phi}_{\ell}(\bm{n}) \leq 0$ and so $a^{\phi}_{\ell,j_{\ell}}(\bm{n}) =0$. If all these components are fully occupied, constraints in \eqref{eqn:constraint:action} and \eqref{eqn:constraint:capacity} force $a^{\phi}_{\ell,j_{\ell}}(\bm{n}) = 1$, which does not violate constraints in \eqref{eqn:constraint:boundary}.

We aim to maximize the energy efficiency of the entire system; specifically, to maximize the ratio of the long-run average job throughput to the long-run average power consumption. Let
\begin{equation}\label{define:L}
\mathfrak{L}^{\phi}=\lim\limits_{T\rightarrow +\infty} \frac{1}{T}\mathbb{E}\int_{0}^{T} \sum\limits_{j\in[J]}\mu_j(N^{\phi}_j(t)) d t
\end{equation}
and
\begin{equation}\label{define:E}
\mathfrak{E}^{\phi}=\lim\limits_{T\rightarrow +\infty} \frac{1}{T}\mathbb{E}\int_{0}^{T} \sum\limits_{j\in[J]}\varepsilon_j(N^{\phi}_j(t)) d t
\end{equation}
represent the long-run average job throughput and the long-run average power consumption, respectively.
Our optimization problem is 
\begin{equation}\label{eqn:obj}
\max\limits_{\phi}~\frac{\mathfrak{L}^{\phi}}{\mathfrak{E}^{\phi}}
\end{equation}
subject to \eqref{eqn:constraint:action}, \eqref{eqn:constraint:capacity} and \eqref{eqn:constraint:boundary}.
Let $\Phi$ represent the set of all policies constrained by \eqref{eqn:constraint:action}, \eqref{eqn:constraint:capacity} and \eqref{eqn:constraint:boundary}.

As in \cite{fu2016asymptotic,fu2020energy}, 
consider an optimal policy $\phi^*$ that maximizes the problem described in \eqref{eqn:obj}, \eqref{eqn:constraint:action}, \eqref{eqn:constraint:capacity} and \eqref{eqn:constraint:boundary} and define a real number
\begin{equation}\label{eqn:e_star}
e^* = \frac{\mathfrak{L}^{\phi^*}}{\mathfrak{E}^{\phi^*}}.
\end{equation}
Following \cite[Theorem 1]{rosberg2014}, if $0<\mathfrak{L}^{\phi}<+\infty$ and $0<\mathfrak{E}^{\phi}<+\infty$, then a policy $\phi\in\Phi$ that maximizes  
\begin{equation}\label{eqn:opt}
\max\limits_{\phi\in\Phi} \mathfrak{L}^{\phi}-e^*\mathfrak{E}^{\phi},~
\end{equation}
subject to \eqref{eqn:constraint:action},~\eqref{eqn:constraint:capacity} and \eqref{eqn:constraint:boundary}
also maximizes the problem defined in \eqref{eqn:obj}, \eqref{eqn:constraint:action}, \eqref{eqn:constraint:capacity} and \eqref{eqn:constraint:boundary}.

The constraints~\eqref{eqn:constraint:boundary} make our problem slightly different from a standard RMABP in the sense that in our problem, rejecting a job has the lowest priority among all the other actions.
Unlike the \emph{uncontrollable} variables constrained by \eqref{eqn:constraint:capacity}, \eqref{eqn:constraint:boundary} guarantees the lowest priority for rejecting a job rather than disables this action.
If $L=1$ and constraints~\eqref{eqn:constraint:boundary} do not exist, then our problem reduces to an RMABP, where the processes $\{N^{\phi}_j(t), t\geq  0\}$ ($j\in[J]$) are parallel restless bandit processes coupled by action constraint \eqref{eqn:constraint:action}.
Note that such an RMABP is no longer applicable to our server farm.

Similar to an RMABP,  addressing our problem requires overcoming the challenge of its large state space, which is exponentially increasing in the number of components and optimal solutions are intractable. In \cite{fu2020energy}, a scalable policy was proposed for a similar problem with simplified $\mu_i(n)\equiv \mu_i$ and $\varepsilon_i(n)\equiv \varepsilon_i$ for all $n=1,2,\ldots,C_i$. This policy was proved to be asymptotically optimal, for which the performance gap to optimality diminishes exponentially in the scale of the problem.  
However,  in this paper, the generalized $\mu_i(n)$ and $\varepsilon_i(n)$ (that is, the multiple power states) prevent the same technique from being applied directly.
For general RMABPs or some relevant problems, the Whittle relaxation technique does not ensure a good, scalable policy that asymptotically approaches optimality.
From \cite{whittle1988restless, weber1990index},  asymptotic optimality relies on two important but non-trivial properties related to the underlying stochastic process: \emph{Whittle indexability}  and the existence of a \emph{global attractor} in the asymptotic regime. These two properties do not necessarily hold in general and remain open questions in the past several decades.

As mentioned in Section~\ref{sec:intro}, the existence of an appropriate global attractor has been discovered in a class of RMABPs \cite{fu2018restless}, including the special case studied in \cite{fu2020energy}. In this paper, the server farm with more general service and power consumption rates does not fall exactly in the scope of \cite{fu2018restless}. 
This requires a new analysis of the entire system. 
As mentioned in Section~\ref{sec:rWork}, Whittle indexability for relatively general RMABPs has been widely studied with provided sufficient conditions for optimizing cumulative rewards/costs \cite{ninomora2001restless,ninomora2002dynamic,ninomora2007dynamic,ninomora2020verification}. 
Nonetheless, the generalized transition and reward rates of the underlying stochastic process and the objective in this paper - maximization of an average reward - prevent these off-the-shelf techniques from being applied directly.
Although the maximization of the average reward of an MDP can usually be translated to the maximization of the cumulative reward of the same process with an attached real-valued criterion, the exact value of this attached criterion cannot be known a priori and has a strong impact on the discussion of Whittle indexability.   
In \cite{fu2020energy}, for the very special server farm with $\mu_i(n)\equiv \mu_i$ and $\varepsilon_i(n)\equiv \varepsilon_i$, the real-valued criterion was directly offset during the analysis of Whittle indexability, which significantly simplified the entire discussion. Whittle indexability ensures the existence of a scalable policy, derived from \emph{Whittle indices},  for a standard RMABP \cite{whittle1988restless}. 
Unfortunately, such an unknown criterion persists in this paper which requires a thorough analysis of the entire system.
We will discuss in Section~\ref{sec:indexability} the indexability of the server farm with the generic $\mu_i(n)$ and $\varepsilon_i(n)$ and in Section~\ref{sec:asym_opt} further results for asymptotically optimal policies.

\section{Whittle Relaxation and Indexability}
\label{sec:indexability}

\subsection{Whittle Relaxation}\label{subsec:relaxation}

Following the idea of \emph{Whittle relaxation} \cite{whittle1988restless}, we relax constraints~\eqref{eqn:constraint:action} and \eqref{eqn:constraint:boundary} to
\begin{equation}\label{eqn:constraint:relax:action}
\lim\limits_{t\rightarrow+\infty}\mathbb{E}\Bigl[\sum\limits_{j\in \mathcal{J}_{\ell}\cup\{j_{\ell}\}}a^{\phi}_{\ell,j}(\bm{N}^{\phi}(t))\Bigr]=1,~\forall \ell\in[L]
\end{equation}
and
\begin{equation}\label{eqn:constraint:relax:boundary_tmp}
\lim\limits_{t\rightarrow+\infty}\mathbb{E}\Bigl[\bar{a}^{\phi}_{\ell}(\bm{N}^{\phi}(t)) + \sum\limits_{j\in\mathcal{J}_{\ell}}I(C_{i_j}-N^{\phi}_j(t))\Bigr] \leq 1,~\forall \ell\in[L],
\end{equation}
respectively. Similarly, \eqref{eqn:constraint:capacity} can be rewritten as 
\begin{equation}\label{eqn:constraint:relax:capacity}
\lim\limits_{t\rightarrow+\infty}\mathbb{E}\left[a^{\phi}_{\ell,j}(N_j^{\phi}(t))\ |\ N_j^{\phi}(t) = C_{i_j}\right]=0,~\forall \ell\in[L], j\in[J].
\end{equation}

Define a special policy $\phi_0$ with $a^{\phi_0}_{\ell,j}(\bm{n})=1$ for all $\ell\in[L]$, $j\in\mathcal{J}_{\ell}$, and $\bm{n}\in\mathscr{N}$. 
Note that $\phi_0\notin \Phi$ because it violates the constraints on action variables.
We apply the $\phi_0$ to all stochastic processes $\{N^{\phi_0}_j(t),t\geq 0\}$ for $j\in[J]$.
Define
\begin{equation*}
A_{\ell} \coloneqq \sum\limits_{j\in\mathcal{J}_{\ell}}\lim\limits_{t\rightarrow+\infty}\mathbb{E}\left[I(C_{i_j}-N^{\phi_0}_j(t))\right],
\end{equation*} 
where $\lim_{t\rightarrow+\infty}\mathbb{E}\left[I(C_{i_j}-N^{\phi_0}_j(t))\right]$ is the proportion of time that  $N^{\phi}_j(t) < C_{i_j}$.  The value $1-A_{\ell}$ represents the blocking probability of job class $\ell$ under the policy $\phi_0$.
Thus, we can further relax \eqref{eqn:constraint:relax:boundary_tmp} as
\begin{equation}\label{eqn:constraint:relax:boundary}
I\bigl(\bar{\alpha}_{\ell}^{\phi}\bigr) \leq I(1- A_{\ell}),~\forall \ell\in[L]
\end{equation}
where $\bar{\alpha}_{\ell}^{\phi}=\lim_{t\rightarrow+\infty}\mathbb{E}\bigl[\bar{a}^{\phi}_{\ell}(\bm{N}^{\phi}(t))\bigr]$.
Equations \eqref{eqn:opt}, \eqref{eqn:constraint:relax:action}, \eqref{eqn:constraint:relax:capacity} and \eqref{eqn:constraint:relax:boundary} comprise a \emph{relaxed} version of our original problem described in \eqref{eqn:opt}, \eqref{eqn:constraint:action}, \eqref{eqn:constraint:capacity} and \eqref{eqn:constraint:boundary}. Define $\tilde{\Phi}$ as the set of all the policies $\phi$ satisfying \eqref{eqn:constraint:relax:action}, \eqref{eqn:constraint:relax:capacity} and \eqref{eqn:constraint:relax:boundary} so that $\Phi\subset\tilde{\Phi}$.

For clarity of presentation, 
\begin{itemize}
\item 
define $\pi^{\phi}_{j}(n)$ as the steady state probability of state $n\in\mathscr{N}_{i_j}$ under policy $\phi\in\Phi$, and define $\bm{\pi}^{\phi}_{j}=(\pi^{\phi}_{j}(n):\ n\in\mathscr{N}_{i_j})$;
\item for $n\in\mathscr{N}_{i_j}$, $j\in[J]\cup\{j_{\ell}:\ell\in[L]\}$, $\ell\in[L]$, and $\phi\in\tilde{\Phi}$,
define 
$$\alpha^{\phi}_{\ell,j}(n)=\lim_{t\rightarrow +\infty}\mathbb{E}\bigl[a^{\phi}_{\ell,j}(\bm{N}^{\phi}(t))\ |\ N^{\phi}_{j}(t)=n\bigr];$$ 
\item for $\bm{\nu}\in\mathbb{R}^{L}$ and $\bm{\omega}\in\mathbb{R}^{LJ}$,
define 
\begin{equation}\label{eqn:reward_func:1}
r^{\phi}_{j,n}(\bm{\nu},\bm{\omega})=\mu_j(n)-e^*\varepsilon_j(n) - \sum\limits_{\ell:\ j\in\mathcal{J}_{\ell}}\nu_{\ell}\alpha^{\phi}_{\ell,j}(n)
\end{equation} 
if $n< C_{i_j}$;
\begin{equation}\label{reward_func:2}
r^{\phi}_{j,n}(\bm{\nu},\bm{\omega})=\mu_j(n)-e^*\varepsilon_j(n) - \sum\limits_{\ell:\ j\in\mathcal{J}_{\ell}}(\nu_{\ell}+\omega_{\ell,j})\alpha^{\phi}_{\ell,j}(n)
\end{equation}
otherwise;
and let 
$\bm{r}^{\phi}_j(\bm{\nu},\bm{\omega})=(r^{\phi}_{j,n}(\bm{\nu},\bm{\omega}):\ n\in\mathscr{N}_{i_j})$; 
\item
let $\bm{\mathcal{i}}^{\phi}=(I(\bar{\alpha}^{\phi}_{\ell}):\ell\in[L])$, $\bm{\mathcal{a}}^{\phi}=(\alpha_{\ell,j_{\ell}}^{\phi}: \ell\in[L])$ and $\bm{I}=(I(1-A_{\ell}):\ell\in[L])$.
\end{itemize}
The Lagrangian dual function of the relaxed problem is:
\begin{equation}\label{eqn:dual}
L(\bm{\nu},\bm{\omega},\pmb{\gamma})
= \max\limits_{\phi\in\tilde{\Phi}}~\sum\limits_{j\in[J]} \bm{\pi}_{j}^{\phi}\cdot \bm{r}^{\phi}_j(\bm{\nu},\bm{\omega})-\bm{\nu}\cdot\bm{\mathcal{a}}^{\phi}-\pmb{\gamma}\cdot {\bm{\mathcal{i}}}^{\phi} + \bm{\nu}\cdot \bm{1} + \pmb{\gamma}\cdot \bm{I},
\end{equation}
with Lagrangian multipliers $\bm{\nu}$, $\bm{\omega}$ and $\pmb{\gamma}$ corresponding to  \eqref{eqn:constraint:relax:action}, \eqref{eqn:constraint:relax:capacity} and \eqref{eqn:constraint:relax:boundary}, respectively.
Here, $\bm{x}\cdot\bm{y}$ is the inner product of vectors $\bm{x}$ and $\bm{y}$.

In the same vein of \cite{whittle1988restless}, for given multipliers $\bm{\nu}$, $\pmb{\gamma}$ and $\bm{\omega}$, the optimal solution of the maximization problem in \eqref{eqn:dual} is also optimal for
\begin{multline}\label{eqn:decompose}
\max\limits_{\phi\in\tilde{\Phi}}~\sum\limits_{j\in[J]} \bm{\pi}_{j}^{\phi}\cdot \bm{r}^{\phi}_j(\bm{\nu},\bm{\omega})-\bm{\nu}\cdot\bm{\mathcal{a}}^{\phi}-\pmb{\gamma}\cdot \bm{\mathcal{i}}^{\phi} 
=\sum\limits_{j\in[J]}\max\limits_{\bm{\alpha}^{\phi}_j\in[0,1]^{L|\mathscr{N}_{i_j}|}}\bm{\pi}_{j}^{\phi}\cdot \bm{r}^{\phi}_j(\bm{\nu},\bm{\omega})\\
+
\sum\limits_{\ell\in[L]}\max\limits_{\alpha^{\phi}_{\ell,j_{\ell}}\in[0,1]}(-\nu_{\ell}\alpha^{\phi}_{\ell,j_{\ell}}-\gamma_{\ell}I(\bar{\alpha}^{\phi}_{\ell})) ,
\end{multline}
where 
$$\bm{\alpha}^{\phi}_{j}=(\alpha^{\phi}_{\ell,j}(n):\ n\in\mathscr{N}_{i_j},~\ell\in[L]).$$
Since the constraints of action variables are now interpreted by the multipliers,
the problem described in \eqref{eqn:decompose} can be decomposed into $J+L$ independent sub-problems:
\begin{itemize}
\item for $j\in[J]$, 
\begin{equation}\label{eqn:subproblem:normal}
\max\limits_{\bm{\alpha}^{\phi}_j\in[0,1]^{L|\mathscr{N}_{i_j}|}}\bm{\pi}_{j}^{\phi}\cdot\bm{r}^{\phi}_j(\bm{\nu},\bm{\omega}),
\end{equation}
where $\bm{x}\cdot\bm{y}$ represents the inner product of vectors $\bm{x}$ and $\bm{y}$, and the steady state distribution $\bm{\pi}^{\phi}_j$ is determined by only the action vector $\bm{\alpha}^{\phi}_j$ associated with the underlying process $\{N^{\phi}_j(t),t\geq 0\}$;
similarly,
\item for $\ell\in[L]$,
\begin{equation}\label{eqn:subproblem:blocking}
\max\limits_{\alpha^{\phi}_{\ell,j_{\ell}}\in[0,1]}(-\nu_{\ell}\alpha^{\phi}_{\ell,j_{\ell}}-\gamma_{\ell}I(\bar{\alpha}^{\phi}_{\ell})) .
\end{equation}
\end{itemize}
Define $\tilde{\Phi}_1$ as the set of policies determined by action variables $\bm{\alpha}^{\phi}_j\in[0,1]^{L|\mathscr{N}_{i_j}|}$ ($j\in[J]$) and $\alpha^{\phi}_{\ell,j_{\ell}}\in[0,1]$ ($\ell\in[L]$).

{\bf Remark} The sub-problems in \eqref{eqn:subproblem:normal} and \eqref{eqn:subproblem:blocking} are independent problems,  each of which has only one-dimensional state space and thus experiences remarkably lower computation time than the original problem.
Nonetheless, 
similar to the discussions in \cite{whittle1988restless} and subsequent work about RMABP \cite{weber1990index, nino2001restless, nino2007dynamic, papadimitriou1999complexity,verloop2016asymptotically, fu2018restless}, non-trivial properties are generally required to establish theoretical connections between these one-dimensional sub-problems and the high-dimensional original problem. 
A detailed survey about RMABP has been provided in Section~\ref{sec:rWork}.

\subsection{Whittle Indexability}\label{subsec:indexability}

For a standard RMABP, 
Whittle \cite{whittle1988restless} proposed the well-known \emph{Whittle index policy} when a non-trivial property related to each bandit process was satisfied.  
This property was later referred to as the Whittle indexability. 
More precisely, for our problem defined in Section~\ref{sec:problem}, when $L=1$, the bandit process $\{N^{\phi}_j(t), t\geq 0\}$, for $j\in[J]$, reduces to a bandit process for a standard RMABP.
In this special case, based on \cite{whittle1988restless}, for each $j\in[J]$, if there exist an optimal solution $\phi^*$ for the problem described in \eqref{eqn:subproblem:normal} and a vector of real numbers $\bm{\upsilon}^*_j = (\upsilon^*_{\ell,j}(n): n\in\mathscr{N}_{i_j}, \ell\in[L])$, satisfying, for all $n\in \mathscr{N}_{i_j}\backslash\{C_{i_j}\}$,
\begin{equation}\label{eqn:indexability}
\alpha^{\phi^*}_{\ell,j}(n)=\left\{
\begin{cases}
1, & \text{if } \nu_{\ell} < \upsilon^*_{\ell,j}(n),\\
a, & \text{if }\nu_{\ell} =\upsilon^*_{\ell,j}(n),\\
0, & \text{otherwise},
\end{cases}\right.
\end{equation}
where $a$ can be any number in $[0,1]$ and $\ell=L=1$, then we say the bandit process $\{N^{\phi}_j(t), t\geq 0\}$ associated with $j\in[J]$ is Whittle indexable and the real number $\upsilon^*_{\ell,j}(n)$ is the \emph{Whittle index} for state $n\in\mathscr{N}_{i_j}$ of the process. If all the bandit processes of an RMABP are Whittle indexable, then the RMABP is Whittle indexable.
Note that, in \eqref{eqn:indexability}, the real number $\upsilon^*_{\ell,j}(n)$ must be independent from $\nu_{\ell}$.
Although the policy $\phi^*$ is optimal for the sub-problem described in \eqref{eqn:subproblem:normal}, which is usually not applicable to the original problem, the Whittle index $\upsilon^*_{\ell,j}(n)$ intuitively represents the marginal reward of taking $a^{\phi}_{\ell,j}(\bm{N}^{\phi}(t))=1$ when $N^{\phi}_j(t) = n$ and hence brings a bird's-eye view of approximating optimality for the original problem. 
Whittle \cite{whittle1988restless} proposed a scalable \emph{index policy} by prioritizing states $n\in\mathscr{N}_{i_j}$ of bandit processes $j\in[J]$ according to the descending order of their Whittle indices, which was  subsequently referred to as the Whittle index policy. 
In \cite{whittle1988restless}, Whittle conjectured asymptotic optimality of the Whittle index policy, and it was proved by Weber and Weiss \cite{weber1990index} under another non-trivial condition - the existence of a global attractor related to the original stochastic process. As indicated in Section~\ref{sec:intro}, asymptotic optimality acts as an important performance guarantee for scalable policies in large-scale optimization problems.

Recall that, in general, bandit processes are not necessarily Whittle indexable, so does an RMABP. 
Please refer to Section~\ref{sec:rWork} for a detailed survey of past studies on Whittle indexability.
Here, we focus on the server farm problem described in Section~\ref{sec:problem}.

\begin{definition}
We say a physical component $j\in[J]$ in cluster $i\in[I]$ is \emph{energy-efficiently unimodal} if,
for any $n_1,n_2\in\mathscr{N}_{i}\backslash\{C_i\}$ with $n_1<n_2$, 
\begin{multline}\label{eqn:unimodal}
\bigl(\mu_i(n_2+1)-\mu_i(n_2)\bigr)\bigl(\varepsilon_i(n_1+1)\mu_i(n_1)-\varepsilon_i(n_1)\mu_i(n_1+1)\bigr)\\\leq \bigl(\mu_i(n_1+1)-\mu_i(n_1)\bigr)\bigl(\varepsilon_i(n_2+1)\mu_i(n_2)-\varepsilon_i(n_2)\mu_i(n_2+1)\bigr)
\end{multline}
\end{definition}
Intuitively, the energy-efficient unimodality implies a mild relationship of the component energy efficiency, $\mu_i(n)/\varepsilon_i(n)$, in different states $n\in\mathscr{N}_i$: there is at most one bump on the curve of $\mu_i(n)/\varepsilon_i(n)$ as $n$ tends from $0$ to $C_i$.
For example, if there exists $n_1\in\mathscr{N}_i\backslash\{C_i,C_i-1\}$ with $\mu_i(n_1)=\mu_i(n_1+1)$, then \eqref{eqn:unimodal} indicates either $\varepsilon_i(n_1) = \varepsilon_i(n_1+1)$ or, for all $n_2=n_1+1,n_1+2,\ldots,C_i-1$, $\mu_i(n_2)=\mu_i(n_2+1)$; that is, the curve of $\mu_i(n)/\varepsilon_i(n)$  is flat as $n$ tends from $n_1$ to $C_i$.
For $n_1<n_2$, if $\mu_i(n_1)<\mu_i(n_1+1)$ and $\mu_i(n_2)=\mu_i(n_2+1)$, then \eqref{eqn:unimodal} holds all the time.
For any $n_1,n_2\in\mathscr{N}_i\backslash\{C_i\}$ with $n_1<n_2$,   $\mu_i(n_1)<\mu_i(n_1+1)$ and $\mu_i(n_2)<\mu_i(n_2+1)$,  from \eqref{eqn:unimodal},   if $\mu_i(n_1+1)/\varepsilon_i(n_1+1) \leq \mu_i(n_1)/\varepsilon_i(n_1)$ then $\mu_i(n_2+1)/\varepsilon_i(n_2+1) \leq \mu_i(n_2)/\varepsilon_i(n_2)$. 
As a consequence, the curve of $\mu_i(n)/\varepsilon_i(n)$  has at most one bump as $n$ tends from $0$ to $C_i$.
The energy-efficient unimodality is only a mild condition because for real-world computing components, such as CPUs and GPUs, energy efficiency increases with processing speed - $\mu_i(n)/\varepsilon_i(n)$ usually increases with $n$ \cite{mei2017survey,takouna2011accurate}.

Consider other examples of energy-efficient unimodality. If $\mu_i(0)=0$ and $\mu_i(n_1)=\mu_i(n_2)$ for any $n_1,n_2\in\mathscr{N}_i\backslash\{0\}$, then, whatever the power consumption rates are, each component in cluster $i$ is energy-efficiently unimodal.
If $\mu_i(n)=\mu n$ ($n\in\mathscr{N}_i$) for some $\mu\in\mathbb{R}_+$ (that is, the jobs on a component in cluster $i$ are processed in the same vein of an $M/M/C_i/C_i$ queue), then the energy-efficient unimodality holds if and only if, for any $n\in\mathscr{N}_i\backslash\{C_i, C_i-1\}$, $2\varepsilon_i(n+1)\leq \varepsilon_i(n)+\varepsilon_i(n+2)$.

\begin{proposition}\label{prop:relaxed_opt_precise}
When the job sizes are exponentially distributed and all computing components are energy-efficiently unimodal, if, for all $\ell\in[L]$, if $\nu_{\ell} = \nu \lambda_{\ell}$ for some $\nu\in\mathbb{R}$, then there exist $H\in\mathbb{R}$ and a policy $\phi^*$ satisfying, 
for  $n\in \mathscr{N}_{i_j}\backslash\{C_{i_j}\}$, $j\in\mathcal{J}_{\ell}$, $\ell\in[L]$, 
\begin{equation}\label{eqn:prop:relaxed_opt_precise}
\alpha^{\phi^*}_{\ell,j}(n)=\left\{
\begin{cases}
1, & \text{if } \nu_{\ell} < \max\limits_{\begin{subarray}~n'=n+1,\\n+2,\ldots,C_{i_j}\end{subarray}}\frac{\lambda_{\ell}}{\mu_{i_j}(n')}\Bigl(\mu_{i_j}(n')-e^* \varepsilon_{i_j}(n')-g^*_j(\bm{\nu})\Bigr),\\
a, & \text{if }\nu_{\ell} =\max\limits_{\begin{subarray}~n'=n+1,\\n+2,\ldots,C_{i_j}\end{subarray}}\frac{\lambda_{\ell}}{\mu_{i_j}(n')}\Bigl(\mu_{i_j}(n')-e^* \varepsilon_{i_j}(n')-g^*_j(\bm{\nu})\Bigr),\\
0, & \text{otherwise},
\end{cases}\right.
\end{equation}
where $a$ can be any value in $[0,1]$ and
\begin{equation}\label{eqn:prop:g_star}
g^*_j(\bm{\nu}) \coloneqq \max\limits_{\phi\in\tilde{\Phi}_1} \lim\limits_{T\rightarrow +\infty} \frac{1}{T}\mathbb{E}\int_0^T \Bigl(\mu_{i_j}(N^{\phi}_j(t)) - e^*\varepsilon^{\phi}_{i_j}(N^{\phi}_j(t)) - \sum\limits_{\ell\in[L]}\nu_{\ell} \alpha^{\phi}_{\ell,j}(N^{\phi}_j(t))\Bigr) dt,
\end{equation}
such that, for all $h>H$,  the policy $\phi^*$ is optimal for the maximization problem in the Lagrangian dual function $L(\bm{\nu},\bm{\eta},\pmb{\gamma})$ of  the relaxed problem. 
\end{proposition}

The proof of Proposition~\ref{prop:relaxed_opt_precise} is provided in Appendix~\ref{app:prop:relaxed_opt_precise}. Here, the real number $g^*_j(\bm{\nu})$ is equal to the maximized average reward gained by the process $\{N^{\phi}_j(t),t\geq 0\}$ for $j\in[J]$, where the reward rate for state $N^{\phi}_j(t)=n$ is $\mu_{i_j}(n) - e^*\varepsilon^{\phi}_{i_j}(n) - \sum\nolimits_{\ell\in[L]}\nu_{\ell} \alpha^{\phi}_{\ell,j}(n)$. For an MDP, the $g^*_j(\bm{\nu})$ in \eqref{eqn:prop:relaxed_opt_precise} is usually referred to as the \emph{attached criterion} used to translate the maximization of the average reward to the maximization of the expected cumulative reward of the same process \cite{ross1992applied}.

Equation~\eqref{eqn:prop:relaxed_opt_precise} has the same form as \eqref{eqn:indexability} except that $g^*_j(\bm{\nu})$ is dependent on $\bm{\nu}$. 
Unlike the simplified case discussed in \cite{fu2020energy}, the $g^*_j(\bm{\nu})$ cannot be offset or expressed in a closed form.
This dependence between $g^*_j(\bm{\nu})$ and $\bm{\nu}$ significantly complicates the analysis of indexability and the computation of the indices and prevents the same technique in \cite{fu2020energy} from being applied directly. Recall that Whittle indexability is sufficient for the existence of Whittle indices, used to construct the well-known Whittle index policy, but it does not necessarily hold for every bandit process.
Although the Whittle index policy or a similar policy based on the Whittle indices has great potential to perform closely to an optimal solution when the system is large, asymptotically optimal policies are not necessarily unique. 
For the purpose of the server farm problem, instead of seeking perfect Whittle indexability, we consider a less stringent property referred to as the \emph{asymptotic indexability}.
\begin{definition}\label{define:asym_indexability}
For $j\in[J]$, if there exist a real vector $\bm{\upsilon}^*_j\coloneqq (\upsilon^*_{\ell,j}(n): \ell\in[L],n\in\mathscr{N}_{i_j})$, a policy $\phi^*$ with action variables $\alpha^{\phi^*}_{\ell,j}(n)$ ($n\in \mathscr{N}_{i_j}\backslash\{C_{i_j}\}$ and $\ell\in[L]$) satisfying \eqref{eqn:indexability}, and $H>0$ such that, for all $h>H$, $\phi^*$ is optimal for the maximization problem in \eqref{eqn:subproblem:normal}, then we say that the process $\{N^{\phi}_j(t), t\geq 0\}$  is asymptotically  indexable with indices $\bm{\upsilon}^*_j$.
\end{definition}

\begin{proposition}\label{prop:relaxed_opt_explicit}
When the job sizes are exponentially distributed, if there exists $H>0$ such that, for all $h>H$, an optimal solution $\phi^*$ for the maximization problem in the Lagrangian dual function $L(\bm{\nu},\bm{\eta},\pmb{\gamma})$ of the relaxed problem 
exists and satisfies \eqref{eqn:prop:relaxed_opt_precise}, then, for any $j\in[J]$, the process $\{N^{\phi}_j(t), t\geq 0\}$ is asymptotically indexable.
\end{proposition}

The proof of Proposition~\ref{prop:relaxed_opt_explicit} is provided in Appendix~\ref{app:prop:relaxed_opt_explicit}. 
When  the process $\{N^{\phi}_j(t), t\geq 0\}$ reduces to a standard bandit process with $L=1$, asymptotic indexability indicates Whittle indexability with sufficiently large $h$.
Similar to the Whittle indices, in a large system, for $\ell\in[L]$, $j\in[J]$, and $n\in\mathscr{N}_{i_j}\backslash\{C_{i_j}\}$, the real number $\upsilon^*_{\ell,j}(n)$ represents the marginal reward of admitting a new job of class $\ell$ when there are $n$ jobs being served by physical component $j$. 
For $\ell\in[L]$, $j\in [J]$ and $n\in\mathscr{N}_{i_j}\backslash\{C_{i_j}\}$, we refer to $\upsilon^*_{\ell,j}(n)$ as the \emph{index} of state $n$ of the process $\{N^{\phi}_j(t), t\geq 0\}$ for job class $\ell$.
For a large system, asymptotic indexability ensures the existence of a \emph{threshold-style policy}, as described in \eqref{eqn:indexability}, that is optimal for the corresponding sub-problem.
Asymptotic indexability plays an important role in proposing an asymptotically optimal and scalable policy for the original problem defined in \eqref{eqn:opt}, \eqref{eqn:constraint:action},~\eqref{eqn:constraint:capacity} and \eqref{eqn:constraint:boundary}. 

\begin{corollary}\label{coro:relaxed_opt_explicit}
When the job sizes are exponentially distributed and all computing components are energy-efficiently unimodal, if, for all $\ell\in[L]$, if $\nu_{\ell} = \nu \lambda_{\ell}$ for some $\nu\in\mathbb{R}$, then, for any $j\in[J]$, the process $\{N^{\phi}_j(t), t\geq 0\}$ is asymptotically indexable.
\end{corollary}
Corollary~\ref{coro:relaxed_opt_explicit} is a straightforward result of Propositions~\ref{prop:relaxed_opt_precise} and \ref{prop:relaxed_opt_explicit} and provides a sufficient condition, related to only the offered traffic, service and energy consumption rates, for the asymptotic indexability.

\subsection{Existence of Indices}\label{subsec:exist_indices}

More importantly, we are interested in the exact values of the indices $\bm{\upsilon}^*$ that impose asymptotic indexability and further lead to a scalable, near-optimal policy for the original problem.

\begin{proposition}\label{prop:asym_opt:relaxed_opt}
When the job sizes are exponentially distributed and all computing components are energy-efficiently unimodal, if, for all $\ell\in[L]$, if $\nu_{\ell} = \nu \lambda_{\ell}$ for some $\nu\in\mathbb{R}$,
then, for any $j\in[J]$, the process $\{N^{\phi}_j(t), t\geq 0\}$ is asymptotically indexable with indices satisfying
\begin{equation}\label{eqn:asym_opt:index}
\upsilon^*_{\ell,j}(n)= \max\limits_{\begin{subarray}~n'=n+1,\\n+2,\ldots,C_i\end{subarray}}\lambda_{\ell}\Bigl(1-\frac{e^* \varepsilon_{i_j}(n') + \Gamma_j(\bm{\upsilon}_j^*(n))}{\mu_{i_j}(n')}\Bigr)
\end{equation}
where $\bm{\upsilon}^*_j(n) = (\upsilon^*_{\ell,j}(n): \ell\in[L])$, and 
\begin{equation}\label{eqn:asym_opt:Gamma}
\Gamma_j(\bm{\upsilon}_j^*(n))\coloneqq \max\limits_{n\in\mathscr{N}_{i_j}\backslash\{C_{i_j}\}}\bm{\pi}^{\psi_j(n)}_j\cdot \bm{r}^{\psi_j(n)}_j(\bm{\upsilon}_j^*(n),\cdot),
\end{equation}
with policy $\psi_j(n)$ satisfying, for all $\ell\in\{\ell'\in[L] | j\in\mathcal{J}_{\ell'}\}$ and $n'\in\mathscr{N}_{i_j}\backslash\{C_{i_j}\}$, $\alpha^{\psi_j(n)}_{\ell,j}(n')=1$ if $n'\leq n$; and $\alpha^{\psi_j(n)}_{\ell,j}(n')=0$ otherwise. 
\end{proposition}
The proposition is proved in Appendix~\ref{app:asym_opt:relaxed_opt}. 
The policy $\psi_j(n)$ satisfies constraints in~\eqref{eqn:constraint:relax:capacity} and thus is not affected by $\bm{\omega}$; we write the reward vector $\bm{r}^{\phi}_j(\bm{\nu},\bm{\omega})$ as $\bm{r}^{\phi}_j(\bm{\nu},\cdot)$ since it becomes independent of the second argument under $\psi_j(n)$.
The proposition is achieved by invoking Proposition~\ref{prop:relaxed_opt_precise} and approximating the value of $g^*_j(\bm{\nu})$ in \eqref{eqn:prop:relaxed_opt_precise}. In particular, if $\mu_i(n)/(\varepsilon_i(n)-\varepsilon_i(0))$ is increasing in $n\in\mathscr{N}_i\backslash\{0\}$, then the $\upsilon^*_{\ell,i}(n)$ satisfying \eqref{eqn:asym_opt:index} and \eqref{eqn:asym_opt:Gamma} reduces to a closed form, $\lambda_{\ell}(1-e^*/\mathcal{r}_i)$ where $\mathcal{r}_i = \mu_i(C_i)/(\varepsilon_i(C_i)-\varepsilon_i(0))$. 
Previous work discussed in \cite{fu2016asymptotic, fu2018restless} fell in the scope of this special case, for which, as mentioned earlier in Section~\ref{subsec:indexability}, the discussion can be  simplified by offsetting the attached criterion $g^*_j(\bm{\nu})$.
We refer to a detailed discussion for Proposition~\ref{prop:asym_opt:relaxed_opt} in Appendix~\ref{app:asym_opt:relaxed_opt}.

The index values $\upsilon^*_{\ell,j}(n)$ mentioned in Proposition~\ref{prop:asym_opt:relaxed_opt} can be obtained by solving Equations~\eqref{eqn:asym_opt:index} and \eqref{eqn:asym_opt:Gamma}. 
More precisely, observing \eqref{eqn:asym_opt:Gamma}, function $\Gamma_j(\bm{\nu})$ is dependent on $j\in[J]$ and $\bm{\nu}\in\mathbb{R}^L$ only through $i_j$ and the sum $\sum_{\ell: i_j\in\mathcal{I}_{\ell}}\nu_{\ell}$, respectively.
For notational convenience, we rewrite $\Gamma_j(\bm{\nu})$ as $\Gamma_j(\eta_j)$ with $\eta_j=\sum_{\ell:i_j\in\mathcal{I}_{\ell}}\nu_{\ell}$, and, for $\eta\in\mathbb{R}$, $i\in[I]$ and any $j\in\mathscr{J}_i$, define $\bar{\Gamma}^h_i(\eta^0)\coloneqq \Gamma_j(h\eta^0)$.
Let $\bar{\Gamma}^{h,\text{EXP}}_i(\eta^0)$  be equal to the value of $\bar{\Gamma}^h_i(\eta^0)$ with assumed exponentially distributed job sizes.
For $n\in\mathscr{N}_{i}\backslash\{C_{i}\}$ and $i\in [I]$, define
\begin{equation}\label{define:f}
f^h_{i,n}(\eta^0) \coloneqq \eta^0+ \hat{\lambda}^0_{i}\min\limits_{\begin{subarray}~n'=n+1,\\n+2,\ldots,C_i\end{subarray}}\Bigl(\frac{\bar{\Gamma}^{h,\text{EXP}}_i(\eta^0) +e^*\varepsilon_{i}(n')}{\mu_{i}(n')}
-1\Bigr)
\end{equation}
where $\hat{\lambda}^0_i=\sum_{\ell:i\in\mathcal{I}_{\ell}}\lambda^0_{\ell}$. 
For $n\in\mathscr{N}_i\backslash\{C_i\}$, $j\in\mathscr{J}_i$, $i\in\mathcal{I}_{\ell}$ and $\ell\in[L]$, given $\eta^0_{i,n}$ satisfying $f^h_{i,n}(\eta^0_{i,n}) = 0$,
we obtain that $\upsilon^*_{\ell,j}(n)= \bar{\upsilon}^*_{\ell,i}(n) \coloneqq h\eta^0_{i,n}\lambda^0_{\ell}/\hat{\lambda}^0_i$. 
Note that, from~\eqref{eqn:asym_opt:index} and \eqref{eqn:asym_opt:Gamma}, for any $i\in[I]$, $\upsilon^*_{\ell,j}(n)$ remains the same for all $j\in\mathscr{J}_i$.
In this context, solving~\eqref{eqn:asym_opt:index} and \eqref{eqn:asym_opt:Gamma} is equivalent to finding a zero point of function $f^h_{i,n}(\eta^0)$. 
Proposition~\ref{prop:asym_opt:relaxed_opt} ensures the existence of such a zero point when the sub-processes are asymptotically indexable. 
In Section~\ref{sec:algo}, in conjunction with proposing a scalable policy for the original server farm problem (described in \eqref{eqn:opt}, \eqref{eqn:constraint:action}, \eqref{eqn:constraint:capacity} and \eqref{eqn:constraint:boundary}), we will provide more detailed steps about computationally calculating the index values $\upsilon^*_{\ell,j}(n)$ that satisfy~\eqref{eqn:asym_opt:index} and \eqref{eqn:asym_opt:Gamma}.

Remarkably, in the general case without assuming asymptotic indexability, a solution also exists for $f^h_{i,n}(\eta^0)=0$, so do the indices $\upsilon^*_{\ell,j}(n)$ satisfying~\eqref{eqn:asym_opt:index} and \eqref{eqn:asym_opt:Gamma}. Consider the following proposition.
\begin{proposition}\label{prop:exist_g}
When job sizes are exponentially distributed, for $h\in\mathbb{N}_+\cup\{+\infty\}$,  $n\in\mathscr{N}_{i}\backslash\{C_{i}\}$, and $i\in [I]$, $f^h_{i,n}(\eta^0)$ is Lipschitz continuous in $\eta^0\in\mathbb{R}$, and there exists $\eta^0_{i,n}\in\mathbb{R}$ such that $f^h_{i,n}(\eta^0_{i,n})=0$. 
\end{proposition}
The proposition is proved in Appendix~\ref{app:exist_g}.

The value $\upsilon^*_{\ell,j}(n)$ is in fact a specific value of the multiplier $\nu_{\ell}$ considered in the process $\{N^{\phi}_j(t),t\geq 0\}$ ($j\in\mathscr{J}_i$) associated with state $n\in\mathscr{N}_i$ of cluster $i$, providing a bird's-eye view of the marginal reward of activating this process (that is, setting $\alpha^{\phi}_{\ell,j}(n)=1$).
For the original problem described in \eqref{eqn:opt}, \eqref{eqn:constraint:action}, \eqref{eqn:constraint:capacity} and \eqref{eqn:constraint:boundary}, we will propose a scalable policy in Section~\ref{sec:algo} that prioritizes components with higher $\upsilon^*_{\ell,j}(n)$ (derived by~\eqref{eqn:asym_opt:index} and \eqref{eqn:asym_opt:Gamma}), and we prove its asymptotic optimality under mild conditions related to the offered traffic with exponentially distributed job sizes and energy efficient unimodality.
When the hypothesis fails, extensive numerical results will be presented in Section~\ref{sec:simulation} to demonstrate the effectiveness of the proposed policy.

\section{Multiple Power Modes with Priorities (MPMP)}\label{sec:algo}
For the original problem described in \eqref{eqn:opt}, \eqref{eqn:constraint:action}, \eqref{eqn:constraint:capacity} and \eqref{eqn:constraint:boundary}, we propose a policy that prioritizes the components and jobs according to the descending order of the indices $\bm{\upsilon}^*\coloneqq (\upsilon^*_{\ell,j}(n): \ell\in[L],j\in[J],n\in\mathscr{N}_{i_j})$ satisfying \eqref{eqn:asym_opt:index} and \eqref{eqn:asym_opt:Gamma} for all $j\in[J]$.
These indices represent marginal rewards of serving the various jobs.  
More precisely, define a policy $\varphi$ satisfying
\begin{equation}\label{eqn:index_policy:1}
a^{\varphi}_{j,\ell}(\bm{N}^{\varphi}(t))=
\left\{\begin{array}{ll}
1,& \text{if } N^{\varphi}_j(t)<C_{i_j} \text{ and }\\
& j= \arg\max\limits_{j'\in\mathcal{J}_{\ell}}\Bigl[\frac{1}{h}\upsilon^*_{\ell,j'}(N^{\varphi}_j(t))\Bigr],\\
0,& \text{otherwise},
\end{array}\right.
\end{equation}
where $\upsilon^*_{\ell,j}(n)$ ($\ell\in[L]$,$j\in[J]$, $n\in\mathscr{N}_{i_j}$) is derived from \eqref{eqn:asym_opt:index} and \eqref{eqn:asym_opt:Gamma},
and
\begin{equation}\label{eqn:index_policy:2}
a^{\varphi}_{\ell,j_{\ell}}(\bm{N}^{\varphi}(t)) = 1 - \sum\limits_{j\in\mathcal{J}_{\ell}}a^{\varphi}_{j,\ell}(\bm{N}^{\varphi}(t)).
\end{equation}
Here, the $\frac{1}{h}$ before $\upsilon^*_{\ell,j'}(N^{\varphi}_j(t))$ is used to keep the value finite for all $h\in\mathbb{N}_+\cup\{+\infty\}$ and will not change the order of the indices $\bm{\upsilon}^*$.
In \eqref{eqn:index_policy:1}, tie-breaking rules can be arbitrary if $\arg\max$ returns more than one argument.
Note that all the theoretical results presented in this paper hold for arbitrary tie-breaking rules.
Nonetheless, beyond the theoretical results, the policies $\varphi$ with different tie-breaking rules can potentially have different performance.
In Section~\ref{sec:simulation}, for the numerical results, we will consider specific tie-breaking rules to complete the simulations and provide numerical setting details. 
Such a policy $\varphi$ is a feasible policy in $\Phi$; that is, the action variables determined by \eqref{eqn:index_policy:1} and \eqref{eqn:index_policy:2} satisfy constraints~\eqref{eqn:constraint:action}, \eqref{eqn:constraint:capacity} and \eqref{eqn:constraint:boundary}.
For implementation in a server farm system, $\varphi$ is scalable with computational complexity at most linear to the number of physical components for each arriving job.
Equation~\eqref{eqn:index_policy:1} indicates that, for an $\ell$-job newly arrived at time $t$, policy $\varphi$ always selects the component $j$ with the highest index value $\upsilon^*_{\ell,j}(N^{\varphi}_j(t))$, among all components $j\in\mathcal{J}_{\ell}$ with at least a vacant slot ($N^{\varphi}_j(t) < C_{i_j}$). If all components $j\in\mathcal{J}_{\ell}$ are fully occupied, then, from \eqref{eqn:index_policy:2}, the virtual component $j_{\ell}$ is selected to reject this job.
In Section~\ref{sec:asym_opt}, we prove in Proposition~\ref{prop:asym_opt} that $\varphi$ is asymptotically optimal under certain conditions.

\subsection{Indices with Known Criterion}\label{subsec:calculate_indices}

Recall that the server farm problem in \eqref{eqn:obj}, \eqref{eqn:constraint:action}, \eqref{eqn:constraint:capacity} and \eqref{eqn:constraint:boundary} has been translated to the problem described in \eqref{eqn:opt}, \eqref{eqn:constraint:action}, \eqref{eqn:constraint:capacity} and \eqref{eqn:constraint:boundary} by introducing a given real number $e^*\in\mathbb{R}$ that satisfies \eqref{eqn:e_star}. 
As mentioned in Section~\ref{subsec:exist_indices}, from Proposition~\ref{prop:exist_g}, with given $e^*\in\mathbb{R}$, the indices $\bm{\upsilon}^*$ satisfying \eqref{eqn:asym_opt:index} and \eqref{eqn:asym_opt:Gamma}  can always be obtained by solving $f^h_{i,n}(\eta^0)=0$ ($n\in\mathscr{N}_i\backslash\{C_i\}$, $i\in[I]$).
In particular, $\upsilon^*_{\ell,j}(n)= \bar{\upsilon}^*_{\ell,i_j}(n) =h\eta^0_{i_j,n}\lambda^0_{\ell}/\hat{\lambda}^0_{i_j}$ where $\eta^0_{i_j,n}$ satisfies $f^h_{i_j,n}(\eta^0_{i_j,n})=0$.
We propose the pseudo-code for such zero points of $f^h_{i,n}$ ($n\in\mathscr{N}_i\backslash\{C_i\}$, $i\in[I]$) in Algorithm~\ref{algo:nu_star}, where 
a bisection method is used to approximate the indices with a precision parameter $\epsilon>0$.
\begin{definition}
For $e\in\mathbb{R}$, $h\in\mathbb{N}_+$, $\ell\in[L]$, $i\in[I]$, $j\in\mathscr{J}_i$, and $n\in\mathscr{N}_i\backslash\{C_i\}$, let $u_{\ell,i}(n,e)$ represent an estimate of $\frac{1}{h}\bar{\upsilon}^*_{\ell,i}(n)$, or equivalently, an estimate of the index $\frac{1}{h}\upsilon^*_{\ell,j}(n)$ of state $n$ of the process $\{N^{\phi}_j(t),t\geq0\}$ for job class $\ell$ with given $e^*=e$, through Algorithm~\ref{algo:nu_star}.
Define $\bm{u}(e)\coloneqq (u_{\ell,i}(n,e): n\in\mathscr{N}_i\backslash\{C_i\},i\in\mathcal{I}_{\ell},\ell\in[L])$.
\end{definition}

If we simplify our problem by allowing $e^*$ to be any given non-negative real number, the resulting objective defined in \eqref{eqn:opt} becomes the maximization of the difference between the job throughput $\mathfrak{L}^{\phi}$ and the power consumption $\mathfrak{E}^{\phi}$ \emph{weighted} by the known  $e^*$. 
Such an objective is popular and has been widely considered and studied in the literature \cite{pedram2012energy,wu2017stochastic,fu2018restless}. 
In this simplified problem, based on later results in Section~\ref{sec:asym_opt},
the policy $\varphi$ with indices $\bm{\upsilon}^*$ approximated by $\bm{u}(e^*)$
is asymptotically optimal if  $L=1$ or Condition~\ref{cond:heavy} holds for an energy-efficiently unimodal system.

\begin{algorithm}[t]
\linespread{0.8}\selectfont

\SetKwFunction{FIndices}{Indices}
\SetKwProg{Fn}{Function}{:}{\KwRet}
\SetKwInOut{Input}{Input}\SetKwInOut{Output}{Output}
\SetAlgoLined
\DontPrintSemicolon
\Input{Real value $e^*$.}
\Output{Vectors $\bm{u}(e^*)$ and $\bm{\eta}^0$.}
\Fn{\FIndices{$e^*$}}{
   $\eta^0_{i,n}\gets 0$ for all $i\in[I]$ and $n\in\mathscr{N}_i$.\;
   \For {$i\in[I]$}{
   	   $\eta_1 = 1$\;
       $\eta_2\gets \hat{\lambda}^0_i\left(1-\frac{e^*(\varepsilon_i(C_i)-\varepsilon_i(0))}{\mu_i(C_i)}\right)$\;
       $\eta_3\gets\eta_2$\;
       \While {$f^h_{i,C_i-1}(\eta_2 - \eta_1) > 0$}{
    	    $\eta_2\gets \eta_2 - \eta_1$\;
		    $\eta_1 \gets 2\eta_1$\;
   		} 
       $\eta_1 \gets \eta_2-\eta_1$\;       
		\tcc{\small $\eta_1$ and $\eta_2$ are the lower and upper bounds of $\eta^0_{i,C_i-1}$.}
	    \For {$n=C_i-1$ {\bf to} $0$}{
     		\While{$\eta_2-\eta_1 < \epsilon$}{
	       		$\eta_{i,n} \gets (\eta_2+\eta_1)/2$\;
	           \uIf {$f^h_{i,n}(\eta_{i,n}) < 0$}{
    	          		$\eta_1 \gets \eta_{i,n}$\;
		        }\uElseIf {$f^h_{i,n}(\eta_{i,n}) > 0$}{
	    	          $\eta_2 \gets \eta_{i,n}$\;
				 }\Else{
	     	          {\bf Break}\;
	      		 }
		    }
	       $\eta^0_{i,n} \gets (\eta_2+\eta_1)/2$\;
           $\eta_1 \gets \eta^0_{i,n}$\;
           $\eta_2 \gets \max\Bigl\{\hat{\lambda}^0_i\bigl(1-\frac{e^*(\varepsilon_i(n)-\varepsilon_i(0))}{\mu_i(n)}\bigr),~\eta_3\Bigr\}$\;
           $\eta_3\gets \eta_2$\;
		}
    }

    $u_{\ell,i}(n,e^*) \gets \eta^0_{i,n}\lambda^0_{\ell}/\hat{\lambda}^0_i$\;
}
\caption{Indices for given $e^*$.}\label{algo:nu_star}
\end{algorithm}

\subsection{Approximating the Unknown Criterion}\label{subsec:ratio}

The value of $e^*$ satisfying \eqref{eqn:e_star} is not known a priori.
For the purpose of this paper, we need to obtain this specific $e^*$, for which the optimal solution maximizing objective \eqref{eqn:opt} also maximizes \eqref{eqn:obj}.
From the definition of $e^*$ in \eqref{eqn:e_star}, we obtain
\begin{equation}\label{eqn:cond:e_star}
\max\limits_{\phi\in\Phi} \frac{1}{h}\Bigl(\mathfrak{L}^{\phi} -e^*\mathfrak{E}^{\phi}\Bigr) = 0,
\end{equation}
where the maximization is subject to \eqref{eqn:constraint:action}, \eqref{eqn:constraint:capacity} and \eqref{eqn:constraint:boundary}, and the $\frac{1}{h}$ is used to keep the value at the left-hand side of \eqref{eqn:cond:e_star} finite for all $h\in\mathbb{N}_+\cup\{+\infty\}$.
Following similar ideas of \cite{rosberg2014}, 
$e^*$ can be approximated by a bisection method with stopping condition \eqref{eqn:cond:e_star}.
Nonetheless, as mentioned in Section~\ref{sec:problem}, 
because of the complexity of solving the maximization in the left-hand side of \eqref{eqn:cond:e_star} (or the right-hand side of \eqref{eqn:e_star}), the optimal solutions are intractable.
We thus resort to scalable techniques that effectively approximate $e^*$.

\begin{definition}
For $e\in\mathbb{R}$, $h\in\mathbb{N}_+$, $\ell\in[L]$, $i\in[I]$, and $n\in\mathscr{N}_i\backslash\{C_i\}$, let $u^*_{\ell,i}(n,e)\in\mathbb{R}$ represent a solution of $f^h_{i,n}\bigl(\frac{\hat{\lambda}^0_i}{\lambda^0_{\ell}}u^*_{\ell,i}(n,e)\bigr)=0$ with substituted $e$ for $e^*$. Let $\bm{u}^*(e)\coloneqq (u^*_{\ell,i}(n,e): \ell\in[L],i\in[I],n\in\mathscr{N}_i\backslash\{C_i\})$.
\end{definition}
From Proposition~\ref{prop:exist_g}, $\bm{u}^*(e)$ also represents a vector of the estimates $\bm{u}(e)$ obtained by Algorithm~\ref{algo:nu_star} in the ideal case with the precision parameter $\epsilon\downarrow 0$.
For the system with scaling parameter $h\in\mathbb{N}_+$, we consider a policy $\psi^h(\bm{\nu},\bm{a}^h,e)\in\tilde{\Phi}_1$ with given $\bm{\nu}\in\mathbb{R}^L$, $e\in\mathbb{R}_0$ and $\bm{a}^h\in\{0,1\}^J$, satisfying, for $\ell\in[L]$, $j\in\mathscr{J}_{\ell}$, and $n\in\mathscr{N}_{i_j}\backslash\{C_{i_j}\}$,
\begin{equation}\label{eqn:index_style_policy:1}
\alpha^{\psi^h(\bm{\nu},\bm{a}^h,e)}_{\ell,j}(n) =\left\{\begin{cases}
1, & \text{if } \sum_{\ell': j\in\mathcal{J}_{\ell'}}\nu_{\ell'}< \sum_{\ell': j\in\mathcal{J}_{\ell'}}u^*_{\ell',i_j}(n,e),\\
a^h_j, & \text{if } \sum_{\ell': j\in\mathcal{J}_{\ell'}}\nu_{\ell'}= \sum_{\ell':j\in\mathcal{J}_{\ell'}}u^*_{\ell',i_j}(n,e),\\
0, &\text{otherwise,}
\end{cases}\right.
\end{equation}
$\alpha^{\psi^h(\bm{\nu},\bm{a}^h,e)}_{\ell,j}(C_{i_j})=0$ for all $\ell\in[L]$ and $j\in\mathscr{J}_{\ell}$, and, for all $\ell\in[L]$, 
\begin{equation}\label{eqn:index_style_policy:2}
\alpha^{\psi^h(\bm{\nu},\bm{a}^h,e)}_{\ell,j_{\ell}} = \max\Bigl\{0,1-\sum\limits_{j\in\mathcal{J}_{\ell}}\sum\limits_{n\in\mathscr{N}_{i_j}}\pi^{\psi^h(\bm{\nu},\bm{a}^h,e)}_j(n)\alpha^{\psi^h(\bm{\nu},\bm{a}^h,e)}_{\ell,j}(n)\Bigr\},
\end{equation}
where $\pi^{\phi}_j(n)$ is defined in Section~\ref{subsec:relaxation} and represents the steady state probability of state $n\in\mathscr{N}_{i_j}$ under a policy $\phi$.
For such a policy $\psi^h(\bm{\nu},\bm{a}^h, e)$, the process $\{N^{\psi^h(\bm{\nu},\bm{a}^h, e)}_j(t), t\geq 0\}$ is a birth-and-death process with state transition rates linear to $h$ and finitely many states, leading to the existence of $\lim\nolimits_{h\rightarrow +\infty} \bm{\pi}^{\psi^h(\bm{\nu},\bm{a}^h,e)}_j$. 
Define, for $\ell\in[L]$, 
\begin{equation}\label{eqn:sum_alpha}
A^{h,\psi^h(\bm{\nu},\bm{a}^h,e)}_{\ell} \coloneqq \sum\limits_{j\in \mathscr{J}_{\ell}\cup \{j_{\ell}\}} \sum\limits_{n\in\mathscr{N}_{i_j}}\pi^{\psi^h(\bm{\nu},\bm{a}^h,e)}_j(n)\alpha^{\psi^h(\bm{\nu},\bm{a}^h,e)}_{\ell,j}(n),
\end{equation}
which is the expected average sum of the action variables, namely the left-hand side of \eqref{eqn:constraint:relax:action}, under the policy $\psi^h(\bm{\nu},\bm{a}^h,e)$ in the asymptotic regime.
Since $A^{\psi^h(\bm{\nu},\bm{a}^h,e)}_{\ell}$ is decreasing in $\bm{\nu}$ and increasing in $\bm{a}^h$, there exist $\bm{\nu}\in\mathbb{R}^L$ and $\bm{a}^h\in\{0,1\}^J$ such that $\lim\nolimits_{h\rightarrow \infty}A^{h,\psi^h(\bm{\nu},\bm{a}^h,e)}_{\ell} = 1$ for all $\ell\in[L]$; that is, constraints~\eqref{eqn:constraint:relax:action} are satisfied by substituting the policy $\psi^h(\bm{\nu},\bm{a}^h,e)$ for $\phi$. More precisely, let $\mathscr{V}$ represent the set of $ (\bm{\nu},\bm{a}^h)$ such that $\lim\nolimits_{h\rightarrow \infty}A^{h,\psi^h(\bm{\nu},\bm{a}^h,e)}_{\ell} = 1$. Define
\begin{equation}
(\bar{\bm{\nu}},\bar{\bm{a}}^h) = \arg\max\limits_{(\bm{\nu},\bm{a}^h)\in\mathscr{V}} \lim\limits_{h\rightarrow \infty}\sum\limits_{j\in\mathcal{J}_{\ell}}\sum\limits_{n\in\mathscr{N}_{i_j}}\pi^{\psi^h(\bm{\nu},\bm{a}^h,e)}_j(n)\alpha^{\psi^h(\bm{\nu},\bm{a}^h,e)}_{\ell,j}(n).
\end{equation}
For $\phi\in\tilde{\Phi}$, $h\in\mathbb{N}_+$ and $e\in\mathbb{R}$, define
\begin{equation}\label{eqn:Gamma_e:limits}
\Gamma^{h,\phi}(e) \coloneqq \frac{1}{h}\sum\limits_{j\in[J]}\bm{\pi}^{\phi}_j\cdot\bm{r}^{\phi,e}_j,
\end{equation}
where $\bm{x}\cdot\bm{y}$ represents the inner product of vectors $\bm{x}$ and $\bm{y}$, and
\begin{equation}\label{eqn:define_r}
\bm{r}^{\phi,e}_j\coloneqq \Bigl(\mu_{i_j}(n)-e\varepsilon_{i_j}(n): n\in\mathscr{N}_{i_j}\Bigr).
\end{equation} 
The policy $\psi^h(\bar{\bm{\nu}},\bar{\bm{a}}^h,e)$ is applicable to the relaxed problem but not necessarily to the original problem. We will discuss in Section~\ref{sec:asym_opt} that, when job sizes are exponentially distributed, for any given real number $e^*$, $\Gamma^{h,\varphi}(e^*)$ converges to $\Gamma^{h,\psi^h(\bar{\bm{\nu}},\bar{\bm{a}}^h,e^*)}(e^*)$ as $h\rightarrow +\infty$. 
In the asymptotic regime, if the policy $\psi^h(\bar{\bm{\nu}},\bar{\bm{a}}^h,e^*)$ is coincidently optimal for the relaxed problem, it must also be optimal for the original problem because $\Gamma^{h,\varphi}(e^*) \leq \max_{\phi\in\Phi}\Gamma^{h,\phi}\leq\max_{\phi\in\tilde{\Phi}}\Gamma^{h,\phi}$. In this case, due to the simplicity of computing $\Gamma^{h,\psi^h(\bar{\bm{\nu}},\bar{\bm{a}}^h,e)}(e)$, we can approximate the value of $e^*$ satisfying \eqref{eqn:e_star} by utilizing the condition in \eqref{eqn:cond:e_star}.

For $e\in\mathbb{R}$, let 
\begin{equation}
\Gamma(e) \coloneqq \lim_{h\rightarrow+\infty}\Gamma^{h,\psi^h(\bar{\bm{\nu}},\bar{\bm{a}}^h,e)}(e).
\end{equation}
\begin{proposition}\label{prop:e_star_equal}
When job sizes are exponentially distributed,  
if, for a given $e\in\mathbb{R}$, 
\begin{equation}\label{eqn:e_star_equal:condition}
\lim\limits_{h\rightarrow\infty} \left\lvert \Gamma^{h,\psi^h(\bar{\bm{\nu}},\bar{\bm{a}}^h,e)}(e) - \max\limits_{\phi\in\tilde{\Phi}}\Gamma^{h,\phi}(e)\right\rvert = 0,
\end{equation}
then $\Gamma(e)$ is Lipschitz continuous and piece-wise linear in $e\in\mathbb{R}$, there exists a unique solution $e_0\in\mathbb{R}$ for 
\begin{equation}\label{eqn:exist_zero_point:gamma}
\Gamma(e_0)=0,
\end{equation} 
and, for the unique $e_0$ satisfying \eqref{eqn:exist_zero_point:gamma},
\begin{equation}\label{eqn:e_star_equal}
e_0 = \lim_{h\rightarrow +\infty}e^* \geq 0.
\end{equation}
\end{proposition}
The proposition is proved in Appendix~\ref{app:prop:e_star_equal}. Note that $e^*$ defined in \eqref{eqn:e_star} is dependent on the scaling parameter $h$. 
The proposition indicates that, if, for any given $e\in\mathbb{R}$, the policy $\psi^h(\bar{\bm{\nu}},\bar{\bm{a}}^h,e)$ is optimal for the relaxed problem described in \eqref{eqn:opt},  \eqref{eqn:constraint:relax:action}, \eqref{eqn:constraint:relax:capacity} and \eqref{eqn:constraint:relax:boundary} in the asymptotic regime (that is, \eqref{eqn:e_star_equal:condition} is satisfied), then $e^*$ can be approximated by the zero point $e_0$ asymptotically.
When \eqref{eqn:e_star_equal:condition} does not hold, $\Gamma(e)$ may be discontinuous at some points and there may not exist $e_0$ satisfying \eqref{eqn:exist_zero_point:gamma}. For clarify, define 
\begin{equation}\label{eqn:define_e0}
e_0 \coloneqq  \inf\{e\in \mathbb{R}~|~\forall e' > e,~\Gamma(e') < 0\}.
\end{equation}
From Proposition~\ref{prop:e_star_equal}, if job sizes are exponentially distributed and \eqref{eqn:e_star_equal:condition} holds, then $e_0$ defined in \eqref{eqn:define_e0} coincides with the unique solution for \eqref{eqn:exist_zero_point:gamma}.

Consider the following situation related to the traffic load.
\begin{condition}\label{cond:heavy}
$A_{\ell}\leq 1$ in the asymptotic regime for all $\ell\in[L]$.
\end{condition} 
Condition~\ref{cond:heavy} implies that the system is in \emph{heavy traffic} with positive blocking probabilities of all job classes under policy $\phi_0$ in the asymptotic regime. 
\begin{lemma}\label{lemma:convergence_relaxed_opt}
When job sizes are exponentially distributed, if the components are energy-efficiently unimodal and Condition~\ref{cond:heavy} or $L=1$, then \eqref{eqn:e_star_equal:condition} holds.
\end{lemma}
The lemma is based on Propositions~\ref{prop:asym_opt:relaxed_opt} and is proved in Appendix~\ref{app:lemma:convergence_relaxed_opt}.

\begin{corollary}\label{coro:e_star_equal}
When job sizes are exponentially distributed,
if the computing components are energy-efficiently unimodal, and Condition~\ref{cond:heavy} holds or $L=1$, then \eqref{eqn:e_star_equal} holds, $e_0$ defined in \eqref{eqn:define_e0} is the unique solution for \eqref{eqn:exist_zero_point:gamma}, and the function $\Gamma(e)$ is Lipschitz continuous and piece-wise linear in $e\in\mathbb{R}$.
\end{corollary}
The corollary is a direct result of Lemma~\ref{lemma:convergence_relaxed_opt} and Proposition~\ref{prop:e_star_equal}.

\begin{algorithm}[t]
\linespread{1}\selectfont

\SetKwFunction{FMin}{Gamma}
\SetKwProg{Fn}{Function}{:}{\KwRet}
\SetKwInOut{Input}{Input}\SetKwInOut{Output}{Output}
\SetAlgoLined
\DontPrintSemicolon
\Input{Real value $e$.}
\Output{$\Gamma(e)$}
\Fn{\FMin{$e$}}{
	$\bar{\bm{v}}\gets \lim_{h\rightarrow \infty} \bm{u}^*(e)$  and 
	$\eta_{i,n} \gets \sum_{\ell:i\in\mathcal{I}_{\ell}} \bar{v}_{\ell,i}(n)$\;
	Rank the state-component (SC) pairs $(n,i)$ ($i\in\mathscr{N}$, $n\in\mathscr{N}_i\backslash\{C_i\}$) according to the descending order of $\eta_{i,n}$; if $\eta_{i,n+1} = \eta_{i,n}$ for some $(i,n)$, rank $(i,n)$ prior to $(i,n+1)$.
    Let $(i_k,n_k)$ represent the $k$th SC pair.\;
    $s_{\ell}\gets 0$ for all $\ell\in[I]$ and $k\gets 1$\;
    $q_i\gets 0$ for all $i\in[I]$ and $u_i\gets 0$ for all $i\in[I]$\;
	$\hat{\lambda}^0_i \gets \sum_{\ell\in \{\ell'\in[L] | i\in\mathcal{I}_{\ell'}\}}\lambda^0_i$ for all $i\in[I]$\;
    \While {$\exists \ell\in[L], s_{\ell} < 1 \And k\leq \sum_{i\in[I]}C_i$}{
		$x\gets \arg\min_{\ell:i_k\in\mathcal{I}_{\ell}}(1-s_{\ell})$\;
		\If {$s_{x} < 1$}{
			\uIf {$M_{i_k}^0\bigl(\mu_{i_k}(n_k+1)-\mu_{i_k}(n_k)\bigr) \geq \hat{\lambda}^0_{i_k}(1-s_x)$}{
					$s_{\ell} \gets s_{\ell}+1-s_x$ for all $\ell$ with $i_k\in\mathcal{I}_{\ell}$\;
					$u_{i_k} \gets \frac{\hat{\lambda}^0_{i_k}(1- s_x)}{M_{i_k}^0\bigl(\mu_{i_k}(n_k+1)-\mu_{i_k}(n_k)\bigr)}$\;
					$q_{i_k} \gets n_k$\;
			}\Else{
				$s_{\ell} \gets s_{\ell}+\frac{M_{i_k}^0\bigl(\mu_{i_k}(n_k+1)-\mu_{i_k}(n_k)\bigr)}{\hat{\lambda}^0_{i_k}}$ for all $\ell$ with $i_k\in\mathcal{I}_{\ell}$\;
				$u_{i_k} \gets 0$ and $q_{i_k} \gets n_k+1$\;
			}
 		}
		$k\gets k+1$\;
   	}
	$\Gamma(e)\gets 0$\;
	\For {$i\in[I]$}{
		$r_1 \gets \bigl(\mu_i(q_i)-e\varepsilon_i(q_i)\bigr)(1-u_i)$\;
	 	$r_2 \gets \bigl(\mu_i(q_i+1)-e\varepsilon_i(q_i+1)\bigr)u_i$\;
		$\Gamma(e)\gets \Gamma(e) +(r_1+r_2)M_i^0$\;
	}
}
\caption{Fitting the value of $\Gamma(e)$.}\label{algo:opt3}
\end{algorithm}

\begin{algorithm}[t]
\linespread{0.8}\selectfont

\SetKwFunction{FPar}{ParameterForEnergyEfficiency}
\SetKwProg{Fn}{Function}{:}{\KwRet}
\SetKwInOut{Input}{Input}\SetKwInOut{Output}{Output}
\SetAlgoLined
\DontPrintSemicolon
\Output{An estimate of $e_0$, denoted by $\bar{e}_0$.}
\Fn{\FPar}{
	$e_1 \gets 0$ and
	$e_2 \gets \sum_{i\in[I]}M_i^0\max_{n\in\mathscr{N}_{i}}\frac{\mu_{i}(n)}{\varepsilon_{i}(n)}$\;
	\While {$e_2-e_1> \epsilon$}{
		$e \gets (e_2+e_1)/2$\;
		$\Gamma\gets$\text{Gamma(e)}
		\tcc*{\small Call Algorithm~\ref{algo:opt3}.}
    	\uIf {$\Gamma< 0$}{
			$e_2\gets e$\;
		}\uElseIf {$\Gamma> 0$} {
 			$e_1\gets e$\;
		}\Else{
 		   {\bf Break}\;
		}
	}
	$\bar{e}_0\gets (e_2+e_1)/2$\;
	\Return $\bar{e}_0$\;
}
\caption{Approximating Parameter $e^*$.}\label{algo:e_star}
\end{algorithm}

Given the simple form of $\psi^h(\bm{\nu},\bm{a}^h,e)$ defined in \eqref{eqn:index_style_policy:1} and \eqref{eqn:index_style_policy:2}, we can directly fit the values of $\bar{\bm{\nu}}\in\mathbb{R}^L$ and $\Gamma(e)$ to satisfy \eqref{eqn:constraint:relax:action}.
A pseudo-code for fitting the values of $\bar{\bm{\nu}}\in\mathbb{R}^L$ and $\Gamma(e)$ is provided in Algorithm~\ref{algo:opt3}.
We can also estimate the value of $e_0$ through a bisection method with the precision parameter $\epsilon>0$, for which a pseudo-code is provided in the Algorithm~\ref{algo:e_star}. Let $\bar{e}_0$ represent such an estimate of $e_0$.
Lemma~\ref{lemma:convergence_relaxed_opt} and Proposition~\ref{prop:e_star_equal} guarantee that the estimate $\bar{e}_0$ is within $[e_0-\epsilon,e_0+\epsilon]$, where $e_0$ is the unique zero point of  $\Gamma(e_0) = 0$, and $e_0$ is equal to $e^*$ in the asymptotic regime, when the computing components are energy-efficiently unimodal and Condition~\ref{cond:heavy} holds or $L=1$.

\begin{definition}
For $e_0$ defined in \eqref{eqn:define_e0},
we construct a policy $\varphi(e_0)$ by substituting $\varphi(e_0)$ and $u^*_{\ell,i_j}(n,e_0)$ for $\varphi$ and $\upsilon^*_{\ell,j}(n)$, respectively,  in \eqref{eqn:index_policy:1} and \eqref{eqn:index_policy:2}. 
We refer to it as the \emph{Multiple Power Modes with Priorities (MPMP)} policy.
\end{definition}
Similar to $\varphi$, defined in \eqref{eqn:index_policy:1} and \eqref{eqn:index_policy:2}, the policy $\varphi(e_0)$ prioritizes physical components according to the descending order of $\bm{u}^*(e_0)$ and approaches $\varphi$ in the asymptotic regime under certain conditions.
This policy is applicable to the original problem, described in \eqref{eqn:opt}, \eqref{eqn:constraint:action}, \eqref{eqn:constraint:capacity} and \eqref{eqn:constraint:boundary}.
From Corollary~\ref{coro:e_star_equal}, if components are energy-efficiently unimodal, and Condition~\ref{cond:heavy} holds or $L=1$,  MPMP asymptotically approaches $\varphi$, which will be proved to be optimal in the asymptotic regime in Proposition~\ref{prop:asym_opt}. 
Recall that, under the same conditions, from Corollary~\ref{coro:e_star_equal}, the estimate $\bar{e}_0$ output by Algorithm~\ref{algo:e_star} is within $[e_0-\epsilon,e_0+\epsilon]$ with $\lvert\Gamma(\bar{e}_0)-\Gamma(e_0)\rvert \leq C\epsilon$, where $\epsilon>0$ is the precision parameter used in Algorithm~\ref{algo:e_star} and $C$ is the Lipschitz constant for $\Gamma(e)$.
For cases without assuming energy-efficient unimodality, Condition~\ref{cond:heavy} or $L=1$, numerical results are presented in Section~\ref{sec:simulation} to demonstrate the effectiveness of MPMP.

In Algorithm~\ref{algo:MPMP}, we provide a pseudo-code of implementing MPMP for each arrived job. The implementation is not unique and can be further speed up by utilizing sensible data structure, such as the \emph{maxheap}. 
In Algorithm~\ref{algo:MPMP}, as an example, we select the component with the smallest label in the tie case for implementing MPMP. As explained in Section~\ref{sec:algo}, all the theoretical results presented in this paper apply to arbitrary tie-breaking rules, but, for simulation purpose, a specific rule is necessary. Although different tie-breaking rules potentially lead to different performance, based on Proposition~\ref{prop:exp_diminishing}, when job sizes are exponentially distributed, MPMP with any tie-breaking rule will quickly approach the same performance as $h$ becomes large.

Recall that, based on the definition of $\varphi(e_0)$ (that is, MPMP) in the above paragraph, the indices used for MPMP is expected to be $\bm{u}^*(e_0)$, which is defined as an ideal output of Algorithms~\ref{algo:nu_star} with the precision parameter $\epsilon\rightarrow 0$. All the theoretical results presented in this paper apply to MPMP, or equivalently, $\varphi(e_0)$, in this ideal case.
For the numerical results in Section~\ref{sec:simulation}, 
with slightly abused notation,
we still refer to the policy output by Algorithms~\ref{algo:nu_star}-\ref{algo:MPMP} as the MPMP policy, when the indices are estimated by Algorithms~\ref{algo:nu_star} and~\ref{algo:e_star} with $\epsilon = 10^{-15}$.

For any $e\in\mathbb{R}$ and its estimate $e+\epsilon$ with a small $\epsilon>0$, the indices $\bm{u}^*(e)$ and $\bm{u}^*(e+\epsilon)$ may be completely different. Nonetheless, even if  $\lVert\bm{u}^*(e)-\bm{u}^*(e+\epsilon)\rVert$ is large, the performance of the resulting policies can still be negligible. 
For $e\in\mathbb{R}$, define policy $\varphi(e)$ by substituting $\varphi(e)$ and $u^*_{\ell,i_j}(n,e)$ for $\varphi$ and $\upsilon^*_{\ell,j}(n)$, respectively,  in \eqref{eqn:index_policy:1} and \eqref{eqn:index_policy:2}.
\begin{proposition}\label{prop:varphi_closeness}
When job sizes are exponentially distributed, for $e\in\mathbb{R}$, if \eqref{eqn:e_star_equal:condition} holds, then, for any $\epsilon>0$, there exists a constant $C>0$ such that
\begin{equation}\label{eqn:prop:varphi_closeness}
\lim\limits_{h\rightarrow+\infty}\Bigl\lvert \Gamma^{h,\varphi(e)}(e)-\Gamma^{h,\varphi(e+\epsilon)}(e+\epsilon)\Bigr\rvert \leq C \epsilon.
\end{equation}
\end{proposition}
The proposition is proved in Appendix~\ref{app:prop:varphi_closeness}. Together with Lemma~\ref{lemma:convergence_relaxed_opt}, if the components are energy-efficiently unimodal, and Condition~\ref{cond:heavy} holds or $L=1$, then \eqref{eqn:prop:varphi_closeness} holds. That is, although $\lVert\bm{u}^*(e)-\bm{u}^*(e+\epsilon)\rVert$ may be  sensitive to $e$, the performance deviation of the resulting policies can still be negligible for sufficiently small $\epsilon$.
In particular, when the estimate $\bar{e}_0$ of $e_0$ is subject to a small $\epsilon>0$, the performance deviation between policies $\varphi(\bar{e}_0)$ and $\varphi(e_0)$, or equivalently the MPMP policy, is also bounded by $C\epsilon$.

\begin{algorithm}[t]
\linespread{0.8}\selectfont
\SetKwFunction{FImp}{ImplementingMPMP}
\SetKwProg{Fn}{Function}{:}{\KwRet}
\SetKwInOut{Input}{Input}\SetKwInOut{Output}{Output}
\SetAlgoLined
\DontPrintSemicolon
\Input{The class label $\ell$ of an arrived job, the indices $u^*_{\ell,j}(n)\in \mathbb{R}$ for $j\in\mathscr{J}_{\ell}$ and $n\in\mathscr{N}_{i_j}$,  and current system state $\bm{N}^{\text{MPMP}}(t)$ upon the arrival.}
\Output{The selected component $j$ in $\mathscr{J}_{\ell}\cup\{j_{\ell}\}$ to accommodate this job.}
\Fn{\FImp}{
    $j\gets j_{\ell}$ and 
    $u\gets -\infty$\;
	\For {$\forall j'\in\mathscr{J}_{\ell}$}{
		\If {$N^{\text{MPMP}}_{j'}(t) < C_{i_{j'}}$ AND $u^*_{\ell,j'}(N^{\text{MPMP}}_{j'}(t)) > u$}{
		    $j\gets j'$ and 
		    $u\gets u^*_{\ell,j'}(N^{\text{MPMP}}_{j'}(t))$\;
		}
	}
 %   \Return $j$;
}
\caption{Implementing MPMP with given indices.}\label{algo:MPMP}
\end{algorithm}

The MPMP policy is scalable and applicable to the original server farm problem.
Let $C=\sum_{i\in[I]}C_i$. 
The computational complexity of computing $\bm{u}(\bar{e}_0)$ is $O\bigl(P_2(P_1C+C\ln C+IL)\bigr)$ where $P_1$ and $P_2$ are the depths of the convergence trees for the bisection processes implemented in Algorithms~\ref{algo:nu_star} and \ref{algo:e_star}, respectively, and only dependent on the precision parameter $\epsilon$.
This complexity is linear in the number of clusters $I$ and the number of job classes $L$, and log-linear in the capacity of each component $C_i$, resulting in a reasonably fast procedure of obtaining the estimated indices $\bm{u}(\bar{e}_0)$. Recall that the estimated indices $\bm{u}(\bar{e}_0)$ are pre-calculated, and the computational complexity for pre-computing $\bm{u}(\bar{e}_0)$ is different from that of implementing MPMP. 
As mentioned earlier in this section, for implementing MPMP in a server farm system, the computational complexity is at most linear to the number of physical components for each arrived job.
Note that fitting the value of $e_0$ provides a straightforward method but not the only way to approximate $e^*$. Based on Proposition~\ref{prop:asym_opt}, any $e=\lim_{h\to +\infty}e^*$ can lead to an asymptotically optimal policy $\varphi(e)$, satisfying \eqref{eqn:index_policy:1} and \eqref{eqn:index_policy:2} with substituted $\varphi(e)$ and $u^*_{\ell,i_j}(n,e)$ for $\varphi$ and $\upsilon^*_{\ell,j}(n)$, respectively.

\section{Asymptotic Optimality} \label{sec:asym_opt}
For given $e\in\mathbb{R}$, we say a policy $\phi\in\Phi$, applicable to the original problem described in \eqref{eqn:opt}, \eqref{eqn:constraint:action}, \eqref{eqn:constraint:capacity} and \eqref{eqn:constraint:boundary} with substituted $e$ for $e^*$, is \emph{asymptotically optimal} if 
\begin{equation}\label{eqn:asym_opt}
\lim\limits_{h\rightarrow\infty} \left\lvert \Gamma^{h,\phi}(e) - \max\limits_{\phi'\in\Phi}\Gamma^{h,\phi'}(e)\right\rvert = 0.
\end{equation}
If a policy is asymptotically optimal, it  approaches optimality for the server farm problem as $h\rightarrow +\infty$.
Recall that, when the parameter $e=e^*$ satisfies \eqref{eqn:e_star}, a policy optimal for the problem described in \eqref{eqn:opt}, \eqref{eqn:constraint:action}, \eqref{eqn:constraint:capacity} and \eqref{eqn:constraint:boundary} is also optimal for the server farm problem described in \eqref{eqn:obj}, \eqref{eqn:constraint:action}, \eqref{eqn:constraint:capacity} and \eqref{eqn:constraint:boundary} aiming to maximize the ratio of the long-run average job throughput to the long-run average power consumption.
\begin{definition}
For any given $e\in\mathbb{R}$, let $\varphi(e)$ represent the policy described in \eqref{eqn:index_policy:1} and \eqref{eqn:index_policy:2} with substituted $\varphi(e)$ and $u^*_{\ell,i_j}(n,e)$ for $\varphi$ and $\upsilon^*_{\ell,j}(n)$, respectively.
\end{definition}
\begin{proposition}\label{prop:asym_opt}
For the problem described in \eqref{eqn:opt}, \eqref{eqn:constraint:action}, \eqref{eqn:constraint:capacity} and \eqref{eqn:constraint:boundary} with given $e\in\mathbb{R}$ substituting for $e^*$, when job sizes are exponentially distributed, if \eqref{eqn:e_star_equal:condition} holds, then the policy $\varphi(e)$ is asymptotically optimal; that is, 
\eqref{eqn:asym_opt} holds with substituted $\varphi(e)$ for $\phi$.
\end{proposition}
The proposition is proved in Appendix~\ref{app:prop:asym_opt}. 
Together with Lemma~\ref{lemma:convergence_relaxed_opt}, when job sizes are exponentially distributed, if the components are energy-efficiently unimodal, and Condition~\ref{cond:heavy} holds or $L=1$, then the policy $\varphi(e)$ is asymptotically optimal.
The problem, described in \eqref{eqn:opt}, \eqref{eqn:constraint:action}, \eqref{eqn:constraint:capacity} and \eqref{eqn:constraint:boundary} with given $e\in\mathbb{R}$ substituting for $e^*$, aims to maximize the difference between the long-run average job throughput and the long-run average energy consumption rate weighted by the given $e$.
Since we assume quite general $e\in\mathbb{R}$ and $\mu_i(n)$ and $\varepsilon_i(n)$ for the component clusters, the average job throughput and energy consumption rate can be directly generalized as the average reward and cost of the process.
This is a popular objective and has been widely used in the literature, such as \cite{pedram2012energy,wu2017stochastic,fu2018restless}.
For such an objective, the simple policy $\varphi(e)$ is asymptotically optimal under the provided conditions.

Recall that, in general, asymptotically optimal policies do not necessarily exist for an RMABP or any similar system. 
Even if we perfectly model our problem as an RMABP or RMABP-like problem, there are no off-the-shelf steps that can lead to a scalable, asymptotically optimal policy.
As mentioned in earlier sections, the existence of scalable, near-optimal policies, usually referred to as index policies, is subject to (Whittle) indexability of all the underlying bandit processes. 
Asymptotic optimality of the proposed index policies, if exist, is further dependent on a non-trivial condition that requires the underlying stochastic process to converge to a fixed point, referred to as a global attractor, almost surely in the asymptotic regime. Neither the indexability nor the global attractor holds in general. 
In Section~\ref{sec:indexability}, we have discussed the (asymptotic) indexability and, to prove Proposition~\ref{prop:asym_opt}, it remains to analyze the global attractor, which fortunately exists for our server farm. 
We refer to Appendix~\ref{app:prop:asym_opt} for a rigorous proof of Proposition~\ref{prop:asym_opt}.

\begin{corollary}\label{coro:asym_opt}
For the problem defined by \eqref{eqn:obj}, \eqref{eqn:constraint:action}, \eqref{eqn:constraint:capacity} and \eqref{eqn:constraint:boundary}, when job sizes are exponentially distributed, if the computing components are energy-efficiently unimodal, and Condition~\ref{cond:heavy} holds or $L=1$, then the policy $\varphi(e_0)$, or equivalently the MPMP policy, is asymptotically optimal; that is, 
\eqref{eqn:asym_opt} holds with substituted $\varphi(e_0)$ for $\phi$.
\end{corollary}
The corollary is a straightforward result of Corollary~\ref{coro:e_star_equal}, Lemma~\ref{lemma:convergence_relaxed_opt} and Proposition~\ref{prop:asym_opt}.

{\bf Remark} 
Recall that the MPMP policy is scalable for the original server farm problem, for which the computational complexity has been discussed at the end of Section~\ref{subsec:ratio}.
Asymptotic optimality implies that such a scalable policy is approaching optimality as the system becomes large.
Unlike previous work in \cite{fu2016asymptotic, fu2018restless}, Proposition~\ref{prop:asym_opt} and Corollary~\ref{coro:asym_opt} apply to systems with a significantly less stringent relationship between power consumption, service rate, and traffic load.
Numerical results are provided in Section~\ref{sec:simulation} for holistic performance analysis.

We now further discuss the relationship between the sub-optimality of the index policy $\varphi(e)$ and the scaling parameter $h$ of the server farm system.
We rank the state-component (SC) pairs $(n,i)$, $n\in\mathscr{N}_i\backslash\{C_{i_j}\}$ and $i\in[I]$, according to the descending order of $\eta_{i,n}$ where $\eta_{i,n}$ satisfies $f^h_{i,n}(\eta_{i,n})=0$ with $f^h_{i,n}$ defined in \eqref{define:f}. Recall that the indices $\upsilon^*_{\ell,j}(n) = \bar{\upsilon}^*_{\ell,i_j}(n) = \lambda^0_{\ell}\eta_{i_j,n}/\hat{\lambda}^0_{i_j}$ for $\ell\in[L]$, $j\in[J]$, and $n\in\mathscr{N}_{i_j}\backslash\{C_{i_j}\}$.
Then, we place the remaining SC pairs $(n,i)$ with $n=C_i$ for all $i\in[I]$ afterwards.
To emphasize this ranking, we refer to  the $k$th SC pair as SC pair $k$, where $k = 1,2,\ldots, K$ for $K=\sum_{i\in[I]}|\mathscr{N}_i|$.
Let $Z^{\phi,h}_{k}(t)$ represent the proportion of processes $\{N^{\phi}_j(t),t\geq0\}$ ($j\in[J]$) that are in the $k$th SC pair at time $t$ under policy $\phi$; that is,
\begin{equation}\label{eqn:define_z}
Z^{\phi,h}_k(t) = \frac{1}{J}\left|\{j\in [J]\ |\ i_j=i_k,N^{\phi}_j(t) = n_k\}\right|,
\end{equation}
where $i_k\in[I]$ and $n_k\in\mathscr{N}_{i_k}$ are the cluster and state labels of SC pair $k$.
Recall that $J$ is defined in Section~\ref{sec:model} and is dependent on $h$. Define $\bm{Z}^{\phi,h}(t)\coloneqq (Z^{\phi,h}_k(t): k\in[K])$.
In this context, for any given $h\in\mathbb{N}_+$, the stochastic process $\{\bm{N}^{\phi}(t),t\geq 0\}$ can be translated to the process $\{\bm{Z}^{\phi,h}(t),t\geq 0\}$.

Consider a server farm that starts with no job (that is, $\bm{N}^{\phi}(0) = \bm{0}$) and, correspondingly, we have $\bm{Z}^{\phi,h}(0) = \bm{z}^0$ for this empty system.
Let $\mathscr{Z}\coloneqq [0,1]^K$ represent a probability simplex. 
\begin{proposition}\label{prop:exp_diminishing}
When job sizes are exponentially distributed, for given $e\in\mathbb{R}$, there exists $\bm{z}^{\varphi(e)}\in\mathscr{Z}$ such that, for any $\delta > 0$, there exist $s>0$ and $H>0$ satisfying, for all $h>H$, 
\begin{equation}\label{eqn:exp_diminishing}
\lim\limits_{T\rightarrow \infty} \frac{1}{T}\int_0^T \mathbb{P}\Bigl\{ \lVert \bm{Z}^{\varphi(e),h}(t)-\bm{z}^{\varphi(e)}\rVert>\delta\Bigr\} dt\leq e^{-sh},
\end{equation}
where $\bm{Z}^{\varphi(e),h}(0)=\bm{z}^0$.
\end{proposition}
The proposition is proved in Appendix~\ref{app:prop:exp_diminishing}. 
Proposition~\ref{prop:exp_diminishing} indicates that, under the policy $\varphi(e)$, the underlying stochastic process $\{\bm{Z}^{\varphi(e),h}(t),t\geq 0\}$ converges to a global attractor $\bm{z}^{\varphi(e)}$ almost surely as $h\rightarrow \infty$, and, more importantly, the deviation between $\bm{Z}^{\varphi(e),h}(t)$ and $\bm{z}^{\varphi(e)}$ is diminishing exponentially in $h$.
Define $\bm{r}(e)\coloneqq (\mu_{i_k}(n_k)-e\varepsilon_{i_k}(n_k): k\in[K])$. Then, for given $h\in\mathbb{N}_+$ and policy $\phi\in\Phi$, the normalized long-run average reward 
\begin{equation}\label{eqn:gamma_to_z}
\Gamma^{\phi,h}(e) = \frac{J}{h} \bm{r}(e)\cdot\lim\limits_{t\rightarrow\infty} \mathbb{E}[\bm{Z}^{\phi,h}(t)]= \sum\limits_{i\in[I]}M_i^0\bm{r}(e)\cdot\lim\limits_{t\rightarrow\infty} \mathbb{E}[\bm{Z}^{\phi,h}(t)],
\end{equation}
where $\bm{x}\cdot\bm{y}$ represents the inner product of vectors $\bm{x}$ and $\bm{y}$.
From Proposition~\ref{prop:exp_diminishing}, when job sizes are exponentially distributed, $\lim\nolimits_{h\rightarrow\infty} \Gamma^{\varphi(e),h}(e) = \sum\nolimits_{i\in[I]}M_i^0 \bm{r}(e)\cdot\bm{z}^{\varphi(e)}$.
If $\bm{z}^{\varphi(e)}$ coincides with an optimal point of the relaxed problem, described in~\eqref{eqn:opt}, \eqref{eqn:constraint:relax:action}, \eqref{eqn:constraint:relax:capacity} and \eqref{eqn:constraint:relax:boundary} with $e$ substituting for $e^*$, in the asymptotic regime, then the policy $\varphi(e)$ is asymptotically optimal. 
Together with asymptotic indexability discussed in Section~\ref{sec:indexability}, Proposition~\ref{prop:exp_diminishing} can lead to Proposition~\ref{prop:asym_opt}.
Apart from asymptotic optimality, Proposition~\ref{prop:exp_diminishing} also implies that, for a system with large $h$, if $\varphi(e)$ is asymptotically optimal, the performance deviation between $\varphi(e)$ and optimality in the asymptotic regime diminishes exponentially as $h\rightarrow \infty$.
This conclusion extends \cite[Proposition 2]{fu2020energy} to the more generalized server farm model discussed in this paper, and, as a straightforward result of Proposition~\ref{prop:exp_diminishing}, when job sizes are exponentially distributed, \eqref{eqn:exp_diminishing} also applies to the MPMP policy by setting $e=e_0$.

\section{Numerical Results}\label{sec:simulation}

In the case without assuming the energy-efficient unimodality or Condition~\ref{cond:heavy}, we numerically demonstrate the effectiveness of MPMP by comparing it with two baseline policies in different circumstances. 
In the results demonstrated in this section, the $95\%$ confidence intervals based on the Student t-distribution are maintained within $3\%$ of the observed mean. In particular, in Section~\ref{subsec:effectiveness}, we consider exponentially distributed job sizes, and, in Section~\ref{subsec:sensitivity}, the energy efficiency of MPMP is examined with other job-size distributions.

\subsection{Effectiveness of MPMP}
\label{subsec:effectiveness}

Here, we discuss the effectiveness of MPMP by comparing it with baseline policies in two scenarios: systems with randomly generated parameters and a realistic server farm with real-world settings.

\subsubsection{Scenario I}\label{subsubsec:case1}
Consider a server farm with ten clusters ($I=10$) and four job classes ($L=4$), where the peak service rate for each cluster (that is, $\mu_i(C_i)$ for $i\in[I]$) is uniformly randomly generated from $[10,15]$ and the capacities $C_i$ are set $5$ for all $i\in[I]$. For cluster $i\in[I]$, we uniformly randomly generate the energy efficiency of its fully-occupied component, $\mu_i(C_i)/\varepsilon_i(C_i)$, from $[0.5,1]$, and set the power consumption of its idle component $\varepsilon_i(0)$ to be $0.3\varepsilon_i(C_i)(0.9-0.1 i)$. 
For other states $n\in\mathscr{N}_i\backslash\{C_i\}$ of cluster $i$, the service and energy consumption rates are obtained by setting $\mu_i(n)=\mu_i(n+1)\frac{n}{n+1}$  and $\varepsilon_i(n) = (\varepsilon_i(n+1)-\varepsilon_i(0)) \sqrt{\frac{n}{n+1}}+\varepsilon_i(0)$, respectively. In this case, a higher service rate indicates higher power consumption and higher energy efficiency, which follows realistic situations in \cite{giri2010increasing,kaur2015energy}.
We take the scaling parameter $h= 10$ and the number of components in each cluster $M^0_i=1$ for all $i\in[I]$.

For job class $\ell\in[L]$, we uniformly randomly generate an integer $\kappa_{\ell}\in[I]$ as the number of clusters involving available components for $\ell$-jobs, and then randomly select $\kappa_{\ell}$ clusters: all the components within the selected clusters are available components for $\ell$-jobs and join the set $\mathcal{J}_{\ell}$. 
Define the \emph{normalized offered traffic} of job class $\ell\in[L]$ as $\rho_{\ell}\coloneqq \lambda_{\ell}/\sum_{j\in \mathcal{J}_{\ell}}\mu_{i_{j}}(C_{i_j})$ and the \emph{relative difference} of policy $\phi_1$ to $\phi_2$ ($\phi_1,\phi_2\in\Phi$) with respect to energy efficiency as
\begin{equation}
\frac{\mathfrak{L}^{\phi_1}/\mathfrak{E}^{\phi_1}-\mathfrak{L}^{\phi_2}/\mathfrak{E}^{\phi_2}}{\mathfrak{L}^{\phi_2}/\mathfrak{E}^{\phi_2}}.
\end{equation}
In Scenario I, we set $\rho_{\ell} = \rho$ for all $\ell\in[L]$ and then compute the arrival rates of different job classes.

\begin{figure}[t]
\centering
\subfigure[]{\includegraphics[width=0.24\linewidth]{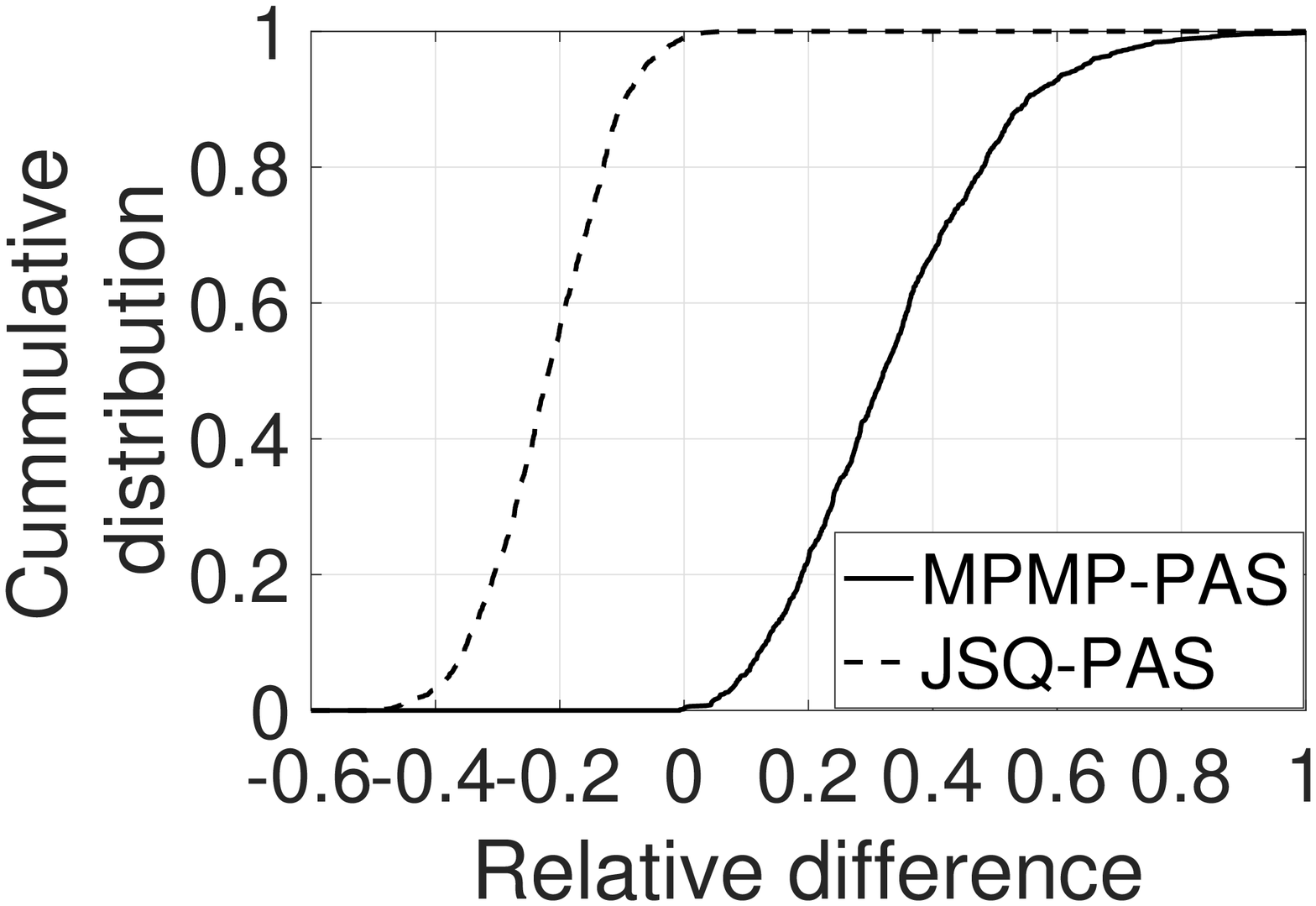}\label{fig:fig1a}}
\subfigure[]{\includegraphics[width=0.24\linewidth]{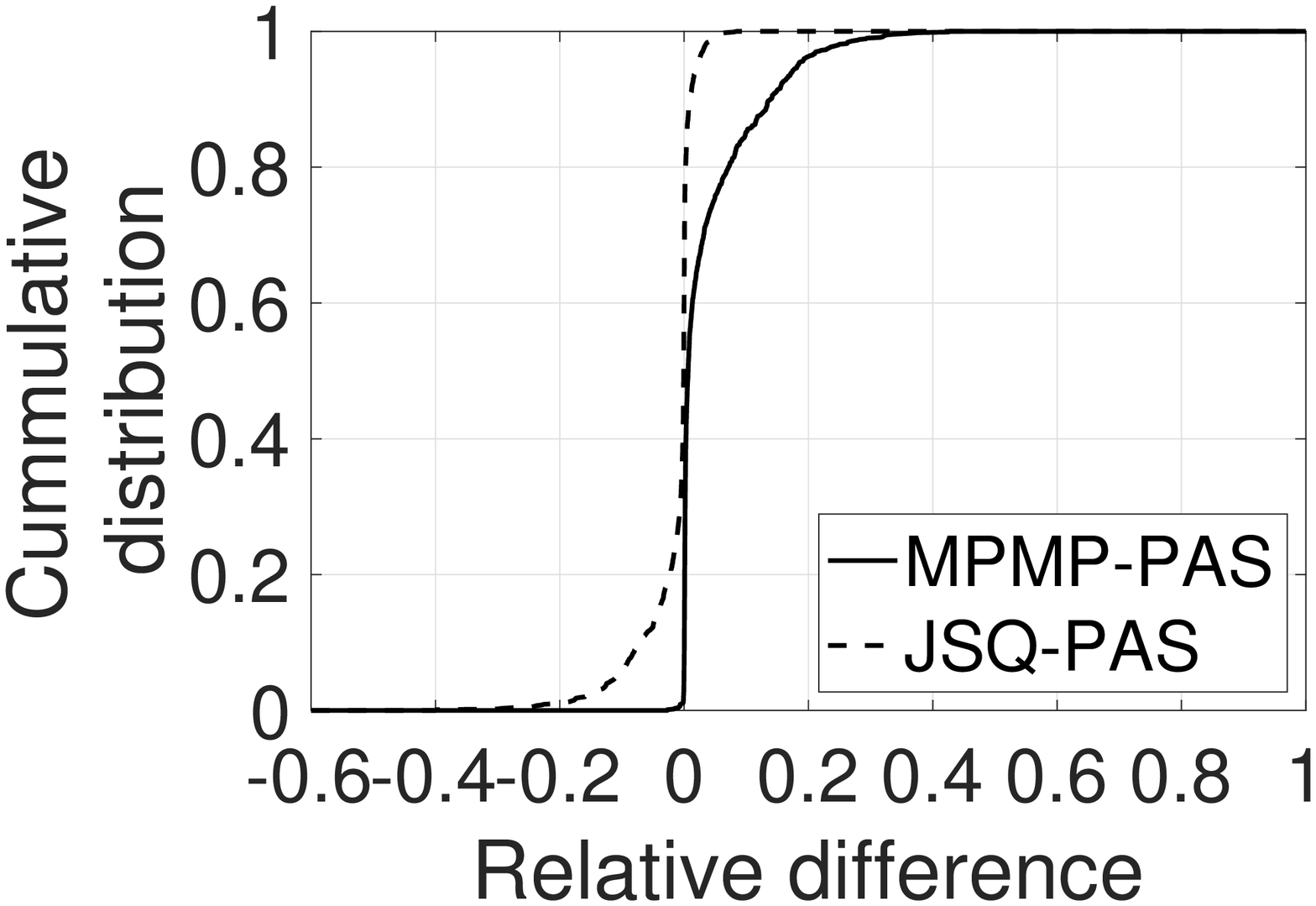}\label{fig:fig1c}}
\caption{Cumulative distribution of the relative difference of MPMP and JSQ to PAS with respect to energy efficiency: (a) $\rho=0.2$; and (b) $\rho=0.5$.}\label{fig:fig1}
\end{figure}

In Figure~\ref{fig:fig1}, we compare the energy efficiency of MPMP to two baseline policies: Join-the-Shortest-Queue (JSQ) and Priorities accounting for Available Servers (PAS). JSQ is a classic job-assignment policy used to balance the distributed workload to all the servers and was proved to maximize the throughput (number of processed jobs) within a given time horizon \cite{winston1977optimality}. PAS was proposed in \cite{fu2020energy} and proved to be asymptotically optimal with respect to energy efficiency in a server farm with only two power modes for each server. As mentioned in Section~\ref{sec:algo},
when all the physical components have only two power modes (the case discussed in \cite{fu2020energy}), MPMP reduces to PAS. 
While, in this paper, we are interested in a more general case with multiple power modes, and PAS fails to capture the features and advantages of intermediate power modes that are associated with neither peak nor idle energy consumption/service rates as illustrated in the following.
In Figure~\ref{fig:fig1}, we demonstrate the cumulative distribution of relative differences for a thousand simulation runs with randomly generated system parameters. 
MPMP outperforms PAS in all cases and achieves $20\%$ higher energy efficiency than that of PAS for around $80\%$ of the simulation runs when $\rho = 0.2$.
Observing the results for $\rho=0.2$ and $0.5$, the advantages of MPMP over PAS are decreasing as the offered traffic becomes heavier. Because heavier offered traffic will impose the underlying stochastic process to stay longer in states with a larger number of fully-occupied physical components, which are running at peak service and energy consumption rates, where the effects of intermediate power modes are mitigated.

\begin{figure}[t]
\centering
\begin{minipage}[]{0.24\textwidth}
\includegraphics[width=\linewidth]{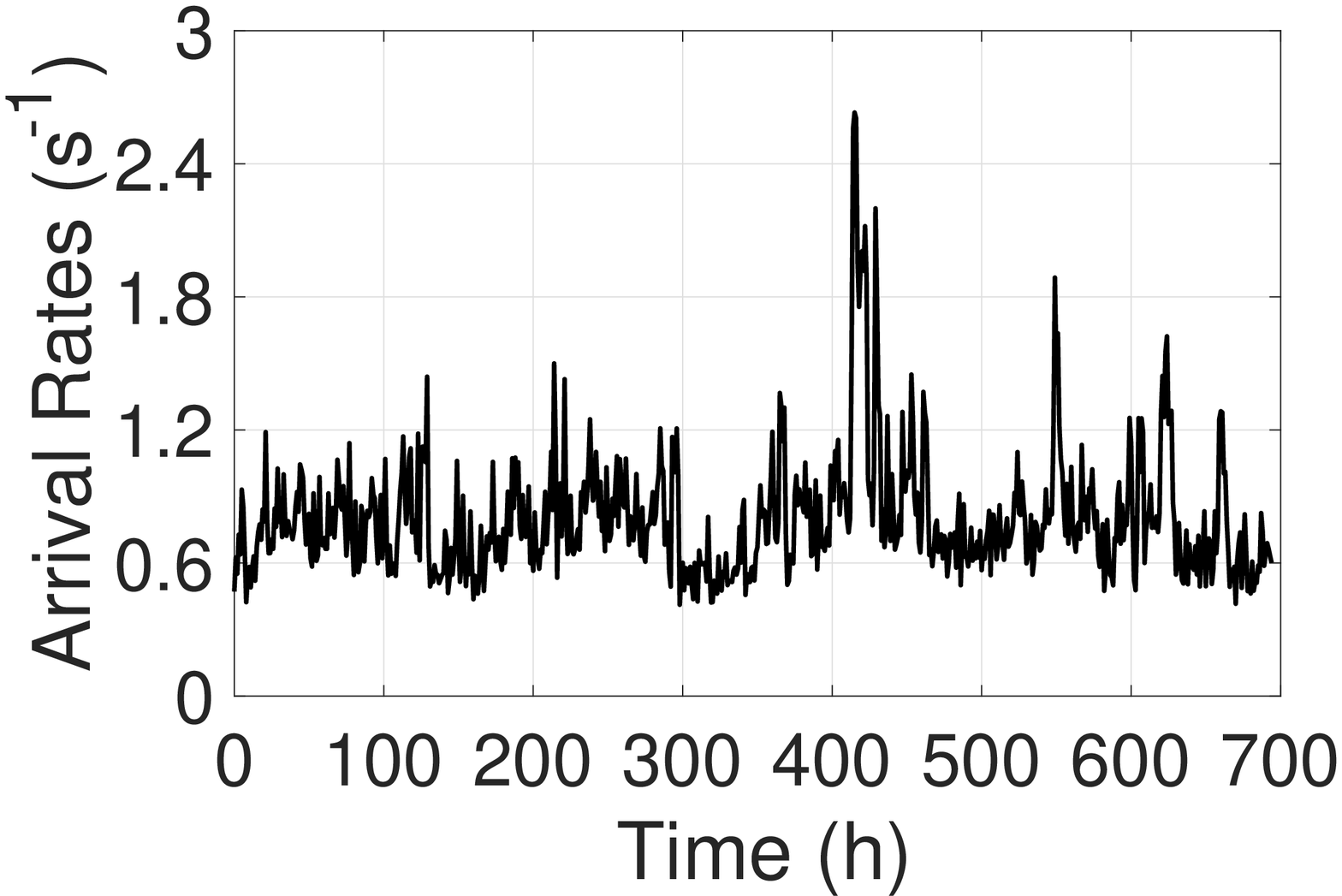}
\caption{Arrival rates of Google trace-logs.\label{fig:google-arrival}}
\end{minipage}
\begin{minipage}[]{0.24\textwidth}
\includegraphics[width=\linewidth]{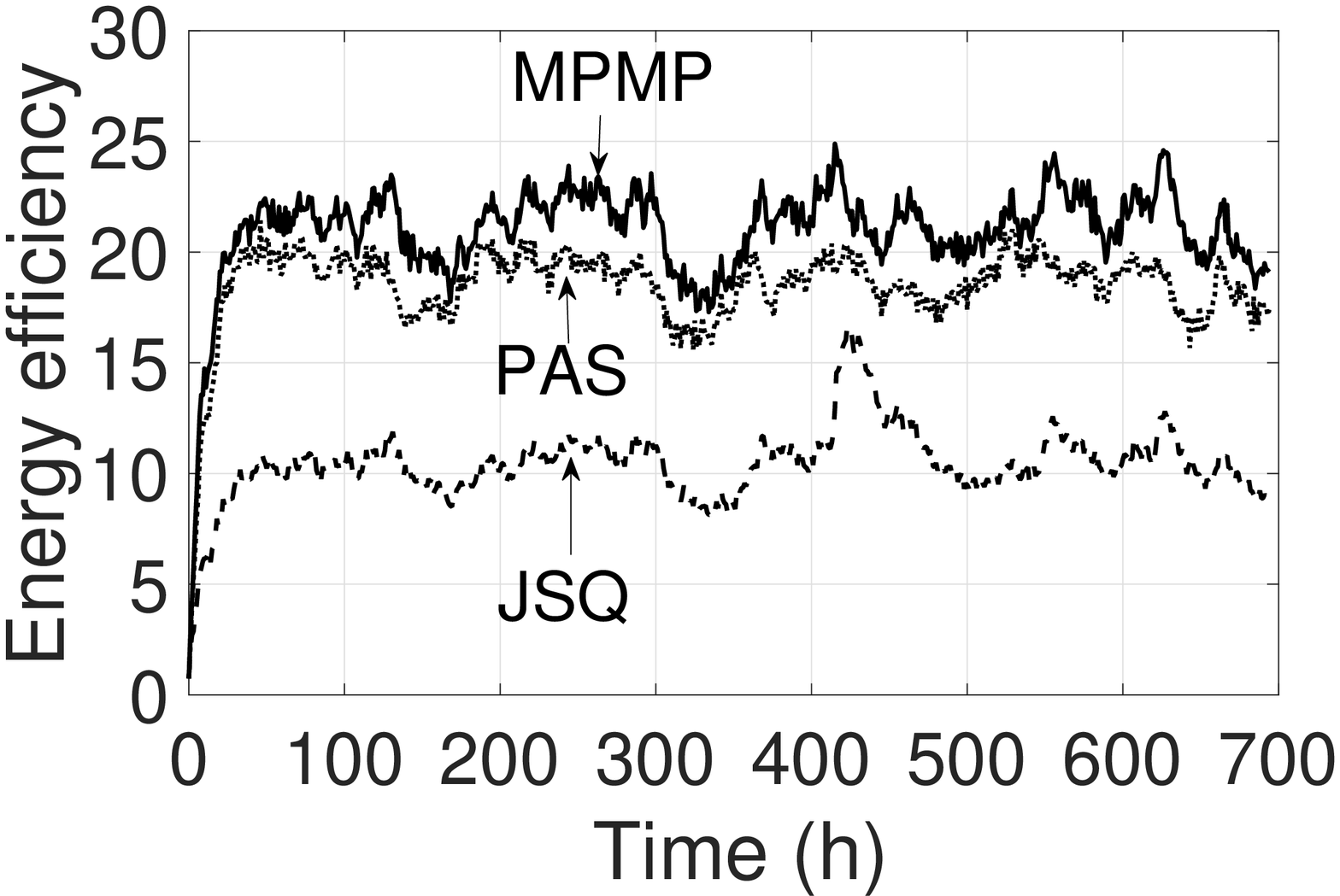}
\caption{Energy efficiency of MPMP, PAS and JSQ.\label{fig:google-TE}}
\end{minipage}
\centering
\begin{minipage}[]{0.24\textwidth}
\includegraphics[width=\linewidth]{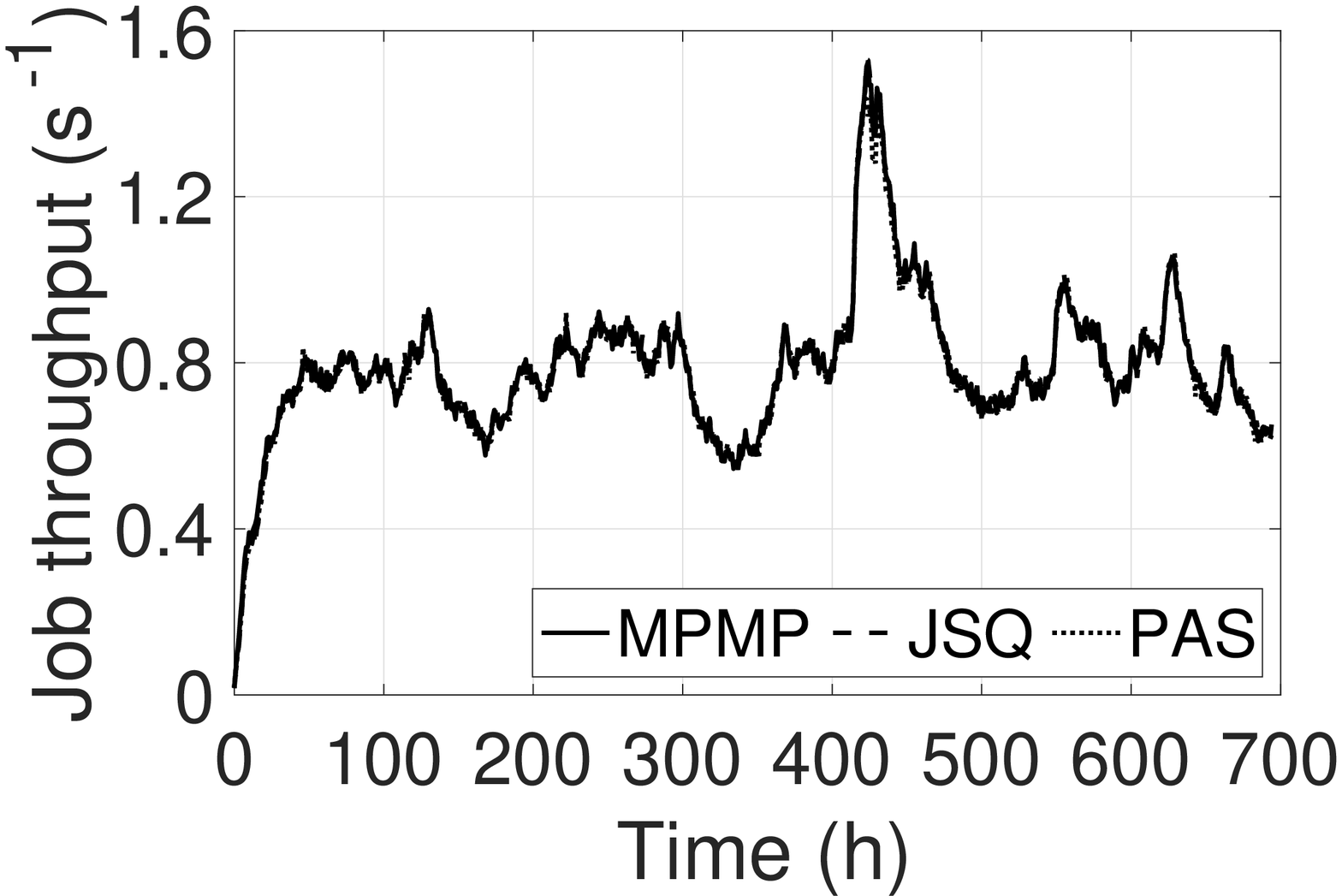}
\caption{Job throughput of MPMP, PAS and JSQ.\label{fig:google-throughput}}
\end{minipage}
\begin{minipage}[]{0.24\textwidth}
\includegraphics[width=\linewidth]{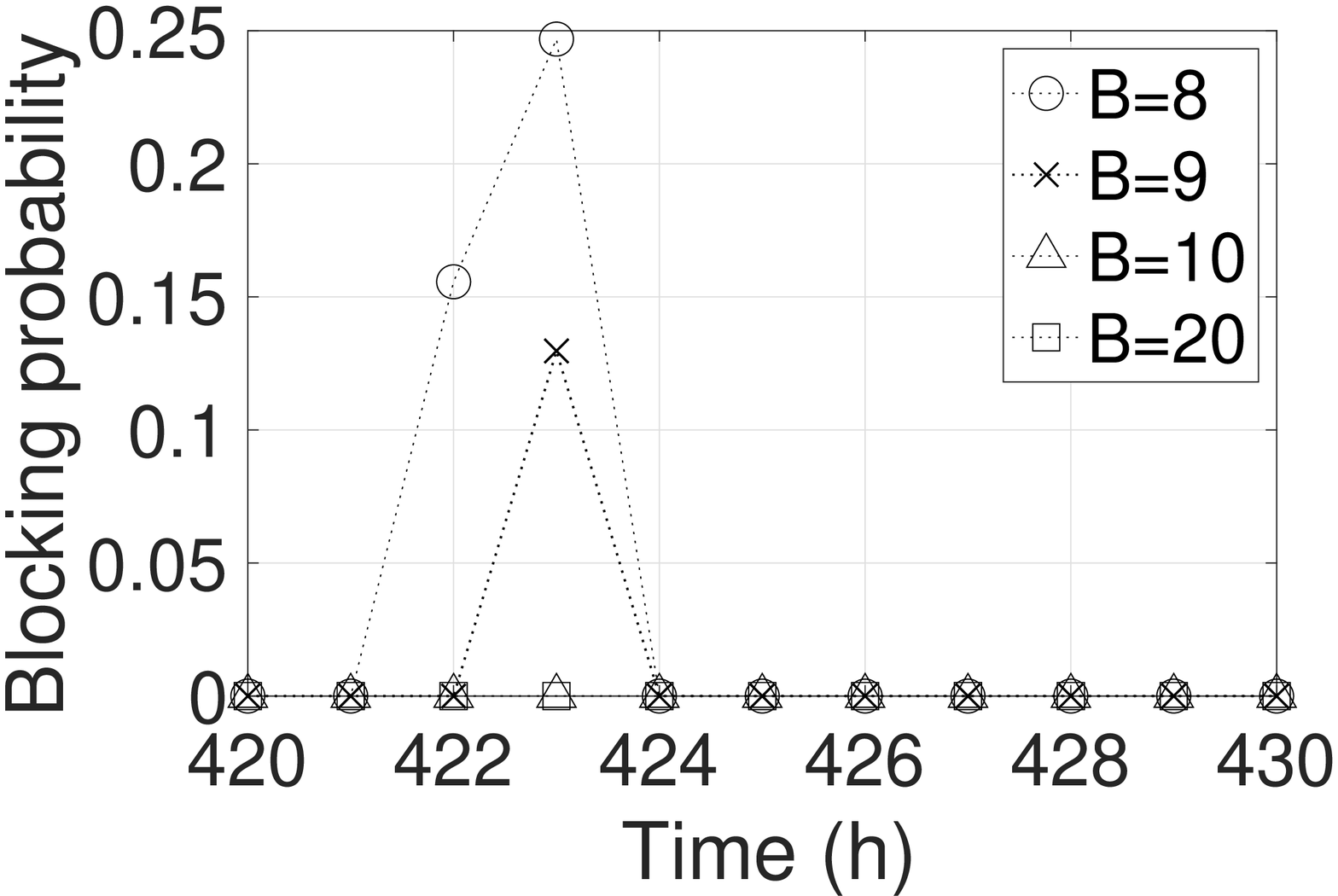}
\caption{Blocking probabilities of MPMP, PAS and JSQ.\label{fig:google-blocking}}
\end{minipage}
\end{figure}

\begin{figure}[t]
\centering
\begin{minipage}[]{0.24\textwidth}
\includegraphics[width=\linewidth]{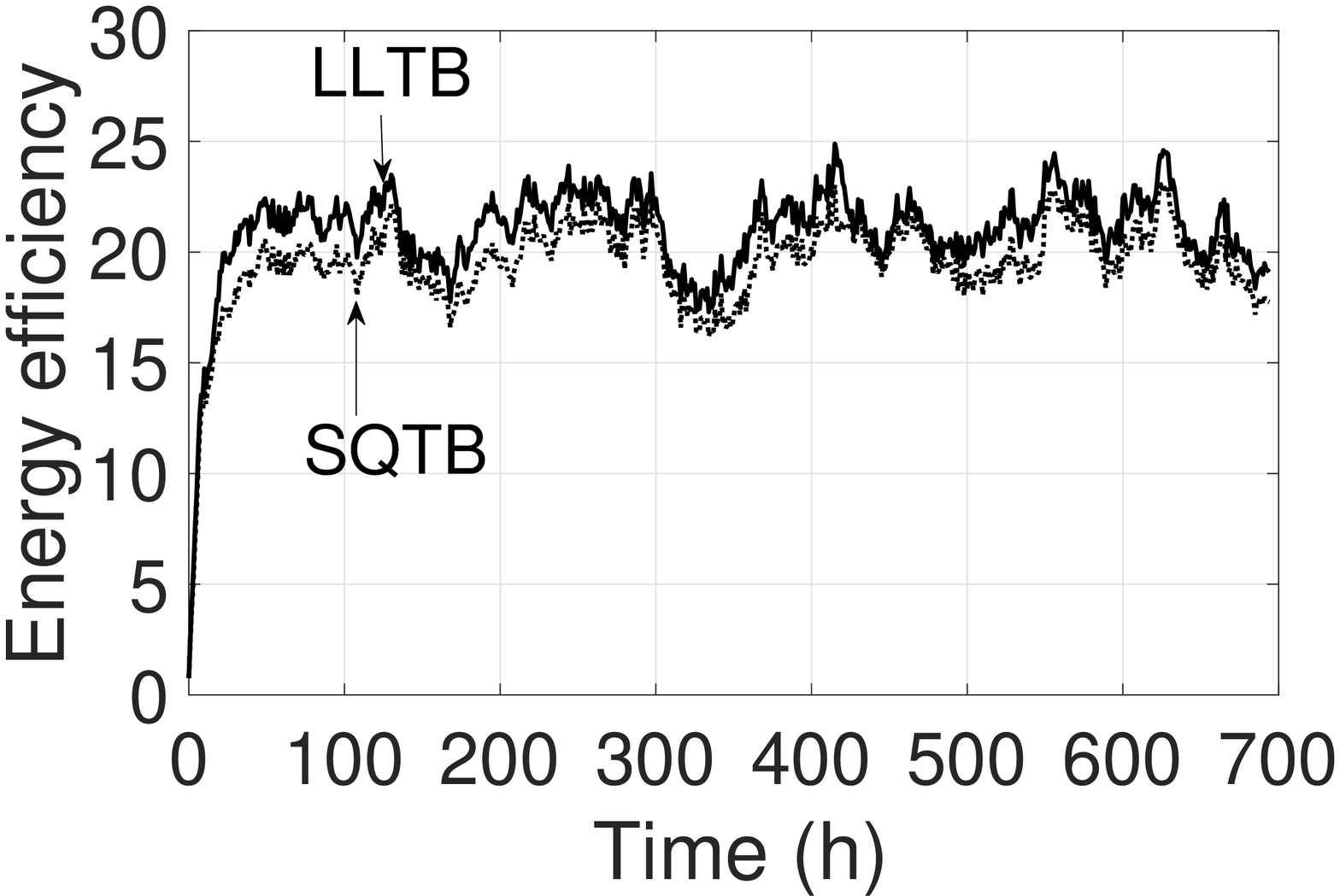}
\caption{Energy efficiency of MPMP with different tie-breaking rules.\label{fig:google-tie-breaking-TE}}
\vspace{2em}
\end{minipage}
\begin{minipage}[]{0.24\textwidth}
\includegraphics[width=\linewidth]{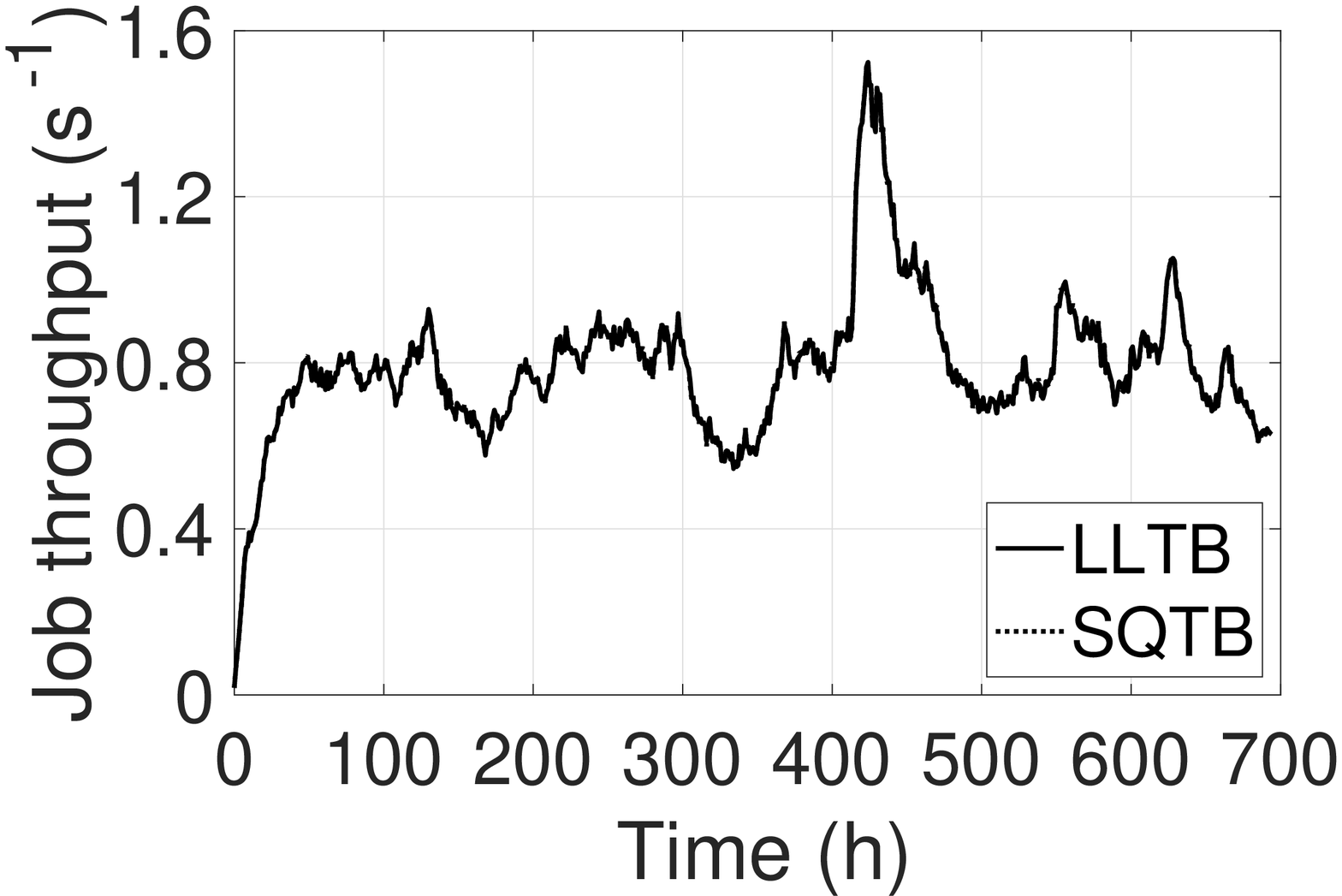}
\caption{Job throughput of MPMP with different tie-breaking rules.\label{fig:google-tie-breaking-T}}
\vspace{2em}
\end{minipage}
\centering
\begin{minipage}[]{0.24\textwidth}
\includegraphics[width=\linewidth]{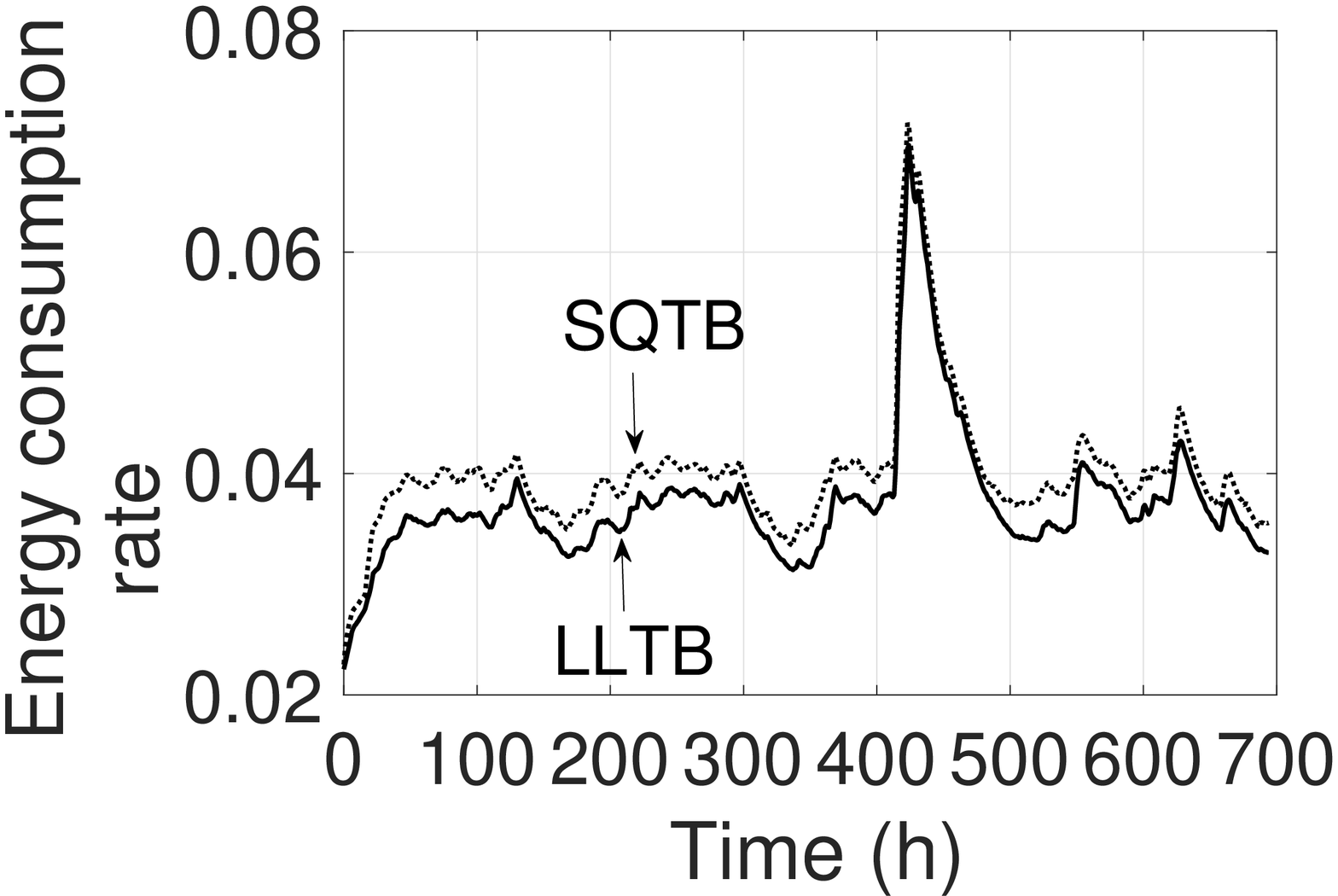}
\caption{Power consumption of MPMP with different tie-breaking rules.\label{fig:google-tie-breaking-E}}
\vspace{1.2em}
\end{minipage}
\begin{minipage}[]{0.24\textwidth}
\includegraphics[width=\linewidth]{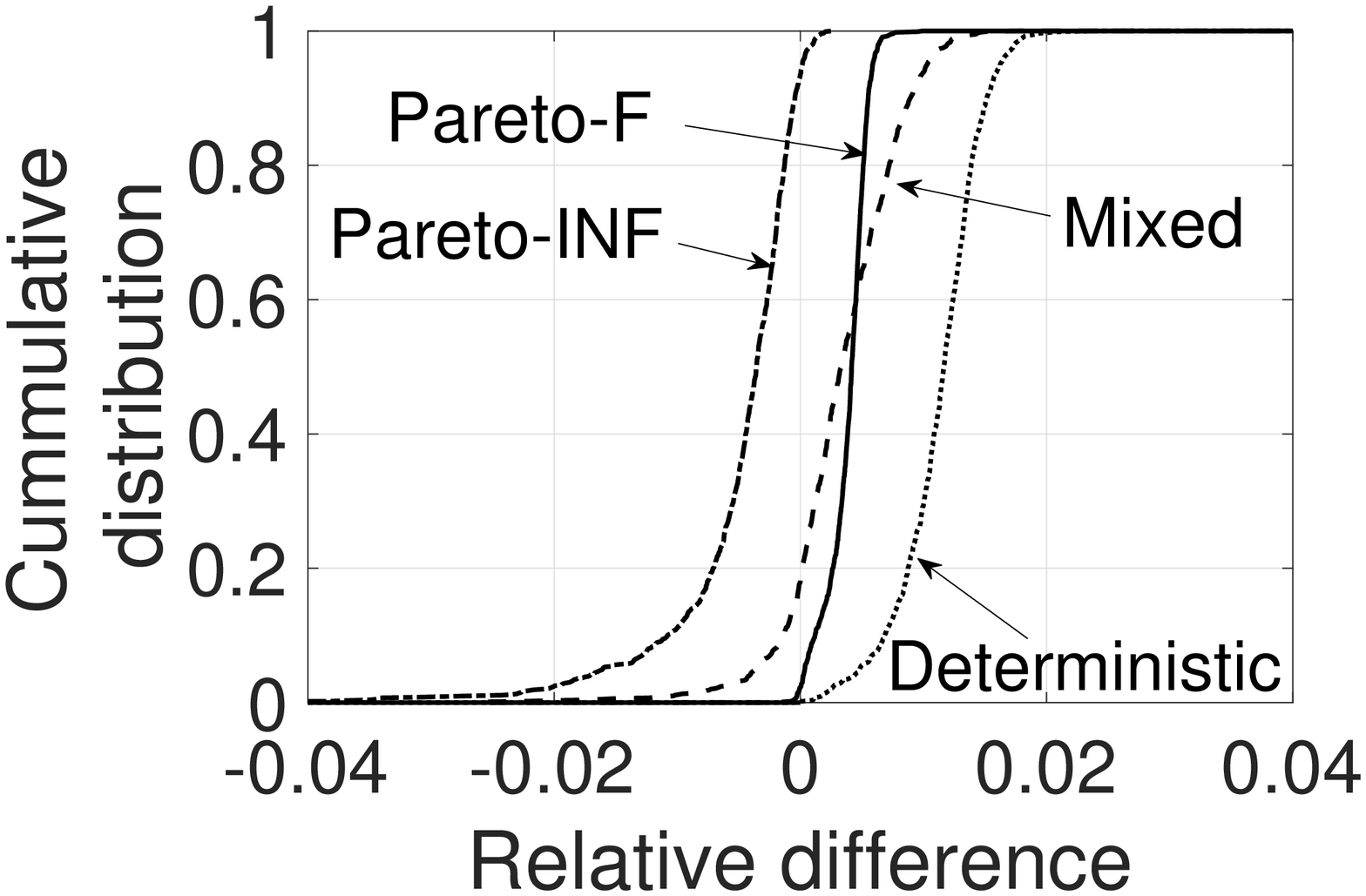}
\caption{Relative difference of
energy efficiency with job size distribution $\mathcal{D}$ from the one
with exponentially distributed job sizes.\label{fig:sensitivity}}
\end{minipage}
\end{figure}

\subsubsection{Scenario II}\label{subsubsec:case2}
Another scenario is considered with Google cluster traces of job arrivals in 2011 \cite{clusterdata:Wilkes2011,clusterdata:Reiss2011}, where there are 12.5 thousand physical components with arriving jobs classified into four groups ($L=4$). The total job arrival rates, estimated as the number of arrived jobs per second, averaged in each hour, are plotted in Figure~\ref{fig:google-arrival}. 
Similar to Scenario I in Section~\ref{subsubsec:case1}, we divide the system into ten clusters ($I=10$), each of which includes 1.25 thousand components; that is, setting $M^0_i=1$ for all $i\in[I]$ with scaling parameter $h=1250$. In this scenario, we no longer assume Poisson arrival processes with given arrival rates $\lambda_{\ell}$ ($\ell\in[L]$) but replace them with the job arrivals of the Google cluster traces. In particular, during busy hours, the total arrival rates of jobs are higher than the total service rate of the server farm. For the simulations discussed in Figures~\ref{fig:google-TE} and \ref{fig:google-throughput}, the capacities of physical components are set to ten ($C_i=10$) for all $i\in[I]$. 
Other detailed settings of the physical components are provided in Appendix~\ref{app:google-settings}.

In Figure~\ref{fig:google-TE}, we present the energy efficiencies of MPMP, PAS, and JSQ, averaged in each hour, for scenario II, where MPMP significantly outperforms PAS and JSQ. In particular, the total energy efficiency of MPMP is around $13\%$ higher than that of PAS.
In Figure~\ref{fig:google-throughput}, we plot the job throughput, normalized by the scaling parameter $h$, of MPMP, PAS, and JSQ for the same simulation discussed in Figure~\ref{fig:google-TE}. Given the clear advantages of MPMP against PAS and JSQ with respect to energy efficiency, they still achieve almost the same job throughputs. It is because the blocking probabilities of the three policies are likely to be negligible, although the capacity for each component is relatively small.

In this context, we further consider the blocking probabilities and present in Figure~\ref{fig:google-blocking} with different component capacities $C_i = B = 8,9,10$ and $20$ ($i\in[I]$). In Figure~\ref{fig:google-blocking}, when $B \leq 9$, we can observe zero blocking probabilities for MPMP except for short time periods of busy hours, when the total arrival rate is higher than the total service rate of the server farm.
When the component capacities are no less than $10$, the blocking probabilities of MPMP and JSQ have become zero; and for $B=20$, the blocking probability of PAS also diminishes. 
As the server farm is sufficiently large with scaling parameter $h=1250$, although the capacity for each physical component is relatively small and there are busy hour time periods with heavy traffic, the blocking probability remains zero for the entire duration of the simulation. 

Recall that MPMP is an index policy that prioritizes physical components with the highest indices.
As mentioned in Section~\ref{sec:algo}, all the theoretical results in this paper apply to the index policy $\varphi(e)$ with arbitrary tie-breaking rules. 
Beyond the theoretical conclusions, in certain cases, different tie-breaking rules may lead to different performance. 
Consider a tie-breaking rule that selects the component with the lowest label, and refer to it as the \emph{Lowest-Label Tie-Breaking} (LLTB).
All the simulations presented in Figures~\ref{fig:fig1}-\ref{fig:google-blocking} are based on LLTB.
In Figure~\ref{fig:google-tie-breaking-TE}, we explore the effects of different tie-breaking rules. 
The settings for the simulations presented in Figure~\ref{fig:google-tie-breaking-TE} are the same as those for Figure~\ref{fig:google-TE}.
Apart from LLTB, we consider another tie-breaking rule that always selects the component with the least number of holding jobs in the tie-breaking cases and refer to it as the \emph{Shortest-Queue Tie-Breaking} (SQTB). In Figure~\ref{fig:google-tie-breaking-TE}, we demonstrate the energy efficiency of MPMP with LLTB and SQTB and observe that LLTB achieves slightly higher energy efficiency than that of SQTB. The total energy efficiency for LLTB is around $5\%$ higher than that of SQTB. 
We also present the job throughput and power consumption of MPMP in Figures~\ref{fig:google-tie-breaking-T} and \ref{fig:google-tie-breaking-E}, respectively. 
In Figure~\ref{fig:google-tie-breaking-T}, SQTB achieves similar job throughput as LLTB;
while, in Figure~\ref{fig:google-tie-breaking-E}, SQTB achieves slightly higher power consumption than LLTB, which implies lower energy efficiency.
Recall that, based on Proposition~\ref{prop:exp_diminishing}, when job sizes are exponentially distributed, MPMP with different tie-breaking rules will approach the same performance as $h$ increases and, for a large system, the performance deviation between different tie-breaking rules diminishes exponentially in $h$. 
The $5\%$ difference in Figure~\ref{fig:google-tie-breaking-TE} between SQTB and LLTB is marginal, providing that the $95\%$ confidence intervals based on the Student t-distribution are $3\%$ of the observed mean.
This paper focuses on scalable and asymptotically optimal policies in large-scale server farms.
A thorough discussion on the effects of different tie-breaking rules in relatively small and practical systems is a fundamentally interesting topic on its own, which is beyond the scope of this paper.

\subsection{Sensitivity}\label{subsec:sensitivity}

The theoretical results presented in this paper and the numerical results presented in Section~\ref{subsubsec:case1} are based on exponential job-size distributions. From past studies \cite{crovella1997self,harchol2013performance}, real-world online applications exhibit job sizes that have heavy-tailed distributions. 
Here, we numerically demonstrate the performance of MPMP considering a heavy-tailed job-size (Pareto) distribution, as well as a mixed version with different job-size distributions for different job classes. 
In particular,
we consider simulations involving three non-exponential job size distributions with unit mean: deterministic, Pareto with shape parameter 2.001 and Pareto with shape parameter 1.98.
We refer to the Pareto distributions with shape parameters 2.001 and 1.98 as Pareto-F and Pareto-INF for short, as they have finite and infinite variances, respectively.
The settings for the simulations discussed in this subsection are the same as those for Figure~\ref{fig:fig1}  except for the job-size distributions.
Apart from the above-mentioned distributions, define a \emph{mixed} case where different job classes have different job-size distributions.

Define $\mathfrak{T}^{\text{MPMP},\mathcal{D}}$ and $\mathfrak{E}^{\text{MPMP},\mathcal{D}}$ as the long-run average job throughput and power consumption under MPMP with job-size distribution $\mathcal{D}$, respectively. 
In Figure~\ref{fig:sensitivity}, we present the cumulative distribution of the relative difference of energy efficiency with job size distribution $\mathcal{D}$ from the one
with exponentially distributed job sizes; that is, the cumulative distribution of
\begin{equation}
\frac{\mathfrak{T}^{\text{MPMP},\mathcal{D}}/\mathfrak{E}^{\text{MPMP},\mathcal{D}} - \mathfrak{T}^{\text{MPMP},exponential}/\mathfrak{E}^{\text{MPMP},exponential}}{\mathfrak{T}^{\text{MPMP},exponential}/\mathfrak{E}^{\text{MPMP},exponential}},
\end{equation}
with $\mathcal{D}=$ deterministic, Pareto-F, Pareto-INF and mixed. For the case of mixed, we set the job-size distributions for the job classes 1-4 as deterministic, exponential, Pareto-F, and Pareto-INF, respectively. 
In Figure~\ref{fig:sensitivity}, we observe that the relative differences of energy efficiency for all the tested simulation runs are within $\pm 3\%$, indicating similar energy efficiencies for tested $\mathcal{D}$ and  the exponential case.
It follows that the energy efficiency of MPMP is not very sensitive to the tested distributions, including the heavy-tailed Pareto-INF.

\section{Conclusions}\label{sec:conclusion}

Relying on the Whittle relaxation technique, we have proposed the MPMP policy that always prioritizes physical components with the highest indices satisfying \eqref{eqn:asym_opt:index}.
The indices are computable through Algorithms~\ref{algo:nu_star} and \ref{algo:e_star} within linear time in the number of clusters $I$ and the number of job classes $L$. It is log-linear in the capacity of each component, and the implementation of the MPMP policy is at most linear in the number of physical components.
It follows that MPMP is scalable as we have demonstrated its applicability to  a large-scale server farm with tens or hundreds of thousands of abstracted servers.

When job sizes are exponentially distributed and all the components are energy-efficiently unimodal, following the ideas of \cite{weber1990index,fu2018restless}, we have proved that, if $L=1$ or Condition~\ref{cond:heavy} holds, MPMP approaches optimality as the scaling parameter $h$ tends to infinity; that is, MPMP approaches optimality as the numbers of components in clusters and the arrival rates of jobs increase proportionately to infinity. The asymptotic optimality is applicable to systems with large and rapidly increasing numbers of storage or/and computing components, where conventional optimization techniques do not apply.  
We have provided an analysis of the entire system, including discussions on the indexability and the global attractor for proving asymptotic optimality in the continuous-time case. 
For a large system, we have proved that the performance deviation between MPMP and an optimal solution in the asymptotic regime diminishes exponentially in the scaling parameter $h$. That is, MPMP will become already close to optimality in a relatively small system.

For the non-asymptotic regime without assuming energy-efficient unimodality and the heavy traffic condition, we have numerically demonstrated the effectiveness of MPMP by comparing it to JSQ and PAS. When the job throughputs are compatible, MPMP has shown substantial advantages against the baseline policies with respect to energy efficiency. We have further investigated the performance of MPMP considering different job-size distributions and demonstrated that MPMP is robust in all the cases we tested.

\appendices

\section{Proof of Proposition~\ref{prop:relaxed_opt_precise}}\label{app:prop:relaxed_opt_precise}

For $\ell\in[L]$, $j\in\mathcal{J}_{\ell}$, let $\omega_{\ell,j}>-\nu_{\ell}$ all the time, such that an optimal solution of sub-problem \eqref{eqn:subproblem:normal}
always has $\alpha^{\phi}_{\ell,j}(C_{i_j})=0 $: constraints \eqref{eqn:constraint:relax:capacity} are satisfied.
In this context, we replace $\bm{r}^{\phi}_{j}(\bm{\nu},\bm{\omega})$ with $\bm{r}^{\phi}_j(\bm{\nu})\coloneqq \bm{r}^{\phi,e^*}_j(\bm{\nu})$ (defined in \eqref{eqn:define_r}) for the remainder of this appendix.

Consider the underlying stochastic process of the problem defined in \eqref{eqn:subproblem:normal} for component $j\in[J]$:
it is a birth-and-death Markov process $\{N^{\phi}_j(t), t\geq 0\}$, for which $N^{\phi}_j(t)$ represents the number of jobs being served by this component at time $t$ and $r^{\phi}_{j,n}(\bm{\nu})$ ($n\in\mathscr{N}_{i_j}$) is the reward rate in state $n$. We refer to this process $\{N^{\phi}_j(t),t\geq 0\}$ as the \emph{sub-process} associated with component $j$ or sub-process $j$.

Let $V^{\phi,\bm{\nu}}_j(n)$ and $T^{\phi,\bm{\nu}}_j(n)$ represent the expected accumulated reward and time, respectively, of sub-process $j$ starting from state $n$ and ending in state $0$ under policy $\phi$, where the reward rate in state $n$ is $r^{\phi}_{j,n}(\bm{\nu})$ and $V^{\phi,\bm{\nu}}_j(0)=T^{\phi,\bm{\nu}}_j(n)\equiv 0$. 
For given $g\in\mathbb{R}$, define $V^{\phi,g,\bm{\nu}}_j(n) = V^{\phi,\bm{\nu}}_j(n)-gT^{\phi,\bm{\nu}}_j(n)$ and $V^{g,\bm{\nu}}_j(n) = \max_{\phi\in\tilde{\Phi}_1}V^{\phi,g,\bm{\nu}}_j(n)$.

\begin{lemma}\label{lemma:relaxed_opt}
When job sizes are exponentially distributed, there exists a $g\in\mathbb{R}$ such that, for any policy $\phi^*\in\tilde{\Phi}_1$ that is optimal for the maximization problem in the Lagrangian dual function $L(\bm{\nu},\bm{\omega},\pmb{\gamma})$ of the relaxed problem,
\begin{equation}\label{eqn:relaxed_opt}
\alpha^{\phi^*}_{\ell,j}(n)=\left\{
\begin{cases}
1, & \text{if } \nu_{\ell} < \lambda_{\ell}\left(V^{g,\bm{\nu}}_j(n+1)-V^{g,\bm{\nu}}_j(n)\right),\\
a, & \text{if }\nu_{\ell} = \lambda_{\ell}\left(V^{g,\bm{\nu}}_j(n+1)-V^{g,\bm{\nu}}_j(n)\right),\\
0, & \text{if }\nu_{\ell} > \lambda_{\ell}\left(V^{g,\bm{\nu}}_j(n+1)-V^{g,\bm{\nu}}_j(n)\right),
\end{cases}\right.
\end{equation}
for recurrent states $n\in \mathscr{N}_{i_j}\backslash\{C_{i_j},0\}$, $j\in\mathcal{J}_{\ell}$, $\ell\in[L]$, where $a$ can be any real number in $[0,1]$.
\end{lemma}

\plainproof
As described in \eqref{eqn:decompose}, the maximization problem in the Lagrangian dual function $L(\bm{\nu},\bm{\omega},\pmb{\gamma})$ of the relaxed problem consists of $J+L$ independent sub-problems:
the $J$ sub-problems described in \eqref{eqn:subproblem:normal} and
the $L$ sub-problems described in \eqref{eqn:subproblem:blocking}
subject to $\alpha^{\phi^*}_{\ell,j}(C_{i_j}) \equiv 0$ ($j\in[J]$,$\ell\in[L]$) and $\alpha^{\phi^*}_{\ell,j}(n)\equiv 0$ ($n\in\mathscr{N}_{i_j}\backslash\{C_{i_j}\}$, $j\neq \mathcal{J}_{\ell}$, $\ell\in[L]$).
Equation \eqref{eqn:relaxed_opt} is obtained by solving sub-problems in \eqref{eqn:subproblem:normal}.

For a constant $g=g^*_j(\bm{\nu})$, $j\in[J]$ and $n\in\mathscr{N}_{i_j}\backslash\{C_{i_j},0\}$, where $g^*_j(\bm{\nu})$ is defined in \eqref{eqn:prop:g_star}, consider the Bellman equation \begin{multline}\label{eqn:relaxed_opt:3}
V^{g,\bm{\nu}}_j(n) = \max\limits_{\bm{\alpha}^{\phi}_{j}(n)\in\{0,1\}^{L}}\Biggl\{\frac{r^{\phi}_{j,n}(\bm{\nu})}{\sum_{\ell:j\in\mathscr{J}_{\ell}}\lambda_{\ell}\alpha^{\phi}_{\ell,j}(n)+\mu_j(n)} + \frac{\sum_{\ell:j\in\mathscr{J}_{\ell}}\lambda_{\ell}\alpha^{\phi}_{\ell,j}(n) V^{g,\bm{\nu}}_j(n+1)}{\sum_{\ell:j\in\mathscr{J}_{\ell}}\lambda_{\ell}\alpha^{\phi}_{\ell,j}(n)+\mu_j(n)}\\
+\frac{\mu_j(n) V^{g,\bm{\nu}}_j(n-1)}{\sum_{\ell:j\in\mathscr{J}_{\ell}}\lambda_{\ell}\alpha^{\phi}_{\ell,j}(n)+\mu_j(n)}\Biggr\},
\end{multline}
where $\bm{\alpha}^{\phi}_{j}(n)=(\alpha^{\phi}_{\ell,j}(n):\ \ell\in[L])$.
For $n=C_{i_j}$, $V^{g,\pmb{\nu}}_j(C_{i_j})$ also satisfies \eqref{eqn:relaxed_opt:3} with $\pmb{\alpha}^{\phi}_j(C_{i_j})\equiv \bm{0}$, and recall that $V^{g,\pmb{\nu}}_j(0) \equiv 0$.
For $n\in\mathscr{N}_{i_j}\backslash\{C_{i_j},0\}$, $j\in\mathcal{J}_{\ell}$, $\ell\in[L]$, equation~\eqref{eqn:relaxed_opt} is obtained by solving \eqref{eqn:relaxed_opt:3}.

Note that the conclusion from \eqref{eqn:relaxed_opt:3} is analyzed in the case with $\alpha^{\phi^*}_{\ell,j}(0) = 1$ for at least one $\ell\in\{\ell'\in[L]: j\in\mathcal{J}_{\ell'}\}$.
If it is not in this case, then sub-process $j$ will stay in state $0$ all the time whatever the action variables of other states will be.
Hence, the policy satisfying \eqref{eqn:relaxed_opt} is still optimal: the lemma is proved.
%\end{proof}
\endproof

Along with similar ideas, define 
$H^{\phi,g,\bm{\nu}}_j $ as the expected accumulated reward of the sub-process $j$ starting with $N^{\phi}_j(0)=0$ until the process enters state $0$ again, where the reward rate for state $n\in\mathscr{N}_{i_j}$ is specified as $r^{\phi}_{j,n}(\bm{\nu})-g$.
We refer to such a sojourn of the sub-process $j$ that starts and ends in state $0$ as the process $\mathcal{P}^{\phi}_j$.
Let $H^{g,\bm{\nu}}_j\coloneqq \max_{\phi\in\tilde{\Phi}_1}H^{\phi,g,\bm{\nu}}_j $. 

\begin{lemma}\label{lemma:relaxed_opt_zero}
When job sizes are exponentially distributed, for any policy $\phi^*\in\tilde{\Phi}_1$ that is optimal to the maximization problem in the Lagrangian dual function $L(\bm{\nu},\bm{\eta},\pmb{\gamma})$ of  the relaxed problem, 
\begin{equation}\label{eqn:relaxed_opt_zero}
\alpha^{\phi^*}_{\ell,j}(0)=\left\{
\begin{cases}
1, & \text{if } \nu_{\ell} < \lambda_{\ell}V^{g^{*}_j(\bm{\nu}),\bm{\nu}}_j(1),\\
a, & \text{if }\nu_{\ell} = \lambda_{\ell}V^{g^{*}_j(\bm{\nu}),\bm{\nu}}_j(1),\\
0, & \text{if }\nu_{\ell} > \lambda_{\ell}V^{g^{*}_j(\bm{\nu}),\bm{\nu}}_j(1),
\end{cases}\right.
\end{equation}
where $a$ can be any real number in $[0,1]$.
\end{lemma}
\vspace{3em}
\plainproof
For $\ell\in\{\ell'\in[L]|j\in\mathcal{J}_{\ell'}\}$, $j\in[J]$, let 
\begin{equation}
\bar{\lambda}_j^{\ell}=\sum_{\ell': j\in\mathcal{J}_{\ell'},\ell'\neq \ell}\lambda_{\ell'}\alpha^{\phi^*}_{\ell',j}(0)~\text{and}~\bar{\nu}_j^{\ell}=\sum_{\ell':j\in\mathcal{J}_{\ell'},\ell'\neq \ell}\nu_{\ell'}\alpha^{\phi^*}_{\ell',j}(0).
\end{equation}
If $\bar{\lambda}^{\ell}_j>0$, we obtain 
\begin{equation}\label{eqn:relaxed_opt_zero:1}
H^{g,\bm{\nu}}_j=\max\limits_{\alpha^{\phi}_{\ell,j}(0)\in\{0,1\}}\left\{\frac{\mu_j(0)-e^*\varepsilon_j(0)-g-\bar{\nu}^{\ell}_j-\nu_{\ell}\alpha^{\phi}_{\ell,j}(0)}{\bar{\lambda}^{\ell}_j+\lambda_{\ell}\alpha^{\phi}_{\ell,j}(0)}+V^{g,\bm{\nu}}_j(1)\right\},
\end{equation}
so that, for any $g\in\mathbb{R}$, and any policy $\phi^*\in\tilde{\Phi}_1$ that is optimal for the maximization problem in the Lagrangian dual function $L(\bm{\nu},\bm{\eta},\pmb{\gamma})$ of  the relaxed problem,
\begin{equation}
\alpha^{\phi^*}_{\ell,j}(0)=\left\{
\begin{cases}
1, & \text{if } \nu_{\ell} < \lambda_{\ell}\left(V^{g,\bm{\nu}}_j(1)-H^{g,\bm{\nu}}_j\right),\\
a, & \text{if }\nu_{\ell} = \lambda_{\ell}\left(V^{g,\bm{\nu}}_j(1)-H^{g,\bm{\nu}}_j\right),\\
0, & \text{if }\nu_{\ell} > \lambda_{\ell}\left(V^{g,\bm{\nu}}_j(1)-H^{g,\bm{\nu}}_j\right),
\end{cases}\right.
\end{equation}
where $a$ can be any real number in $[0,1]$.
Also, if $\bar{\lambda}^{\ell}_j > 0$, then $H^{g^*,\bm{\nu}}_j = 0$ and hence \eqref{eqn:relaxed_opt_zero} holds.

It remains to prove \eqref{eqn:relaxed_opt_zero} when $\bar{\lambda}^{\ell}_j=0$.
Let $\bar{r}^{\phi}_j$ represent the average reward received when the process $\mathcal{P}^{\phi}_{j}$ is in the states $n \in \mathscr{N}_{i_j}\backslash \{0\}$. 
From Lemma~\ref{lemma:relaxed_opt}, $\bar{r}^{\phi^*}_j$ with any optimal solution $\phi^*\in\tilde{\Phi}_1$ is independent from the value of $\bm{\alpha}^{\phi}_{j}(0)$. 
When $\bar{\lambda}^{\ell}_j=0$, $\alpha^{\phi^*}_{\ell,j}(0) = 1$ if 
\begin{equation*}
R_j(0) \leq \frac{1}{1+\frac{\lambda_{\ell}}{\mu_j(1)}}\Bigl(R_j(0)-\nu_{\ell}+\frac{\lambda_{\ell}}{\mu_{j}(1)}\bar{r}^{\phi^*}_j\Bigr),
\end{equation*}
where $R_j(n)=\mu_j(n)-e^*\varepsilon_j(n)$ for $n\in\mathscr{N}_{i_j}$; that is,
\begin{equation}\label{eqn:relaxed_opt_zero:2}
\nu_{\ell} \leq \frac{\lambda_{\ell}}{\mu_{j}(1)}(\bar{r}^{\phi^*}_j - R_j(0)).
\end{equation}

Because $g^*_j(\bm{\nu}) \geq R_j(0)$ (based on its definition in \eqref{eqn:prop:g_star}) and $V^{g^*_j(\bm{\nu}),\bm{\nu}}_j(1)=\lambda_{\ell}(\bar{r}^{\phi^*}_j-g^*_j(\bm{\nu}))/\mu_j(1)$ (Bellman equation), if $\nu_{\ell}\leq \lambda_{\ell}V^{g^*_j(\bm{\nu}),\bm{\nu}}_j(1)$ then \eqref{eqn:relaxed_opt_zero:2} holds: $\alpha^{\phi^*}_{\ell,j}(0)=1$.

We then show that $\nu_{\ell}\leq \lambda_{\ell}V^{g^*_j(\bm{\nu}),\bm{\nu}}_j(1)$ is necessary for $\alpha^{\phi^*}_{\ell,j}(0)=1$.
If $\alpha^{\phi^*}_{\ell,j}(0)=1$, \eqref{eqn:relaxed_opt_zero:2} holds, and
\begin{equation}\label{eqn:relaxed_opt_zero:3}
g^*_j(\bm{\nu}) = \frac{R_j(0)-\nu_{\ell}+\lambda_{\ell}\bar{r}^{\phi^*}_j/\mu_j(1)}{1+\lambda_{\ell}/\mu_j(1)}.
\end{equation}
By substituting \eqref{eqn:relaxed_opt_zero:3} in \eqref{eqn:relaxed_opt_zero:2}, we obtain $\nu_{\ell}\leq \lambda_{\ell}V^{g^*_j(\bm{\nu}),\bm{\nu}}_j(1)$.
This proves the lemma.
\endproof

\proof{Proposition~\ref{prop:relaxed_opt_precise}}
Consider a policy $\phi^*\in\tilde{\Phi}_1$ satisfying \eqref{eqn:relaxed_opt} and \eqref{eqn:relaxed_opt_zero}. There exists another policy $\phi_1$ such that, for $n\in\mathscr{N}_{i_j}\backslash\{C_{i_j}\}$,
\begin{equation}
\alpha^{\phi_1}_j(n) = \left\{\begin{cases}
\alpha^{\phi^*}_j(n), &\text{if } n < m,\\
0, & \text{otherwise},
\end{cases}\right.
\end{equation}
where $m = \min\bigl\{m'\in\mathscr{N}_{i_j}| \alpha^{\phi^*}_j(m')=0\bigr\}$. 
Since all states $n>m$ under policy $\phi^*$ are transient, policy $\phi_1$ leads to the same stationary distribution as $\phi^*$, which is optimal for the maximization in the Lagrangian dual function $L(\bm{\nu},\bm{\omega},\pmb{\gamma})$.

For state $n\in\mathscr{N}_{i_j}\backslash\{C_{i_j}\}$ ($j\in[J]$), let $m^*_j(n) \in \{n,n+1,\ldots,C_{i_j}-1\}$ represent the state such that $\alpha^{\phi^*}_j(m^*_j(n)+1) = 0$, and, if $m^*_j(n)\geq 1$, $\alpha^{\phi^*}_j(n+1)=\ldots=\alpha^{\phi^*}_j(m^*_j(n))=1$.
If $\alpha^{\phi^*}_j(n')=0$ for all $n' \geq n+1$, then $m^*_j(n) =n$.
From Lemmas~\ref{lemma:relaxed_opt} and \ref{lemma:relaxed_opt_zero}, for $n\in\mathscr{N}_{i_j}\backslash\{C_{i_j}\}$, we obtain \begin{equation}\label{eqn:relaxed_opt_precise:1}
V^{g^*_j(\bm{\nu}),\bm{\nu}}_j(n+1)-V^{g^*_j(\bm{\nu}),\bm{\nu}}_j(n)=
\left\{\begin{cases}
\frac{R_j(n+1)-g^*_j(\bm{\nu})}{\mu_{i_j}(n+1)}, & \text{if } m^*_j(n)=n,\\
\frac{R_j(n+1)-g^*_j(\bm{\nu})}{\mu_{i_j}(n+1)} + \sum\limits_{k=1}^{m^*_j(n)-n}\prod\limits_{\iota=1}^{k}\frac{\hat{\lambda}_{i_j}}{\mu_{i_j}(n+\iota)}\Bigl(\frac{R_j(n+k+1)-g^*_j(\bm{\nu})}{\mu_{i_j}(n+k+1)} - \nu\Bigr), & \text{otherwise},
\end{cases}\right.
\end{equation}
where $\hat{\lambda}_{i_j}\coloneqq \hat{\lambda}^0_{i_j}h =\sum_{\ell:i_j\in\mathcal{I}_{\ell}}\lambda_{\ell}$, and $g^*_j(\bm{\nu})$ is a given real number satisfying \eqref{eqn:prop:g_star}. Recall that $V^{g^*_j(\bm{\nu}),\bm{\nu}}_j(0)\equiv 0$. 
Here, we briefly explain how \eqref{eqn:relaxed_opt_precise:1} is achieved for the case with $m^*_j(n) > n$.
From the Bellman equation~\eqref{eqn:relaxed_opt:3}, the definition of $m^*_j(n)$, and Lemmas~\ref{lemma:relaxed_opt} and \ref{lemma:relaxed_opt_zero}, we obtain that, for $n'\leq m^*_j(n)-1$, 
\begin{multline}\label{eqn:final:review:1}
    V^{g^*_j(\pmb{\nu}),\pmb{\nu}}_j(n'+1) = \frac{R_j(n'+1)-\nu \hat{\lambda}_{i_j}-g^*_j(\pmb{\nu})}{\hat{\lambda}_{i_j}+\mu_{i_j}(n'+1)}+\frac{ \hat{\lambda}_{i_j}}{\hat{\lambda}_{i_j}+\mu_{i_j}(n'+1)}V^{g^*_j(\pmb{\nu}),\pmb{\nu}}_j(n'+2)\\ + \frac{ \mu_{i_j}(n'+1)}{\hat{\lambda}_{i_j}+\mu_{i_j}(n'+1)}V^{g^*_j(\pmb{\nu}),\pmb{\nu}}_j(n');
\end{multline}
and 
\begin{equation}\label{eqn:final:review:2}
    V^{g^*_j(\pmb{\nu}),\pmb{\nu}}_j(m^*_j(n)+1) = \frac{R_j(m^*_j(n)+1)-g^*_j(\pmb{\nu})}{\mu_{i_j}(m^*_j(n)+1)}+V^{g^*_j(\pmb{\nu}),\pmb{\nu}}_j(m^*_j(n)).
\end{equation}
Equations \eqref{eqn:final:review:1} and \eqref{eqn:final:review:2} are equivalent to
\begin{multline}\label{eqn:final:review:3}
    V^{g^*_j(\pmb{\nu}),\pmb{\nu}}_j(n'+1)-V^{g^*_j(\pmb{\nu}),\pmb{\nu}}_j(n')\\
    = \frac{R_j(n'+1)-g^*_j(\pmb{\nu})}{\mu_{i_j}(n'+1)} + \frac{ \hat{\lambda}_{i_j}}{\mu_{i_j}(n'+1)}\Bigl(V^{g^*_j(\pmb{\nu}),\pmb{\nu}}_j(n'+2)-V^{g^*_j(\pmb{\nu}),\pmb{\nu}}_j(n'+1) - \nu\Bigr),
\end{multline}
and 
\begin{equation}\label{eqn:final:review:4}
    V^{g^*_j(\pmb{\nu}),\pmb{\nu}}_j(m^*_j(n)+1)-V^{g^*_j(\pmb{\nu}),\pmb{\nu}}_j(m^*_j(n))
    = \frac{R_j(m^*_j(n)+1)-g^*_j(\pmb{\nu})}{\mu_{i_j}(m^*_j(n)+1)},
\end{equation}
respectively.
Let $\chi_{n'}\coloneqq V^{g^*_j(\pmb{\nu}),\pmb{\nu}}_j(n'+1)-V^{g^*_j(\pmb{\nu}),\pmb{\nu}}_j(n')$ for $n' = n,n+1,\ldots,m^*_j(n)$.
Based on \eqref{eqn:final:review:3}, 
\begin{multline}\label{eqn:final:review:5}
    \chi_n = \frac{R_j(n+1)-g^*_j(\pmb{\nu})}{\mu_{i_j}(n+1)} - \frac{\hat{\lambda}_{i_j}}{\mu_{i_j}(n+1)}\nu + \frac{\hat{\lambda}_{i_j}}{\mu_{i_j}(n+1)} \chi_{n+1}\\
    =\frac{R_j(n+1)-g^*_j(\pmb{\nu})}{\mu_{i_j}(n+1)} - \frac{\hat{\lambda}_{i_j}}{\mu_{i_j}(n+1)}\nu + \frac{\hat{\lambda}_{i_j}}{\mu_{i_j}(n+1)}\biggl(\frac{R_j(n+2)-g^*_j(\pmb{\nu})}{\mu_{i_j}(n+2)}- \frac{\hat{\lambda}_{i_j}}{\mu_{i_j}(n+2)}\nu + \frac{\hat{\lambda}_{i_j}}{\mu_{i_j}(n+2)}\chi_{n+2}\biggr)\\
    = \frac{R_j(n+1)-g^*_j(\pmb{\nu})}{\mu_{i_j}(n+1)}
    - \sum\limits_{k=n+1}^{m^*_j(n)}\biggl(\prod\limits_{\iota=n+1}^k\frac{\hat{\lambda}_{i_j}}{\mu_{i_j}(\iota)}\biggr)\nu
    +\sum\limits_{k=n+1}^{m^*_j(n)-1}\biggl(\prod\limits_{\iota=n+1}^k\frac{\hat{\lambda}_{i_j}}{\mu_{i_j}(\iota)}\biggr)\frac{R_j(k+1)-g^*_j(\pmb{\nu})}{\mu_{i_j}(k+1)}
    \\+\biggl(\prod\limits_{\iota=n+1}^{m^*_j(n)}\frac{\hat{\lambda}_{i_j}}{\mu_{i_j}(\iota)}\biggr)\chi_{m^*_j(n)}.
\end{multline}
Plugging \eqref{eqn:final:review:4} in \eqref{eqn:final:review:5}, we obtain
\begin{equation}
    \chi_n = \frac{R_j(n+1)-g^*_j(\pmb{\nu})}{\mu_{i_j}(n+1)}
    - \sum\limits_{k=n+1}^{m^*_j(n)}\biggl(\prod\limits_{\iota=n+1}^k\frac{\hat{\lambda}_{i_j}}{\mu_{i_j}(\iota)}\biggr)\nu
    +\sum\limits_{k=n+1}^{m^*_j(n)}\biggl(\prod\limits_{\iota=n+1}^k\frac{\hat{\lambda}_{i_j}}{\mu_{i_j}(\iota)}\biggr)\frac{R_j(k+1)-g^*_j(\pmb{\nu})}{\mu_{i_j}(k+1)},
\end{equation}
which leads to \eqref{eqn:relaxed_opt_precise:1} for the case with  $n < m^*_j(n)$.

If 
\begin{equation}\label{eqn:relaxed_opt_precise:2}
\nu \geq \max\limits_{\begin{subarray}~n'=n+1,\\n+2,\ldots,C_{i_j}\end{subarray}}\frac{R_j(n')-g^*_j(\bm{\nu})}{\mu_{i_j}(n')},
\end{equation}
where $R_j(n) = \mu_j(n)-e^*\varepsilon_j(n)$ for $n\in\mathscr{N}_{i_j}$, then, from \eqref{eqn:relaxed_opt_precise:1}, 
\begin{equation*}
V^{g^*_j(\bm{\nu}),\bm{\nu}}_j(n+1)-V^{g^*_j(\bm{\nu}),\bm{\nu}}_j(n)\leq \frac{R_j(n+1)-g^*_j(\bm{\nu})}{\mu_{i_j}(n+1)} \leq \nu; 
\end{equation*}
together with Lemmas~\ref{lemma:relaxed_opt} and \ref{lemma:relaxed_opt_zero}, $\alpha^{\phi^*}_j(n) = 0=\alpha^{\phi_1}_j(n)$.
It remains to prove that there exists a $H\in\mathbb{R}$ such that, for all $h>H$ and $n\in\mathscr{N}_{i_j}\backslash\{C_{i_j}\}$, if \eqref{eqn:relaxed_opt_precise:2} does not hold, $\alpha^{\phi_1}_j(n) = 1$; equation \eqref{eqn:prop:relaxed_opt_precise} is then led by Lemmas~\ref{lemma:relaxed_opt} and \ref{lemma:relaxed_opt_zero} and the continuity of $V^{g^*_j(\bm{\nu}),\bm{\nu}}_j(n+1)-V^{g^*_j(\bm{\nu}),\bm{\nu}}_j(n) - \nu$ in $\nu$, where recall that $\bm{\nu} = \bm{\lambda}\nu$.

For $n=C_{i_j}-1$, if \eqref{eqn:relaxed_opt_precise:2} does not hold, then, from Lemmas~\ref{lemma:relaxed_opt} and \ref{lemma:relaxed_opt_zero}, 
\begin{equation*}
V^{g^*_j(\bm{\nu}),\bm{\nu}}_j(C_{i_j})-V^{g^*_j(\bm{\nu}),\bm{\nu}}_j(C_{i_j}-1) = \frac{R_j(C_{i_j})-g^*_j(\bm{\nu})}{\mu_{i_j}(C_{i_j})} > \nu;
\end{equation*}
that is, $\alpha^{\phi^*}_j(C_{i_j}-1)=1$ and \eqref{eqn:prop:relaxed_opt_precise} holds for $n=C_{i_j}-1$.

We prove the remaining case with $n < C_{i_j}-1$ through iterations. Assume that, for sufficiently large $h$, \eqref{eqn:prop:relaxed_opt_precise}  holds for $n+1,n+2,\ldots,C_{i_j}-1$.
If $m^*_j(n)>n$, then, 
\begin{equation}\label{eqn:relaxed_opt_precise:2.5}
\nu\leq \frac{R_j(m^*_j(n)+1)-g^*_j(\bm{\nu})}{\mu_{i_j}(m^*_j(n)+1)}.
\end{equation}
Also,
by \eqref{eqn:relaxed_opt_precise:1} and Lammas~\ref{lemma:relaxed_opt} and \ref{lemma:relaxed_opt_zero},  if $m^*_j(n)>n$, and 
\begin{equation}\label{eqn:relaxed_opt_precise:3}
\nu < \frac{\frac{R_j(n+1)-g^*_j(\bm{\nu})}{\mu_{i_j}(n+1)}+\sum\limits_{k=1}^{m^*_j(n)-n}\prod\limits_{\iota=1}^k\frac{\hat{\lambda}_{i_j}}{\mu_{i_j}(n+\iota)}\frac{R_j(n+k+1)-g^*_j(\bm{\nu})}{\mu_{i_j}(n+k+1)}}{1+\sum\limits_{k=1}^{m^*_j(n)-n}\prod\limits_{\iota=1}^k\frac{\hat{\lambda}_{i_j}}{\mu_{i_j}(n+\iota)}},
\end{equation}
then $\alpha^{\phi^*}_j(n)=1$.
Let $\mathcal{v}_j(n)$ represent the value of the right-hand side of \eqref{eqn:relaxed_opt_precise:3}.
If \eqref{eqn:relaxed_opt_precise:2} does not hold and component $j$ is energy-efficiently unimodal,
for any $\sigma>0$, there exits $H\in\mathbb{R}$ such that, for all $h>H$, 
\begin{equation}
\mathcal{v}_j(n) \left\{
\begin{cases}
\geq \frac{R_j(m^*j(n)+1)-g^*_j(\bm{\nu})}{\mu_{i_j}(m^*_j(n)+1)}-\sigma, &\text{if }
\frac{R_j(m^*j(n)+1)-g^*_j(\bm{\nu})}{\mu_{i_j}(m^*_j(n)+1)} = \max\limits_{\begin{subarray}~n'=n+1, n+2,\\\ldots,m^*_j(n)+1\end{subarray}}\frac{R_j(n')-g^*_j(\bm{\nu})}{\mu_{i_j}(n')},\\
> \frac{R_j(m^*j(n)+1)-g^*_j(\bm{\nu})}{\mu_{i_j}(m^*_j(n)+1)}, &\text{otherwise}.
\end{cases}\right.
\end{equation}
In the former case, when \eqref{eqn:relaxed_opt_precise:2} does not hold, 
there exists $\sigma>0$ such that 
\begin{equation}
\mathcal{v}_j(n) \geq\frac{R_j(m^*j(n)+1)-g^*_j(\bm{\nu})}{\mu_{i_j}(m^*_j(n)+1)}-\sigma > \nu; 
\end{equation}
and, in the latter case, because of \eqref{eqn:relaxed_opt_precise:2.5}, $\mathcal{v}_j(n) > \nu$.
Accordingly, $\alpha^{\phi_1}_j(n)=\alpha^{\phi^*}_j(n)=1$.
If $m^*_j(n)=n$, and \eqref{eqn:relaxed_opt_precise:2} does not hold, then, from Bellman equation,
\begin{equation}
\frac{R_j(m^*_j(n)+1)-g^*_j(\bm{\nu})}{\mu_{i_j}(m^*_j(n)+1)} = V^{g^*_j(\bm{\nu}),\bm{\nu}}_j(n+1)-V^{g^*_j(\bm{\nu}),\bm{\nu}}_j(n).
\end{equation}
Together with Lemmas~\ref{lemma:relaxed_opt} and \ref{lemma:relaxed_opt_zero}, we obtain $\alpha^{\phi^*}_j(n)=1$.
This proves Proposition~\ref{prop:relaxed_opt_precise}.
\endproof

\section{Proof of Proposition~\ref{prop:relaxed_opt_explicit}}
\label{app:prop:relaxed_opt_explicit}

For $j\in[J]$ and $n\in\mathscr{N}_{i_j}\backslash\{C_{i_j}\}$, define, for $\nu\in\mathbb{R}$,
\begin{equation}\label{define:bar_f}
\bar{f}_{j,n}(\nu) \coloneqq \nu + \min\limits_{\begin{subarray}~n'=n+1,\\n+2,\ldots,C_{i_j}\end{subarray}}\Bigl(\frac{g^*_j(\nu \bm{\lambda})}{\mu_j(n')}-\frac{R_j(n')}{\mu_j(n')}\Bigr).
\end{equation}

\begin{lemma}\label{lemma:monotonicity}
For any $j\in[J]$, $n\in\mathscr{N}_{i_j}\backslash\{C_{i_j}\}$, there exists $\nu\in\mathbb{R}$ such that $\bar{f}_{j,n}(\nu)=0$. In particular, if job sizes are exponentially distributed and all computing components are energy-efficiently unimodal, there exists $H\in\mathbb{R}$ and $\upsilon^h_{j,n}\in\mathbb{R}$ such that, for all $h>H$, 
\begin{equation}\label{eqn:monotonicity}
\bar{f}_{j,n}(\nu) \left\{\begin{cases}
>0, &\text{if }\nu > \upsilon^h_{j,n},\\
=0, &\text{if }\nu = \upsilon^h_{j,n},\\
<0, &\text{otherwise}.
\end{cases}\right.
\end{equation}
\end{lemma}
\plainproof
From the definition, $g^*_j(\nu \bm{\lambda})$ is piece-wise linearly decreasing and continuous in $\nu\in\mathbb{R}$. From the definition in \eqref{define:bar_f}, $\bar{f}_{j,C_{i_j}-1}(\nu) $ is piece-wise linear and continuous in $\nu\in\mathbb{R}$. For $\nu\in\mathbb{R}$, 
\begin{equation}\label{eqn:monotonicity:1}
\frac{d^-}{d\nu}\bar{f}_{j,C_{i_j}-1}(\nu) 
=1 + \frac{d^-}{d\nu} \frac{g^*_j(\nu \bm{\lambda})}{\mu_j(C_{i_j})}
\geq 1 - \frac{\hat{\lambda}_j}{\mu_j(C_{i_j})}\pi_j,
\end{equation}
where $\frac{d^-}{d\nu}$ takes the left derivative, $\pi_j$ is the time proportion that $N^{\phi}_j(t) < C_{i_j}$ under a policy $\phi$, for which $\alpha^{\phi}_{\ell,j}(n)=1$ for all $n\in\mathscr{N}_{i_j}\backslash\{C_{i_j}\}$, $\ell\in\{\ell'\in[L]|j\in\mathcal{J}_{\ell'}\}$. We refer to such a policy $\phi$ as the \emph{fully activated (FA) policy}.
Under the FA policy, consider the underlying Markov chain of sub-process $j$:
\begin{equation}\label{eqn:monotonicity:2}
\frac{\hat{\lambda}_j}{\mu_j(C_{i_j})} \pi_j = \frac{\sum_{n=1}^{C_{i_j}-1}\prod_{k=1}^{n}\frac{\hat{\lambda}_j}{\mu_j(k)}\cdot \frac{\hat{\lambda}_j}{\mu_j(C_{i_j})}+\frac{\hat{\lambda}_j}{\mu_j(C_{i_j})}}{\sum_{n=1}^{C_{i_j}-1}\prod_{k=1}^{n}\frac{\hat{\lambda}_j}{\mu_j(k)}+1 + \prod_{k=1}^{C_{i_j}}\frac{\hat{\lambda}_j}{\mu_j(k)}}.
\end{equation}
Let $d$ be the difference between the numerator and denominator of the right-hand side of \eqref{eqn:monotonicity:2}. By calculation, $d< 0$.
That is, $\frac{d^-}{d\nu}\bar{f}_{j,C_{i_j}-1}(\nu) > 0$ for every point $\nu\in\mathbb{R}$ with existing $\bar{f}_{j,C_{i_j}-1}(\nu) $. 
Together with the piece-wise linearity of $\bar{f}_{j,C_{i_j}-1}(\nu)$, $\bar{f}_{j,C_{i_j}-1}(\nu)$ is monotonically increasing in $\nu\in\mathbb{R}$.  Since $\bar{f}_{j,C_{i_j}-1}(\nu)$ is also continuous in $\nu\in\mathbb{R}$ and $\bar{f}_{j,C_{i_j}-1}(\nu)$ tends to $\pm \infty$ as $\nu\rightarrow \pm\infty$, there exists $\nu$ such that $\bar{f}_{j,C_{i_j}-1}(\nu)=0$.
For any $n\in\mathscr{N}_{i_j}\backslash\{C_{i_j},C_{i_j}-1\}$, from the definition in \eqref{define:bar_f}, for any $\nu\in\mathbb{R}$,
$\bar{f}_{j,n}(\nu) \leq \bar{f}_{j,C_{i_j}-1}(\nu)$.
Since $\bar{f}_{j,n}(\nu)$ tends to $\pm\infty$ as $\nu\rightarrow \pm\infty$ and with the continuity in $\nu\in\mathbb{R}$, there exists a $\nu\in\mathbb{R}$ such that $\bar{f}_{j,n}(\nu)=0$ for any $n\in\mathscr{N}_{i_j}\backslash\{C_{i_j}, C_{i_j}-1\}$. Let $\upsilon_{j,n}$ represent the value of such a $\nu$ with $\bar{f}_{j,n}(\nu)=0$ for $n\in\mathscr{N}_{i_j}\backslash\{C_{i_j}\}$.

We then discuss the uniqueness of the zero point $\upsilon_{j,n}$.
When job sizes are exponentially distributed and the components are energy-efficiently unimodal,
for $n\in\mathscr{N}_{i_j}\backslash\{C_{i_j}\}$ and $\nu\in\mathbb{R}$, if $\bar{f}_{j,n}(\nu)\geq 0$, then, by Proposition~\ref{prop:relaxed_opt_precise}, there exists $H$ such that, for all $h>H$, an optimal solution $\phi^*$ exists and satisfies that, for any $n'\geq n+1$, $\alpha^{\phi^*}_j(n') = 0$. For such $\nu$ with $\bar{f}_{j,n}(\nu)\geq 0$, 
\begin{equation}
\frac{d^-}{d\nu}\bar{f}_{j,n}(\nu) \geq 1 - \frac{\hat{\lambda}_j}{\mu_j(n+1)}\sum\limits_{n'=0}^n\pi^{\phi^*}_j(n') > 0,
\end{equation} 
where $\pi^{\phi^*}_j(n')$ is the stationary distribution of sub-process $j$ under policy $\phi^*$, under which the sub-process achieves the maximal average reward $g^*(\nu\bm{\lambda})$.
Accordingly, the zero point $\upsilon_{j,n}$ is unique and  equation \eqref{eqn:monotonicity} is achieved by setting $\upsilon^h_{j,n}=\upsilon_{j,n}$.
\endproof

\proof{Proposition~\ref{prop:relaxed_opt_explicit}}
It is a straightforward result of Proposition~\ref{prop:relaxed_opt_precise} and Lemma~\ref{lemma:monotonicity} by substituting $\upsilon^*_{\ell,j}(n) = \lambda_{\ell}\upsilon^h_{j,n}$.
\endproof

\section{Proof of Proposition~\ref{prop:asym_opt:relaxed_opt}}\label{app:asym_opt:relaxed_opt}
\proof{Proposition~\ref{prop:asym_opt:relaxed_opt}}
From Proposition~\ref{prop:relaxed_opt_explicit}, for any $j\in[J]$, if $\bm{\nu}=\nu\bm{\lambda}$ for some $\nu\in\mathbb{R}$, then there exist a policy $\psi_j(n)\in\tilde{\Phi}_1$ satisfying, for all $\ell\in\{\ell'\in[L]|j\in\mathscr{J}_{\ell}\}$ and $n'\in\mathscr{N}_{i_j}\backslash\{C_{i_j}\}$, $\alpha_{\ell,j}^{\psi_j(n)}(n') = 1$; and $\alpha_{\ell,j}^{\psi_j(n)}(n') = 0$ otherwise, and $H>0$ such that, for all $h>H$, $\psi_j(n)$ is optimal for the problem described in \eqref{eqn:subproblem:normal}. In other words, for given $j\in[J]$, $n\in\mathscr{N}_{i_j}$ and multipliers $\bm{\nu}=\nu\bm{\lambda}$, there exists $H>0$ such that, for all $h>H$, $\Gamma_j(\bm{\nu}) = g^*_j(\bm{\nu})$. 

Substituting $\Gamma_j(\bm{\upsilon}^*_j(n))$ for $g^*_j(\bm{\upsilon}^*_j(n))$ in \eqref{eqn:monotonicity}, together with Proposition~\ref{prop:relaxed_opt_precise}, we prove the proposition. In particular, the indices $\upsilon^*_{\ell,j}(n) = \lambda_{\ell}\upsilon^h_{j,n}$ where $\upsilon^h_{j,n}$ is the zero point for $\bar{f}_{j,n}(\upsilon^h_{j,n}) = 0$.
\endproof

\section{Proof of Proposition~\ref{prop:exist_g}}\label{app:exist_g}
\label{app:exist_g}

\proof{Proposition~\ref{prop:exist_g}}
From the definition, for given $h\in\mathbb{N}_+$, $i\in[I]$, and $n\in\mathscr{N}_i\backslash\{C_i\}$, $f^h_{i,n}(\eta^0)$ is piece-wise linear and continuous in $\eta^0\in\mathbb{R}$, and
\begin{equation}
\Bigl\lvert \frac{d^-}{d\eta^0} f^h_{i,n}(\eta^0) \Bigr\rvert \leq 1 + \frac{\hat{\lambda}^0}{\mu_i(n+1)}\Bigl\lvert \frac{d^-}{d\eta^0} \bar{\Gamma}^{h,\text{EXP}}_i(\eta^0)\Bigr\rvert
=1 +  \frac{h\hat{\lambda}^0}{\mu_i(n+1)} \sum\limits_{n'=n^+}^1 \prod_{n''=n^+}^{n'}\frac{\mu_i(n'')}{h\hat{\lambda}^0}\pi^{\psi_j(n^+-1)}_j(n^+),
\end{equation}
where $j$ is any element in $\mathscr{J}_i$, $n^+\in\mathscr{N}_i \backslash\{0\}$ is the state such that the state $n^+-1$ maximizes the right-hand side of \eqref{eqn:asym_opt:Gamma}, $\pi_j^{\psi_j(n^+-1)}(n^+)$ is the steady state distribution of state $n^+$ under the policy $\psi_j(n^+-1)$. It follows that
\begin{equation}\label{eqn:prop:exist_g:1}
\Bigl\lvert \frac{d^-}{d\eta^0} f^h_{i,n}(\eta^0) \Bigr\rvert \leq 1 +  \frac{1}{\mu_i(n+1)} \sum\limits_{n'=n^+}^1 \frac{\prod_{n''=n^+}^{n'}\mu_i(n'')}{(h\hat{\lambda}^0)^{n^+-n'}}\pi^{\psi_j(n^+-1)}_j(n^+),
\end{equation}
for which the right-hand side is bounded with some finite constant $C\in\mathbb{R}_+$, and the $j$ is any element in $\mathscr{J}_i$. 
Similarly, for the limit case with $h\rightarrow +\infty$, \eqref{eqn:prop:exist_g:1} holds and the right-hand side is still bounded with some finite constant.
Thus, $f^h_{i,n}(\eta^0)$ is Lipschitz continuous for given $h\in\mathbb{N}_+\cup\{+\infty\}$, $i\in[I]$, and $n\in\mathscr{N}_i\backslash\{C_i\}$.

Observing that, for given $h\in\mathbb{N}_+$ and any $j\in\mathscr{J}_i$, $\bar{f}_{j,n}(\nu) = \frac{1}{\hat{\lambda}^0}f^h_{i,n}(\hat{\lambda}^0\nu)$. From Lemma~\ref{lemma:monotonicity}, there exists $\eta^0\in\mathbb{R}$ such that $f^h_{i,n}(\eta^0) = 0$. 
Since, for  any $h\in\mathbb{N}_+\cup\{+\infty\}$, $f^h_{i,n}(\eta^0)$ is Lipchitz continuous in $\eta^0\in\mathbb{R}$, is bounded with bounded $\eta^0\in\mathbb{R}$ and approaches $\pm\infty$ as $\eta^0\rightarrow \pm\infty$, there exists $\eta^0\in\mathbb{R}$ such that $f^h_{i,n}(\eta^0) = 0$ in the limit case with $h\rightarrow +\infty$. This proves the proposition.
\endproof

\section{Proof of Proposition~\ref{prop:e_star_equal}}
\label{app:prop:e_star_equal}
Define, for $e\in\mathbb{R}$,
\begin{equation}\label{eqn:global_attractor}
\bm{z}^{\psi(e)} \coloneqq \lim\limits_{h\rightarrow +\infty}\lim\limits_{t\rightarrow +\infty}\mathbb{E}[\bm{Z}^{\psi^h(\bar{\bm{\nu}},\bar{\bm{a}}^h,e),h}(t)],
\end{equation}
where $\bm{Z}^{\phi}(t)$ has been defined in \eqref{eqn:define_z}, and the existence of the limit is ensured by the existence of $\lim_{h\rightarrow +\infty} \bm{\pi}^{\psi^h(\bar{\bm{\nu}},\bar{\bm{a}},e)}_j$ for all $j\in[J]$.

\proof{Proposition~\ref{prop:e_star_equal}}
From \cite[Propositions 4 and 5]{fu2018restless}, we obtain that,
when the job sizes are exponentially distributed, for any $\delta>0$,
\begin{equation}\label{eqn:convergence}
\lim\limits_{h\rightarrow +\infty} \lim\limits_{T\rightarrow +\infty} \frac{1}{T} \mathbb{P}\Bigl\{\lVert\bm{Z}^{\varphi(e),h}(t)- \bm{z}^{\psi(e)}\rVert>\delta\Bigr\} = 0,
\end{equation}
where $\bm{Z}^{\varphi,h}(0)=\bm{z}^0$  and the policy $\varphi(e)\in\Phi$ is the index policy defined in \eqref{eqn:index_policy:1} and \eqref{eqn:index_policy:2} with substituted $\varphi(e)$ and $\bm{u}^*(e)$ for $\varphi$ and $\bm{\upsilon}^*$, respectively.
From \eqref{eqn:gamma_to_z} and \eqref{eqn:convergence}, for any $\delta>0$, there exists $H>0$ such that, for all $h>H$,
\begin{equation}\label{app:prop:convergence:0}
 \Gamma^{h,\psi^h(\bar{\bm{\nu}},\bar{\bm{a}}^h,e)}(e) - \delta\leq \Gamma^{h,\varphi(e)}(e) \leq \max\limits_{\phi\in\Phi}\Gamma^{h,\phi}(e),
\end{equation}
where the second inequality comes from $\varphi(e)\in \Phi$.
Together with \eqref{eqn:e_star_equal:condition}, for any $\delta>0$, there exists $H>0$ such that, for all $h>H$,
\begin{equation}\label{eqn:prop:convergence:1}
\max\limits_{\phi\in\tilde{\Phi}}\Gamma^{h,\phi}(e)-\delta\leq \max\limits_{\phi\in\Phi}\Gamma^{h,\phi}(e)\leq \max\limits_{\phi\in\tilde{\Phi}}\Gamma^{h,\phi}(e),
\end{equation}
where the second inequality is based on $\Phi\subset \tilde{\Phi}$.

Recall that, for any $j\in[J]$ and $e\in\mathbb{R}$, the existence of $\lim_{h\rightarrow+\infty} \bm{\pi}^{\psi^h(\bar{\bm{\nu}},\bar{\bm{a}}^h,e)}_j$ leading to the existence of $\lim_{h\rightarrow+\infty} \Gamma^{h,\psi^h(\bar{\bm{\nu}},\bar{\bm{a}},e)}(e)$. When \eqref{eqn:e_star_equal:condition} holds, $\lim_{h\rightarrow+\infty} \max_{\phi\in\tilde{\Phi}}\Gamma^{h,\phi}(e)$ also exists.
Based on \eqref{eqn:prop:convergence:1} and \eqref{eqn:e_star_equal:condition},  for any $e\in\mathbb{R}$,
\begin{equation}\label{eqn:prop:convergence:2}
\lim_{h\rightarrow+\infty} \max_{\phi\in\Phi}\Gamma^{h,\phi}(e)=\lim_{h\rightarrow+\infty} \Gamma^{h,\psi^h(\bar{\bm{\nu}},\bar{\bm{a}}^h,e)}(e) = \lim_{h\rightarrow+\infty} \max_{\phi\in\tilde{\Phi}}\Gamma^{h,\phi}(e).
\end{equation}

Let $\Gamma^{h,*}(e) \coloneqq \max_{\phi\in\tilde{\Phi}}\Gamma^{h,\phi}(e)$, which is piece-wise linear, continuous and decreasing in $e\in\mathbb{R}$.
Since $\lim_{h\rightarrow+\infty}\Gamma^{h,*}(e)$ exists for all given $e\in\mathbb{R}$ and $\lim_{h\rightarrow+\infty}\Gamma^{h,*}(e)$ tends to $\pm \infty$ as $e \rightarrow \mp \infty$, there exists a solution $e\in\mathbb{R}$ for $\lim_{h\rightarrow+\infty}\Gamma^{h,*}(e)=0$. 
Together with \eqref{eqn:e_star_equal:condition}, there exists a zero point $e\in\mathbb{R}$ such that $\lim_{h\rightarrow+\infty}\Gamma^{h,\psi^h(\bar{\bm{\nu}},\bar{\bm{a}}^h,e)}(e)=0$ and $\lim_{h\rightarrow+\infty}\Gamma^{h,\psi^h(\bar{\bm{\nu}},\bar{\bm{a}}^h,e)}(e)$ is continuous in $e\in\mathbb{R}$.

Let $e_0$ represent a specific real number such that $\lim_{h\rightarrow+\infty}\Gamma^{h,\psi^h(\bar{\bm{\nu}},\bar{\bm{a}}^h,e_0)}(e_0)=0$, and, for $h\in\mathbb{N}_+$ and $e\in\mathbb{R}$, let $\phi^h(e)$ represent an optimal solution such that $\Gamma^{h,\phi^h(e)}(e) = \max_{\phi\in\Phi}\Gamma^{h,\phi}(e)$.
From \eqref{eqn:prop:convergence:2}, for any $\delta>0$, there exists $H>0$ such that, for all $h>H$,
\begin{equation}\label{eqn:prop:convergence:3}
\lvert\Gamma^{h,\phi^h(e_0)}(e_0)\rvert = \lvert\mathfrak{L}^{\phi^h(e_0)}-e_0\mathfrak{E}^{\phi^h(e_0)}\rvert < \delta.
\end{equation}
That is, for any $\delta>0$, there exists $H>0$ such that, for all $h>H$, $e_0-\delta <  \frac{\mathfrak{L}^{\phi^h(e_0)}}{\mathfrak{E}^{\phi^h(e_0)}} < e_0+\delta$,
or equivalently, 
\begin{equation}\label{eqn:prop:convergence:4}
\lim_{h\rightarrow +\infty} \frac{\mathfrak{L}^{\phi^h(e_0)}}{\mathfrak{E}^{\phi^h(e_0)}}  = e_0.
\end{equation} 
Based on \eqref{eqn:prop:convergence:3}, for any $\phi\in\Phi$ and $\delta>0$, there exists $H>0$ such that, for all $h>H$,
\begin{equation}
\mathfrak{L}^{\phi}-e_0\mathfrak{E}^{\phi} \leq \mathfrak{L}^{\phi^h(e_0)}-e_0\mathfrak{E}^{\phi^h(e_0)} < \delta.
\end{equation}
It follows that, for any $\phi\in\Phi$ and $\delta>0$, there exists $H>0$ such that, for all $h>H$, $\frac{\mathfrak{L}^{\phi}}{\mathfrak{E}^{\phi}} < e_0+\delta$; or equivalently, for any $\phi\in\Phi$,
\begin{equation}\label{eqn:prop:convergence:5}
\lim\limits_{h\rightarrow +\infty} \frac{\mathfrak{L}^{\phi}}{\mathfrak{E}^{\phi}} \leq e_0.
\end{equation}
From \eqref{eqn:prop:convergence:4} and \eqref{eqn:prop:convergence:5} and the definition of $e^*$ in \eqref{eqn:e_star}, we obtain that 
\begin{equation}\label{eqn:prop:convergence:6}
\lim\limits_{h\rightarrow +\infty} e^* = \lim\limits_{h\rightarrow +\infty} \max\limits_{\phi\in\Phi} \frac{\mathfrak{L}^{\phi}}{\mathfrak{E}^{\phi}} = e_0.
\end{equation}
This proves \eqref{eqn:e_star_equal}.

We then discuss the uniqueness of the zero point $e_0$ for $\lim_{h\rightarrow+\infty}\Gamma^{h,\psi^h(\bar{\bm{\nu}},\bar{\bm{a}}^h,e_0)}(e_0)=0$.
If there exists $e_1\neq e_0$ satisfying $\lim_{h\rightarrow+\infty}\Gamma^{h,\psi^h(\bar{\bm{\nu}},\bar{\bm{a}}^h,e_1)}(e_1)=0$, then, from the above discussion, 
\begin{equation}
 e_1=\lim\limits_{h\rightarrow +\infty} \max\limits_{\phi\in\Phi} \frac{\mathfrak{L}^{\phi}}{\mathfrak{E}^{\phi}} = e_0.
\end{equation}
Hence, $e_0$ is the unique solution for $\lim_{h\rightarrow+\infty}\Gamma^{h,\psi^h(\bar{\bm{\nu}},\bar{\bm{a}}^h,e_0)}(e_0)=0$. 

Recall that, for given $h\in\mathbb{N}_+$, $\Gamma^{h,*}(e) \coloneqq \max_{\phi\in\tilde{\Phi}} \Gamma^{h,\phi}(e)$, which is piece-wise linear in $e\in\mathbb{R}$ with the left derivative 
\begin{equation}
\Bigl\lvert\frac{d^-}{de}\Gamma^{h,*}(e)\Bigr\rvert < \Bigl\lvert\frac{1}{h}\sum_{j\in\mathcal{J}_{\ell}}\bm{\pi}^{\phi^*}_j\cdot \bm{\varepsilon}_{i_j} \Bigr\rvert \leq \max\limits_{\begin{subarray}~i\in[I],\\n\in\mathscr{N}_i\end{subarray}}\varepsilon_i(n) \leq C,
\end{equation}
where $C\in\mathbb{R}_+$ is a constant independent from $h$, $\phi^*$ is an optimal policy satisfying $\Gamma^{h,\phi^*}(e) = \Gamma^{h,*}(e)$, and $\bm{\varepsilon}_i \coloneqq (\varepsilon_i(n): n\in\mathscr{N}_i)$ for $i\in[I]$.
That is, for given $h\in\mathbb{N}_+$, $\Gamma^{h,*}(e)$ is Lipschitz continuous in $e\in\mathbb{R}$ with a bounded Lipschitz constant independent from $h$.

For any given $e\in\mathbb{R}$, from \eqref{eqn:e_star_equal:condition}, the function $\Gamma(e) = \lim_{h\rightarrow +\infty} \max_{\phi\in\tilde{\Phi}} \Gamma^{h,\phi}(e) = \lim_{h\rightarrow+\infty} \Gamma^{h,*}(e) $. 
It follows that, for any $e\in\mathbb{R}$,
\begin{equation}
\lim\limits_{\Delta\rightarrow 0} \Bigl\lvert \Gamma(e+\Delta)-\Gamma(e)\Bigr\rvert 
= \lim\limits_{\Delta\rightarrow 0} \Bigl\lvert \lim\limits_{h\rightarrow +\infty} \bigl(\Gamma^{h,*}(e+\Delta) - \Gamma^{h,*}(e)\bigr)\Bigr\rvert
\leq \lim\limits_{\Delta\rightarrow 0} C\Delta = 0,
\end{equation}
where $C$ is the bounded Lipschitz constant for $\Gamma^{h,*}(e)$. That is, $\Gamma(e)$ is continuous in $e\in\mathbb{R}$. Similarly,
\begin{equation}
\lim\limits_{\Delta\uparrow 0} \frac{1}{\Delta}\Bigl\lvert \Gamma(e+\Delta)-\Gamma(e)\Bigr\rvert 
=\lim\limits_{\Delta\uparrow 0} \frac{1}{\Delta}\Bigl\lvert \lim\limits_{h\rightarrow +\infty} \bigl(\Gamma^{h,*}(e+\Delta) - \Gamma^{h,*}(e)\bigr)\Bigr\rvert 
\leq \lim\limits_{\Delta\uparrow 0} \frac{1}{\Delta}\cdot C\Delta = C,
\end{equation}
leading to the Lipschitz continuity of $\Gamma(e)$.
This proves the proposition.
\endproof

\section{Proof of Lemma~\ref{lemma:convergence_relaxed_opt}}\label{app:lemma:convergence_relaxed_opt}

\proof{Lemma~\ref{lemma:convergence_relaxed_opt}}
The lemma can be proved by showing that constructing $\bm{\nu}$, $\pmb{\gamma}$, and $\phi^*$, where $\bm{\alpha}^{\phi^{*}}_{j}$ and $\bar{\alpha}_{\ell}^{\phi^*}$ maximize the the objectives defined in \eqref{eqn:subproblem:normal} and \eqref{eqn:subproblem:blocking};  and showing that such $\bm{\nu}$, $\pmb{\gamma}$, and $\phi^*$ achieve the complementary slackness for the relaxed problem defined by \eqref{eqn:opt}, \eqref{eqn:constraint:relax:action}, \eqref{eqn:constraint:relax:capacity} and \eqref{eqn:constraint:relax:boundary}.

We start with the case with Condition~\ref{cond:heavy}.
Let 
\begin{enumerate}
\item $\nu_{\ell}/\lambda_{\ell} = \nu$ where 
\begin{equation}\label{eqn:for_heavy_cond}
\nu =\min\left\{0, \min\limits_{j\in\mathcal{J}_{\ell}}\frac{\upsilon^*_{\ell,j}(C_{i_j}-1)}{\lambda_{\ell}}\right\}
\end{equation} 
with $\upsilon^*_{\ell,j}(n)$ ($n\in\mathscr{N}_{i_j}$, $j\in\mathcal{J}_{\ell}$, $\ell\in[L]$) given by \eqref{eqn:asym_opt:index};
and
\item $\gamma_{\ell} = - \lambda_{\ell}\nu$, $\ell\in[L]$.
\end{enumerate}
In this context, from Proposition~\ref{prop:relaxed_opt_explicit}, for sufficiently large $h$, there is an optimal policy $\phi^{*}$ for the problem defined in \eqref{eqn:subproblem:normal}, satisfying $\alpha^{\phi^{*}}_{\ell,j}(n) = 1$ for all $n\in\mathscr{N}_{i_j}\backslash\{C_{i_j}\}$, $j\in\mathcal{J}_{\ell}$. 
Let $\bar{\alpha}^{\phi^*}_{\ell}=1-A_{\ell}$, 
where $0\leq 1-A_{j} \leq 1$ under Condition~\ref{cond:heavy}.
In other words, constraints~\eqref{eqn:constraint:relax:action} and \eqref{eqn:constraint:relax:boundary} are satisfied with equality.
Recall that Constraint~\eqref{eqn:constraint:relax:capacity} has been guaranteed by setting $\bm{\eta}$ to be sufficiently large: the complementary slackness of the relaxed problem is achieved.
Hence, $\phi^{*}$ will also maximize the primal problem defined by \eqref{eqn:opt}, \eqref{eqn:constraint:relax:action}, \eqref{eqn:constraint:relax:capacity} and \eqref{eqn:constraint:relax:boundary}.

We now consider the case with $L=1$. 
If $A_{\ell} \leq 1$, then this is a special case for $L\geq 1$ with Condition~\ref{cond:heavy}.
It remains to discuss the case with $L=1$ and $A_{\ell} >1$.
From Proposition~\ref{prop:relaxed_opt_explicit}, for sufficiently large $h$, there is an optimal policy
$\phi^{*}$ that maximizes the problem defined in
\eqref{eqn:subproblem:normal}, for which $\bm{\alpha}^{\phi^*}_{j}$ is determined by \eqref{eqn:indexability} and \eqref{eqn:asym_opt:index}.
Since $A_{\ell}> 1$, there must exist a $\nu_{\ell}$ such that 
\begin{equation*}
\nu_{\ell} \geq
\min_{\begin{subarray}~j\in\mathcal{J}_{\ell}\\ n\in\mathscr{N}_{i_j}\end{subarray}}\upsilon^*_{\ell,j}(n),
\end{equation*}
where the $\ell$ is the only element in $[L]$,
and
$\sum_{j\in\mathcal{J}_{\ell}}\bm{\pi}^{\phi^{*}}_{j}\cdot \bm{\alpha}^{\phi^{*}}_{j}
=1$.  

If $A_{\ell}> 1$, then let
$\gamma_{\ell} > \max\left\{0,-\nu_{\ell}\right\}$, which guarantees that policy $\phi^*$ with $\bar{\alpha}_{\ell}^{\phi^*} < 0$ maximizes 
the objective defined in
\eqref{eqn:subproblem:blocking}.
Then, $I (\bar{\alpha}^{\phi^*}_{\ell})=0=I(1-A_{\ell})$.
Constraints \eqref{eqn:constraint:relax:action} and
\eqref{eqn:constraint:relax:boundary} achieve equality.  
The complementary
slackness conditions are also satisfied under this setting.
The lemma is then proved.
\endproof

\section{Proof of Proposition~\ref{prop:varphi_closeness}}
\label{app:prop:varphi_closeness}
\proof{Proposition~\ref{prop:varphi_closeness}}
Along similar lines as the proof of Proposition~\ref{prop:e_star_equal} in Appendix~\ref{app:prop:e_star_equal}, based on \eqref{eqn:e_star_equal:condition}, we obtain \eqref{app:prop:convergence:0} and \eqref{eqn:prop:convergence:2}. In other words, for any $e\in\mathbb{R}$,
\begin{equation}
\lim\limits_{h\rightarrow +\infty} \Gamma^{h,\varphi(e)}(e) 
=\lim\limits_{h\rightarrow +\infty} \Gamma^{h,\psi^h(\bar{\bm{\nu}},\bar{\bm{a}}^h,e)}(e)
=\Gamma(e).
\end{equation}
From Proposition~\ref{prop:e_star_equal}, if \eqref{eqn:e_star_equal:condition} holds, $\Gamma(e)$ is Lipschitz continuous in $e\in\mathbb{R}$. It follows with \eqref{eqn:prop:varphi_closeness} and proves the proposition.
\endproof

\section{Proof of Proposition~\ref{prop:asym_opt}}\label{app:prop:asym_opt}
\proof{Proposition~\ref{prop:asym_opt}}
The proposition is proved by invoking \cite[Propositions 4 and 5]{fu2018restless}, which ensures the existence of a global attractor for the process $\{\bm{Z}^{\varphi(e),h}(t),t\geq 0\}$, where $\bm{Z}^{\varphi(e),h}(t)$ has been defined in \eqref{eqn:define_z} and the index policy $\varphi(e)$ has been defined in \eqref{eqn:index_policy:1} and \eqref{eqn:index_policy:2} with substituted  $\varphi(e)$ and $u^*_{\ell,i_j}(n,e)$ for $\varphi$ and $\upsilon^*_{\ell,j}(n)$, respectively.
More precisely, from \cite[Propositions 4 and 5]{fu2018restless}, 
for any $\delta > 0$, we obtain \eqref{eqn:convergence} with $\bm{Z}^{\varphi(e),h}(t) = \bm{z}^0$. The $\bm{z}^{\psi(e)}$ is defined in \eqref{eqn:global_attractor} and is the global attractor for the process $\{\bm{Z}^{\varphi(e),h}(t),t\geq 0\}$.
That is, together with Lemma~\ref{lemma:relaxed_opt}, we obtain
\begin{equation}
\lim\limits_{h\rightarrow+\infty} \Gamma^{h,\varphi(e)}(e) = \lim\limits_{h\rightarrow+\infty} \max\limits_{\phi\in\tilde{\Phi}}\Gamma^{h,\phi}(e),
\end{equation}
which, since $\varphi(e)\in \Phi$ and $\Phi\subset \tilde{\Phi}$, indicates \eqref{eqn:asym_opt} with substituted $\varphi(e)$ for $\phi$. It proves the proposition.
\endproof

\section{Proof of Proposition~\ref{prop:exp_diminishing}}\label{app:prop:exp_diminishing}
\proof{Proposition~\ref{prop:exp_diminishing}}
Along similar lines as the proof of \cite[Proposition 2]{fu2020energy}, by invoking \cite[Theorem 4.1 in Chapter 7 \& Theorem 3.3 in Chapter 3]{freidlin2012random}, there exists a deterministic process $\bm{z}(t)$ taking values in $\mathscr{Z}$ such that,
for any $\delta>0$, there exist $s>0$ and $H > 0$ satisfying, for all $h>H$, 
\begin{equation}
\lim\limits_{T\rightarrow +\infty}\frac{1}{T}\int_0^T\mathbb{P}\Bigl\{\Bigl\lVert \bm{Z}^{h,\varphi(e)}(t) -\bm{z}(t)\Bigr\rVert > \delta \Bigr\} dt\leq e^{-sh}.
\end{equation} 
Together with \eqref{eqn:convergence}, we obtain that, for any $\delta>0$, there exist $s>0$ and $H > 0$ such that, for all $h>H$, \eqref{eqn:exp_diminishing} is achieved. It proves the proposition.
\endproof

\section{Simulation Settings for Section~\ref{subsubsec:case2}}\label{app:google-settings}
Consider a system with $I=10$ clusters, where the service and energy consumption rates of components are instances of pseudo-random variables:
\begin{itemize}
\item $\mu_1(C_1) = 2.425$, $\varepsilon_1(0)=0.0655$, $\frac{\mu_1(C_1)}{\varepsilon_1(C_1)-\varepsilon_1(0)}=13.699$
\item $\mu_2(C_2)=1.620$, $\varepsilon_2(0)=0.0333$, $\frac{\mu_2(C_2)}{\varepsilon_2(C_2)-\varepsilon_2(0)}=15.338$
\item $\mu_3(C_3)=1.758$, $\varepsilon_3(0)=0.0315$, $\frac{\mu_3(C_3)}{\varepsilon_3(C_3)-\varepsilon_3(0)}=14.845$
\item $\mu_4(C_4)=1.600$, $\varepsilon_4(0)=0.0189$, $\frac{\mu_4(C_4)}{\varepsilon_4(C_4)-\varepsilon_4(0)}=18.562$
\item $\mu_5(C_5)=1.728$, $\varepsilon_5(0)=0.0116$, $\frac{\mu_5(C_5)}{\varepsilon_5(C_5)-\varepsilon_5(0)}=26.225$
\item $\mu_6(C_6)=1.668$, $\varepsilon_6(0)=0.0069$, $\frac{\mu_6(C_6)}{\varepsilon_6(C_6)-\varepsilon_6(0)}=33.166$
\item $\mu_7(C_7)=2.390$, $\varepsilon_7(0)=0.0055$, $\frac{\mu_7(C_7)}{\varepsilon_7(C_7)-\varepsilon_7(0)}=43.127$
\item $\mu_8(C_8)=2.116$, $\varepsilon_8(0)=0.0026$, $\frac{\mu_8(C_8)}{\varepsilon_8(C_8)-\varepsilon_8(0)}=51.306$
\item $\mu_9(C_9)=2.416$, $\varepsilon_9(0)=0.0011$, $\frac{\mu_9(C_9)}{\varepsilon_9(C_9)-\varepsilon_9(0)}=70.356$
\item $\mu_{10}(C_{10})=2.224$, $\varepsilon_{10}(0)=0$, $\frac{\mu_{10}(C_{10})}{\varepsilon_{10}(C_{10})-\varepsilon_{10}(0)}=97.625$
\end{itemize}
and, for all clusters $i\in[I]$, the service and energy consumption rates for states $n\in\mathscr{N}_i\backslash \{0,C_i\}$ are given by $\mu_i(n) = \frac{n}{n+1}\mu_i(n+1)$ and $\varepsilon_i(n)=\frac{n}{n+1}(\varepsilon_i(n+1)-\varepsilon_i(0))+\varepsilon_i(0)$, respectively.
In particular, the unit of service rates of all the clusters is $10^{-10} $ number of jobs per second: they are normalized to be sufficiently small that we can observe a positive number of blocked jobs in the simulations for Figure~\ref{fig:google-blocking} and the heavy traffic condition can be achieved during peak hours for the scenario discussed in Section~\ref{subsubsec:case2}.
There are $L=4$ classes of jobs and, for all the classes $\ell\in[L]$, the sets of clusters able to serve an $\ell$-job are
$\mathcal{I}_1 = \{1,5,6,10\}$, $\mathcal{I}_2 = \{1,2,3,4,5,7,8,9\}$, $\mathcal{I}_3 = \{1,6,7,10\}$, and $\mathcal{I}_4 = \{2\}$.

\section{Table of Key Notation}\label{app:symbols}
We provide a table of important notation in Table~\ref{table:symbols}.
\begin{table}[!ht]
	\caption{Important Symbols}\label{table:symbols}
	\begin{tabular}{p{1.4cm}p{7cm}p{1.4cm}p{6.9cm}}
	    \hline
		\multicolumn{4}{l}{Sets}\\
		\hline
		$\mathbb{R},\mathbb{R}_+,\mathbb{R}_0$& Sets of all real numbers, positive real numbers and non-negative real numbers, respectively&
		$\mathbb{N},\mathbb{N}_+,\mathbb{N}_0$& Sets of all integers, positive integers, and non-negative integers, respectively\\
		$[N]$& Set $\{1,2,\ldots,N\}$ given $N\in\mathbb{N}_+$&
		$\mathscr{J}_{\ell}$&Set of available components for jobs in class $\ell\in[L]$\\
		$\mathscr{N}_i$&Set of all states for the bandit process associated with a component in cluster $i\in[I]$, also referred to as the state space of the bandit process&
		$\mathscr{N}$&Set of all states for the entire system, also referred to as the state space of the entire system\\
		$\Phi$& Set of policies that are applicable to the original problem&
		$\tilde{\Phi}$& Set of policies that are applicable to the relaxed problem\\
		$\tilde{\Phi}_1$ & Set of policies that are determined by action variables $\bm{\alpha}_j^{\phi}\in[0,1]^{L|\mathscr{N}_{i_j}|}$ ($j\in[J]\cup\{j_1,j_2,\ldots,j_L\}$) &\\
		\hline
		\multicolumn{4}{l}{Real Numbers and Vectors}\\
		\hline
		$I$& Number of component clusters &
		$J$& Number of physical components\\
		$L$& Number of job classes &
		$\lambda_{\ell}$&Arrival rate of jobs in class $\ell\in[L]$\\
		$j_{\ell}$& Label of the virtual component for job class $\ell\in[L]$&
		$C_i$& Maximal number of jobs that can be served simultaneously by a component in cluster $i\in[I]$\\
		$\mu_i(n)$& Service rate of a component in cluster $i$ when there are $n$ jobs being served by the component &
		$\varepsilon_i(n)$& Power consumption rate of a component in cluster $i$ when there are $n$ jobs being served by the component\\		
		$h\in\mathbb{N}_+$&Scaling parameter of the entire system&
		$M_i^0\cdot h$& Number of components in cluster $i\in[I]$, given the scaling parameter $h\in\mathbb{N}_+$\\
		$\lambda_{\ell}^0$& Arrival rate of jobs in class $\ell$ when $h=1$ and $\lambda_{\ell}^0  = \lambda_{\ell}/h$&
		$i_j$& Label of the cluster that component $j\in[J]$ belongs to\\
		$a^{\phi}_{\ell,j}(\bm{n})$& Action variable, for the original problem, for job class $\ell\in[L]$ of the bandit process associated with component $j\in[J]\cup\{j_{\ell}\}$  when the entire server farm system is in state $\bm{n}\in\mathscr{N}$ under policy $\phi\in\Phi$&
        $\alpha^{\phi}_{\ell,j}(n)$ & Action variable, for the relaxed problem, for job class $\ell\in[L]$ of the bandit process associated with component  $j\in[J]\cup\{j_{\ell}\}$ when the bandit process is in state $n\in\mathscr{N}_{i_j}$ under policy $\phi\in\tilde{\Phi}$ \\
        $\mathfrak{L}^{\phi}$ & Long-run average job throughput under policy $\phi$ &
        $\mathfrak{E}^{\phi}$ & Long-run average power consumption rate under policy $\phi$\\
        $e^*$ & Given real number used to translate objective in \eqref{eqn:obj} to \eqref{eqn:opt} &
        $\bm{\nu}$ & Vector of Lagrangian multipliers for the relaxed problem \\
        $\bm{\upsilon}*$ & Vector of indices &
        $\bm{u}(e)$ & Output of Algorithm~\ref{algo:nu_star} when plug-in $e$ as the given value for $e^*$\\
        $\bm{u}^*(e)$ & Vector of solutions $u^*_{\ell,i}(n,e)$ of $f^h_{i,n}\bigl(\frac{\hat{\lambda}^0_i}{\lambda^0_{\ell}}u^*_{\ell,i}(n,e)\bigr)=0$ for all $\ell\in[L]$, $i\in[I]$ and $n\in\mathscr{N}_i\backslash\{C_i\}$& &\\
        \hline
		\multicolumn{4}{l}{Random Variables}\\
		\hline
		$N_j^{\phi}(t)$ & State variable of the bandit process associated with component $j\in[J]$ under policy $\phi$, representing the number of jobs being served by component $j$ at time $t$ under policy $\phi$ &
		$\bm{N}^{\phi}(t)$ & State variable of the entire system under policy $\phi$, defined as $\bm{N}^{\phi}(t) = (N^{\phi}_j(t):j\in[J])$\\
		\hline
	\end{tabular}
\end{table}

\bibliographystyle{references/IEEEtran}
\bibliography{references/IEEEabrv,references/bib1}

\begin{IEEEbiography}[{\includegraphics[width=1in,height=1.25in,clip,keepaspectratio]{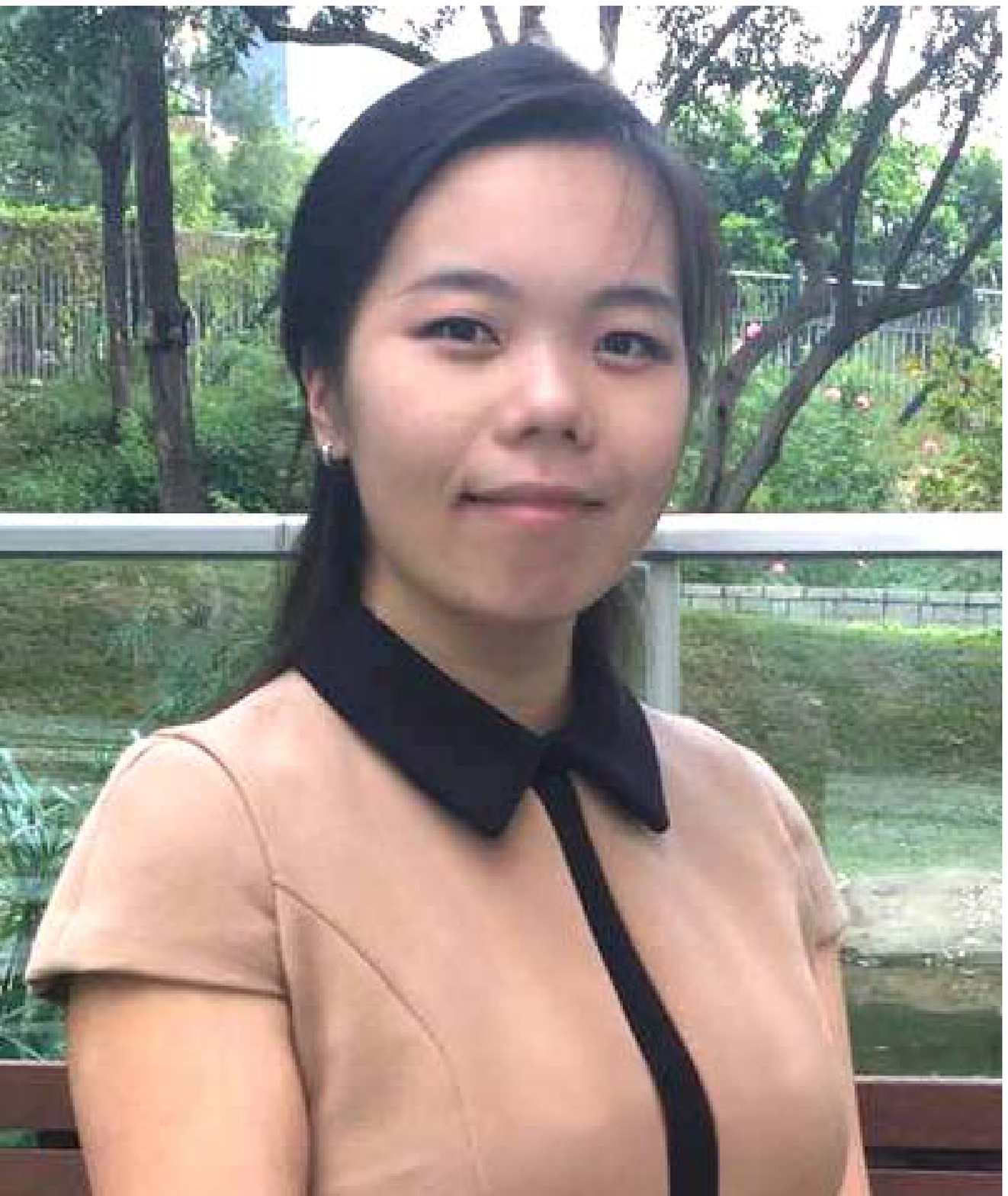}}]{Jing Fu} (S'15--M'16) 
received the B.Eng. degree in
computer science from Shanghai Jiao Tong University,
Shanghai, China, in 2011, and the Ph.D. degree
in electronic engineering from the City University
of Hong Kong in 2016. She has been with the
School of Mathematics and Statistics, the University
of Melbourne as a Post-Doctoral Research Associate
from 2016 to 2019. She is now with RMIT University, Australia, as a lecturer. 
Her research interests
now include energy-efficient networking/scheduling,
resource allocation in large-scale networks, semi-
Markov/Markov decision processes, restless multiarmed
bandit problems, stochastic optimization.
\end{IEEEbiography}

\begin{IEEEbiography}[{\includegraphics[width=1in,height=1.25in,clip,keepaspectratio]{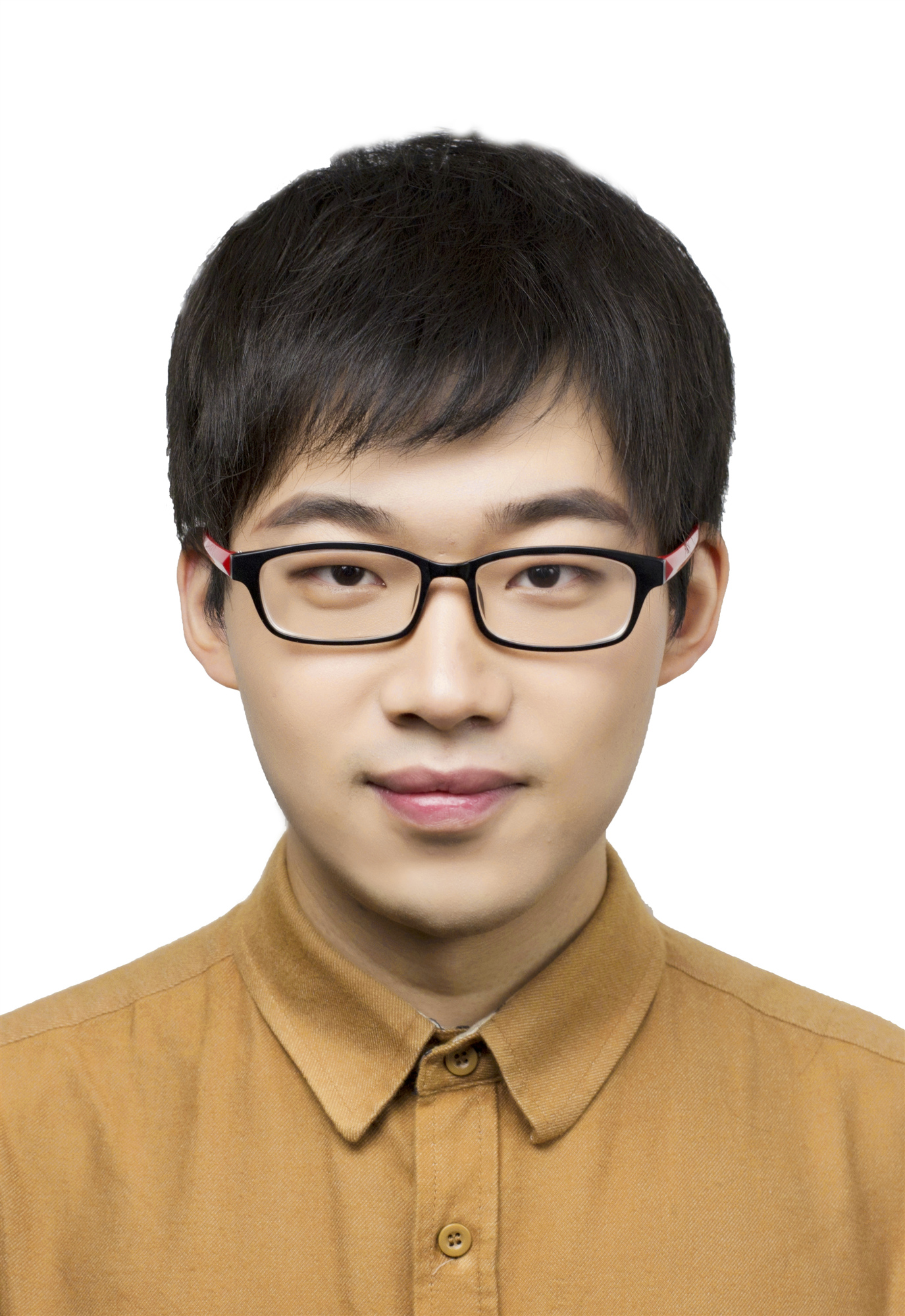}}]{Xinyu Wang } 
received the B.Eng. degree in Electronic Science and Technology from Zhejiang University, Hangzhou, and the M.Sc. degree in Electronic Information Engineering from City University of Hong Kong, in 2016 and 2018 respectively. He is currently a Ph.D. candidate with the Department of Electrical Engineering, City University of Hong Kong, China. His current research interests are path planning and resource allocation in telecommunication networks.
\end{IEEEbiography}

\begin{IEEEbiography}[{\includegraphics[width=1in,height=1.25in,clip,keepaspectratio]{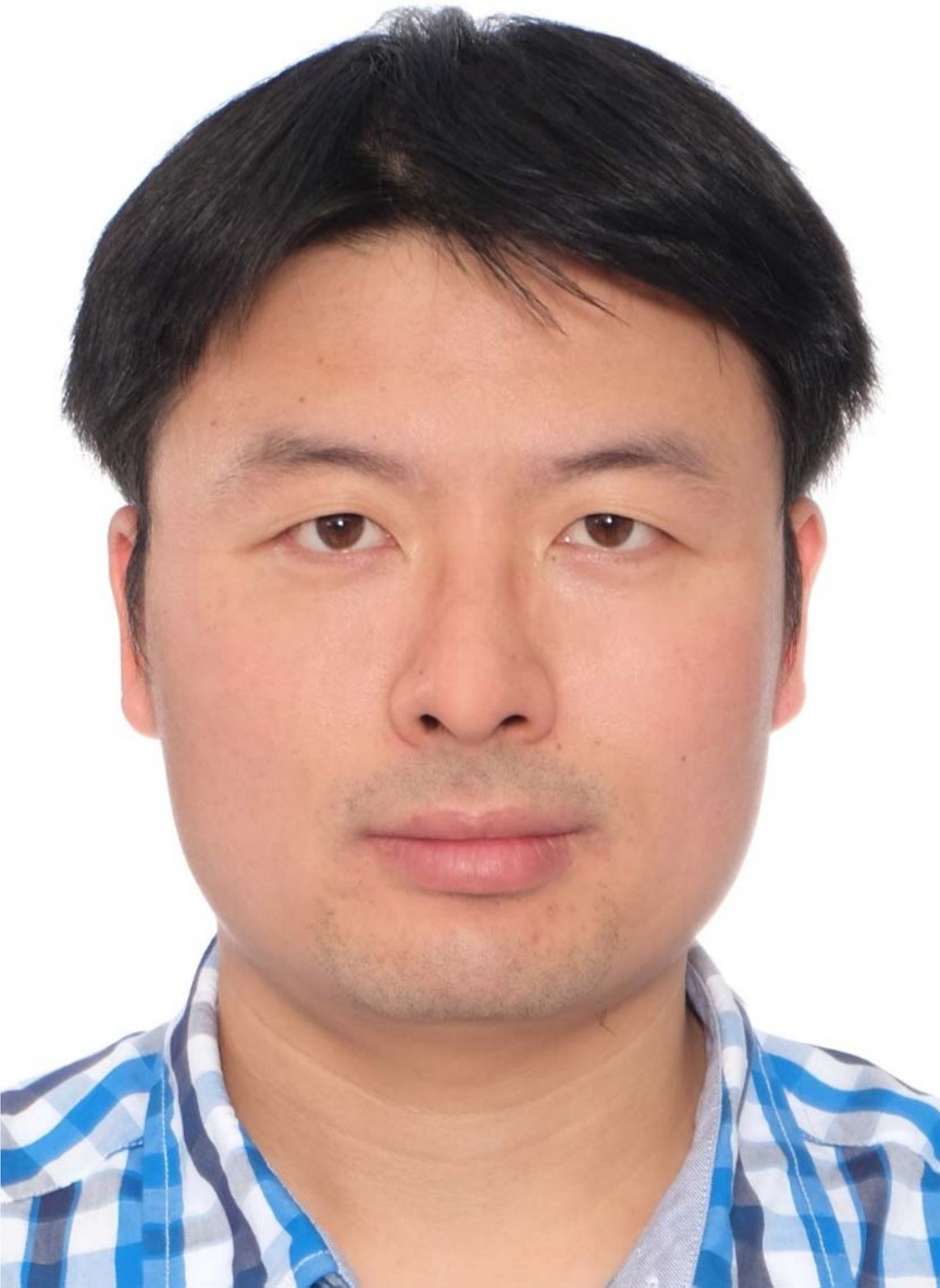}}]{Zengfu Wang } 
received the B.Sc. degree in applied mathematics, the M.Sc. degree in control theory and control engineering, the Ph.D. degree in control science and engineering from Northwestern Polytechnical University, Xi'an, in 2005, 2008, and 2013 respectively. From 2014 to 2017, he was a Lecturer with Northwestern Polytechnical University, where he is currently an Associate Professor. From November 2014 to November 2015, he was a Postdoctoral Research Fellow with the Department of Electronic Engineering, City University of Hong Kong, Hong Kong, China. From December 2019 to December 2020, he was a Visiting Researcher with the Faculty of Electrical Engineering, Mathematics and Computer Science, Delft University of Technology, Delft, The Netherlands. His research interests include path planning, discrete optimization, and information fusion.
\end{IEEEbiography}

\begin{IEEEbiography}[{\includegraphics[width=1in,height=1.25in,clip,keepaspectratio]{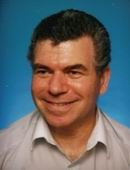}}]{Moshe Zukerman} (M’87–SM’91–F’07– LF’20)
received a B.Sc. degree in industrial engineering and management, an M.Sc. degree in operations research from the Technion – Israel Institute of Technology, Haifa, Israel, and a Ph.D. degree in engineering from the University of California, Los Angeles, in 1985. He was an independent consultant with the IRI Corporation and a Postdoctoral Fellow with the University of California, Los Angeles, in 1985–1986. In 1986–1997, he was with Telstra Research Laboratories (TRL), first as a Research Engineer and, in 1988–1997, as a Project Leader. He also taught and supervised graduate students at Monash University in 1990–2001. During 1997-2008, he was with The University of Melbourne, Victoria, Australia. In 2008 he joined the City University of Hong Kong (CityU) as a Chair Professor of Information Engineering, and a team leader. 
Between December 2020 and September 2022 he served as the Acting Chief Technology Officer of CityU.
He has over 300 publications in scientific journals and conference proceedings. He has served on various editorial boards such as Computer Networks, IEEE Communications Magazine, IEEE Journal of Selected Areas in Communications, IEEE/ACM Transactions on Networking, and Computer Communications.  
\end{IEEEbiography}

\end{document}